\documentclass[11pt,reqno]{amsart}

\usepackage[usenames,dvipsnames]{color}
\usepackage{hyperref}
\hypersetup{
    colorlinks=true,       			
    linkcolor=blue,          			
    citecolor=Magenta,        		
    filecolor=magenta,      			
    urlcolor=cyan           			
}

\usepackage{amsmath,amssymb,hyperref}
\usepackage[all]{xy}

\newtheorem*{thm}{Theorem}
\newtheorem*{prop}{Proposition}
\newtheorem*{lem}{Lemma}
\newtheorem*{cor}{Corollary}

\newenvironment{pf}{\paragraph{{\sc Proof}}}{\qed\par\medskip}
\theoremstyle{definition}

\newtheorem*{rem}{Remark}

\numberwithin{equation}{section}

\newcommand{\lp}{\left(}
\newcommand{\rp}{\right)}

\newcommand{\ga}{\mathfrak{a}}

\renewcommand{\b}{\mathfrak{b}}
\newcommand{\g}{\mathfrak{g}}
\newcommand{\h}{\mathfrak{h}}
\newcommand{\gl}{\mathfrak{gl}}
\newcommand{\Lg}{\mathfrak{g}}

\newcommand{\Lsl}{\mathfrak{sl}}
\newcommand{\Lgl}{\mathfrak{gl}}

\newcommand{\Sym}{\mathfrak{S}}

\newcommand{\U}{\mathcal{U}}

\newcommand{\D}{\mathcal{D}}
\newcommand{\E}{\mathcal{E}}
\newcommand{\F}{\mathcal{F}}
\newcommand{\K}{\mathcal{K}}
\newcommand{\G}{\mathcal{G}}
\newcommand{\V}{\mathcal{V}}
\newcommand{\HH}{\mathcal{H}}
\newcommand{\J}{\mathcal{J}}

\newcommand{\calN}{\mathcal{N}}
\newcommand{\PP}{\mathcal{P}}

\newcommand{\Cs}{\mathcal{S}}

\newcommand{\calZ}{\mathcal Z}

\newcommand{\C}{\mathbb{C}}
\newcommand{\nC}{\mathbb{C}^{\times}}
\newcommand{\N}{\mathbb{N}}
\newcommand{\Q}{\mathbb{Q}}
\newcommand{\R}{\mathbb{R}}
\newcommand{\Z}{\mathbb{Z}}

\newcommand{\fg}{{\mathbf g}}
\newcommand{\fh}{{\mathbf h}}
\newcommand{\fr}{{\mathbf r}}

\newcommand{\Th}{\Theta}

\newcommand {\wh}[1]{\widehat{#1}}
\newcommand {\ol}[1]{\overline{#1}}
\newcommand {\ul}[1]{\underline{#1}}
\newcommand {\eexp}{\operatorname{ch}}
\newcommand {\ddeg}{\operatorname{deg}}
\newcommand {\Ker}{\operatorname{Ker}}
\newcommand{\End}{\operatorname{End}}
\newcommand {\aand}{\qquad\text{and}\qquad}

\newcommand {\fd}{finite--dimensional }
\newcommand {\lhs}{left--hand side }
\newcommand {\rhs}{right--hand side }
\newcommand {\wrt}{with respect to }

\newcommand{\ds}{\displaystyle}
\newcommand{\wt}[1]{\widetilde{#1}}

\newcommand{\zed}{\mathfrak{z}}

\newcommand {\Omit}[1]{}

\newcommand {\ad}{\operatorname{ad}}
\newcommand {\Ad}{\operatorname{Ad}}
\newcommand {\gr}{\operatorname{gr}}
\newcommand {\cloop}{U(L\g)}

\newcommand {\hloop}{U_\hbar(L\Lg)}
\newcommand {\hloopgl}[1]{U_\hbar(L\gl_{#1})}
\newcommand {\Uhg}{U_\hbar\g}
\newcommand {\Yhg}{Y_\hbar(\g)}
\newcommand{\cYhg}{\widehat{Y_{\hbar}\g}}
\newcommand {\Yhgl}[1]{Y_\hbar(\gl_{#1})}
\newcommand {\current}{U(\g[s])}

\newcommand {\isom}{\stackrel{\sim}{\rightarrow}}

\newcommand {\cloopvar}{U(\g[z,z^{-1}])}
\newcommand {\currentvar}{U(\g[s])}
\newcommand{\cbin}[3]{\left(\begin{array}{c} #1\\ #2\end{array}\right)_{#3}}
\newcommand {\bin}[3]{\left[\begin{array}{c}#1\\#2\end{array}\right]_{#3}}

\newcommand {\Fd}{Finite--dimensional }
\newcommand {\ffd}{_{\scriptscriptstyle{\operatorname{fd}}}}
\newcommand {\bfw}{{\mathbf w}}
\newcommand {\Rep}{\operatorname{Rep}}
\newcommand {\eg}{{\it e.g.,} }
\newcommand {\ie}{{\it i.e.,} }

\newcommand {\bfA}{{\mathbf A}}
\newcommand {\bfI}{{\mathbf I}}
\newcommand {\bfJ}{{\mathbf J}}
\newcommand {\Hom}{\operatorname{Hom}}
\newcommand {\fing}{finitely--generated }
\renewcommand {\H}{\mathcal H}

\newcommand {\fact}[2]{#1} 

\newcommand {\dd}{\mathbf{d}}
\newcommand {\Sdd}{{\Sym(\dd)}}
\newcommand {\Sddp}{{\Sym(\dd')}}
\newcommand {\Td}{\operatorname{Td}}
\newcommand {\td}{\operatorname{td}}
\newcommand {\Id}{\operatorname{id}}

\newcommand {\hloopb}[1]{U_\hbar(L\b_{#1})}
\newcommand {\hlooph}{U_\hbar(L\h)}
\newcommand {\hloopi}[1]{U_\hbar(L\Lsl_2^{#1})}

\newcommand {\Yhb}[1]{Y_\hbar(\b_{#1})}
\newcommand {\Yhh}{Y_\hbar(\h)}
\newcommand {\Yhgi}[1]{Y_\hbar(\Lsl_2^{#1})}

\newcommand {\tor}[1]{U_\hbar^{\operatorname{\scriptscriptstyle{tor}}}(#1)}

\renewcommand {\>}{\rangle}

\begin{document}

\title{Yangians and quantum loop algebras}
\author[S. Gautam]{Sachin Gautam}
\address{Mathematics Department,
Columbia University,
2990 Broadway,
New York, NY 10027}
\address{Department of Mathematics,
Northeastern University,
360 Huntington Avenue,
Boston, MA 02115.}
\email{sachin@math.columbia.edu}
\author[V. Toledano Laredo]{Valerio Toledano Laredo}
\address{Department of Mathematics,
Northeastern University,
360 Huntington Avenue,
Boston MA 02115.}
\email{V.ToledanoLaredo@neu.edu}
\thanks{Both authors were supported by NSF grants DMS--0707212 and DMS--0854792}
\begin{abstract}
Let $\g$ be a complex, semisimple Lie algebra. Drinfeld showed that
the quantum loop algebra $\hloop$ of $\g$ degenerates to the Yangian
$\Yhg$. We strengthen this result by constructing an explicit algebra
homomorphism $\Phi$ from $\hloop$ to the completion of $\Yhg$ with
respect to its grading. We show moreover that $\Phi$ becomes an
isomorphism when $\hloop$ is completed \wrt its evaluation ideal.
We construct a similar homomorphism for $\Lg=\Lgl_n$, and show
that it intertwines the actions of $\hloopgl{n}$ and $\Yhgl{n}$ on the
equivariant $K$--theory and cohomology of the variety of $n$--step
flags in $\C^d$ constructed by Ginzburg--Vasserot.
\end{abstract}
\maketitle

\setcounter{tocdepth}{1}
\tableofcontents

\section{Introduction}

\subsection{} 

The present paper is motivated by, and lays the groundwork for a proof of the
trigonometric monodromy conjecture formulated by the second author in \cite
{valerio3}. Let $\g$ be a complex, semisimple Lie algebra, $G$ the corresponding
connected and simply--connected Lie group, $H\subset G$ a maximal torus
and $W$ the corresponding Weyl group. In \cite{valerio3} a flat, $W$--equivariant
connection $\wh{\nabla}_C$ was constructed on $H$ which has logarithmic
singularities on the root subtori of $H$ and values in any \fd representation
of the Yangian $\Yhg$. By analogy with the description of the monodromy
of the rational Casimir connection of $\Lg$ obtained in \cite{valerio4,valerio0},
it was conjectured in \cite{valerio3} that the monodromy of the trigonometric
Casimir connection $\wh{\nabla}_C$ is described by the action of the affine
braid group of $\g$ arising from the quantum Weyl group operators of the
quantum loop algebra $\hloop$. This raises in particular the problem of
relating \fd representations of $\Yhg$ and $\hloop$.

\subsection{}\label{ss:belief} 

Since their construction by Drinfeld \cite{drinfeld-qybe,drinfeld-quantum-groups},
these affine quantum groups have been extensively studied from several perspectives
(see, \eg \cite[chap. 12] {chari-pressley}, \cite{chari-hernandez} and references
therein), and are widely believed to share the same \fd representation theory.
This belief is corroborated in part by the following facts
\begin{enumerate}
\item\label{it:deg} The quantum loop algebra $\hloop$ degenerates to the Yangian
$\Yhg$. Specifically, if $\hloop$ is filtered by the powers of the evaluation ideal at $
z=1$, its associated graded is isomorphic to $\Yhg$ \cite{drinfeld-quantum-groups,
guay-degeneration}.
\item\label{it:polys} \Fd simple modules over $\hloop$ are parametrised by
$\bfI$--tuples of (Drinfeld) polynomials $\{P_i(u)\}_{i\in\bfI}$ satisfying $P_i
(0)=1$, where $\bfI$ is the set of vertices of the Dynkin diagram of $\g$ \cite
{chari-pressley-qaffine}. Similarly, \fd simple modules over $\Yhg$ are classified
by $\bfI$--tuples of monic polynomials \cite{drinfeld-yangian-qaffine,tarasov-1984,
tarasov-1985,chari-pressley-yangian}.
\item\label{it:grep} If $\g$ is simply laced, there exists, for every $\bfw\in\N^\bfI$,
a Steinberg variety $Z(\bfw)$ endowed with an action of $GL(\bfw)\times\nC$
(here $GL(\bfw) = \prod_{i\in\bfI} GL_{w_i}$), and algebra homomorphisms
\begin{align*}
\Psi_U : \hloop 	&\rightarrow K^{GL(\bfw)\times \nC}(Z(\bfw))\\
\Psi_Y : \Yhg 	&\rightarrow H^{GL(\bfw)\times \nC}(Z(\bfw))
\end{align*}
The variety $Z(\bfw)$ and the homomorphism $\Psi_U$ were constructed by 
Nakajima \cite{nakajima-qaffine}, while $\Psi_Y$ was constructed by Varagnolo
\cite{varagnolo-yangian}.
\end{enumerate}

\subsection{}

The above results go some way towards relating the categories of \fd representations
of $\hloop$ and $\Yhg$. For example, exponentiating the roots of Drinfeld polynomials
yields, via \eqref{it:polys}, a surjective map $\exp^*$ between the set of isomorphism
classes of irreducible \fd modules of $\Yhg$ and those of $\hloop$--modules. If $\g$
is simply laced, the geometric realisations \eqref{it:grep} imply further that $\exp^*$
preserves the dimensions of these representations \cite{varagnolo-yangian}.

Despite these results however, and to the best of our knowledge, no natural
relation between the categories of \fd representations of $\hloop$ and $\Yhg$ is
known. Part of the difficulty in exploiting, say, the geometric realisations to pursue
this question lies in the fact that the homomorphisms $\Psi_U,\Psi_Y$ are neither
injective nor surjective \cite{nakajima-qaffine}. Moreover, although these realisations
yield all irreducible representations, the categories $\Rep\ffd(\hloop)$ and $\Rep
\ffd(\Yhg)$ are not semisimple.

\subsection{} 

The aim of the present paper is to clarify the relation between $\hloop$
and $\Yhg$. We do so by constructing an explicit algebra homomorphism
\[\Phi:\hloop\longrightarrow\wh{\Yhg}\]
where $\wh{\Yhg}$ is the completion of $\Yhg$ \wrt its $\N$--grading,
and show that it induces an isomorphism of completed algebras. 
We also show that $\Phi$ exponentiates the roots of Drinfeld polynomials,
though we defer the study of the corresponding pull--back functor
\[F=\Phi^*:\Rep\ffd(\Yhg)\to\Rep\ffd(\hloop)\]
to the sequel of this paper \cite{sachin-valerio-2}.


To state out results more precisely, recall that $\hloop$ and $\Yhg$ are
deformations of the loop and current algebras  $U(\g[z,z^{-1}])$ and $U
(\g[s])$ respectively. Denote by
$$\hlooph,\hloopb{\pm}\subset\hloop\aand \Yhh,\Yhb{\pm}\subset\Yhg$$
the subalgebras deforming $U(\h[z,z^{-1}])$, $U(\b_\pm[z,z^{-1}])$ and
$U(\h[s])$, $U(\b_\pm[s])$ respectively, where $\h\subset\g$ is the Lie
algebra of $H$ and $\b_\pm\subset\g$ are the opposite Borel subalgebras
corresponding to a choice $\{\alpha_i\}_{i\in\bfI}$ of simple roots of $\g$.
For any $i\in\bfI$, let $\Lsl_2^i\subset\g$ be the corresponding
3--dimensional subalgebra and denote by
$$\hloopi{i}\subset\hloop\aand
\Yhgi{i}\subset\Yhg$$
the subalgebras which deform $U(\Lsl_2^i[z,z^{-1}])$ and $U(\Lsl_2^i[s])$
respectively. Then, the main result of this paper is the following
\begin{thm}\label{th:main intro}
There exists an explicit algebra homomorphism $\Phi:\hloop\to\wh{\Yhg}$
with the following properties
\begin{enumerate}
\item $\Phi$ is defined over $\Q[[\hbar]]$.
\item $\Phi$ induces an isomorphism $\wh{\hloop}\to\wh{\Yhg}$, where $\wh
{\hloop}$ is the completion of $\hloop$ \wrt the ideal of $z=1$.
\item $\Phi$ induces Drinfeld's degeneration of $\hloop$ to $\Yhg$.
\item $\Phi$ restricts to a homomorphism $\hlooph\to\wh{\Yhh}$ which induces
the exponentiation of roots on Drinfeld polynomials.
\item $\Phi$ restricts to a homomorphism $\hloopb{\pm}\to\wh{\Yhb{\pm}}$.
\item $\Phi$ restricts to a homomorphism $\hloopi{i}\to\wh{\Yhgi{i}}$ for any $i\in\bfI$.
\end{enumerate}
\end{thm}

\noindent
It is interesting to note that Theorem \ref{th:main intro} stands in stark contrast
with the analogous \fd situation. Indeed, if $\g\ncong\Lsl_2$, no explicit isomorphisms
are known between the quantum group $\Uhg$ and the undeformed enveloping
algebra $U\g[[\hbar]]$ \cite[\S 6.4]{chari-pressley}. Moreover, if $\g\ncong\Lsl_2$,
no algebra isomorphism $\Uhg\to U\g[[\hbar]]$ maps $U_\hbar\Lsl_2^i$ to $U\Lsl
_2^i[[\hbar]]$ for every $i\in\bfI$ \cite[Prop. 3.2]{valerio0}.

\subsection{}\label{ss:the Phi}

The homomorphism $\Phi$ has the following form.
Let $\{E_{i,k},F_{i,k},H_{i,k}\}_{i\in\bfI,k\in\Z}$ be the loop generators of $\hloop$
and $\{x^\pm_{i,m},\xi_{i,m}\}_{i\in\bfI,m\in\N}$ those of $\Yhg$ (see \cite
{drinfeld-yangian-qaffine} and Section \ref{section-definitions} for definitions).
Then,
\begin{align*}
\Phi(H_{i,0})&=d_i^{-1}t_{i,0}\\
\Phi(H_{i,r})&=\frac{\hbar}{q_i-q_i^{-1}}\sum_{m\geq 0}t_{i,m}\frac{r^m}{m!}\\
\Phi(E_{i,k})&=e^{k\sigma_i^+}\sum_{m\geq 0} \fact{g^+_{i,m}}{m!}\,x^+_{i,m}\\
\Phi(F_{i,k})&=e^{k\sigma_i^-}\sum_{m\geq 0} \fact{g^-_{i,m}}{m!}\,x^-_{i,m}
\end{align*}

\noindent
In the formulae above, $r\in\Z^*$, $k\in\Z$, $q=e^{\hbar/2}$ and $q_i=q^{d_i}$,
where the $d_i$ are the symmetrising integers for the Cartan matrix of $\g$. The
$\{t_{i,m}\}_{i\in\bfI,m\in\N}$ are an alternative set of generators of the commutative subalgebra
$\Yhh\subset\Yhg$ generated by the elements $\{\xi_{i,m}\}_{i\in\bfI,m\in\N}$.
They are defined in \cite{levendorskii} by equating the generating functions
$$\hbar\sum_{m\geq 0}t_{i,m}u^{-m-1}=\log(1+\hbar\sum_{m\geq 0}\xi_{i,m}u^{-m-1})$$
The elements $\{g^\pm_{i,m}\}_{i\in\bfI,m\in\N}$ lie in the completion of
$\Yhh$, and are constructed as follows. Consider the formal power series
\[G(v)=\log\left(\frac{v}{e^{v/2} - e^{-v/2}}\right)\in v\Q[[v]]\]
and define $\gamma_i(v)\in\wh{Y^0}[[v]]$ by
\[\gamma_i(v) = 
\hbar\sum_{r\geq 0}\frac{t_{i,r}}{r!}\left(-\frac{d}{dv}\right)^{r+1}G(v)\]
Then,
\begin{equation}\label{eq:the solution}
\sum_{m\geq 0}\fact{g_{i,m}^\pm}{m!}v^m=
\left(\frac{\hbar}{q_i-q_i^{-1}}\right)^{1/2}
\exp\left(\frac{\gamma_i(v)}{2}\right)
\end{equation}

\noindent
Finally, $\sigma_i^\pm$ are the homomorphisms of the subalgebras $\Yhb
{\pm}\subset\Yhg$ generated by $\{\xi_{j,r},x_{j,r}^\pm\}_{j\in\bfI,r\in\N}$,
which fix the $\xi_{j,r}$ and act on the remaining generators as the shifts
$x_{j,r}^\pm\to x^\pm_{j,r+\delta_{ij}}$.

Note that the formulae connecting the generators $\{H_{i,k}\}$ of $\hloop$
and $\{t_{i,m}\}$ of $\Yhh$ essentially coincide with those connecting the
generators of $\h[z,z^{-1}]$ and $\h[s]$.

\subsection{} 

The above formulae apply equally well when $\g$ is a symmetrisable
Kac--Moody algebra. Our proof of Theorem \ref{th:main intro} shows
that they define a homomorphism from the quantum affinization $\wh
{\Uhg}$ of the quantum group $\Uhg$ \cite{jing,nakajima-qaffine}, to
the completion of the Yangian $\Yhg$, provided the following holds
\begin{enumerate}
\item the entries of the Cartan matrix of $\g$ satisfy $a_{ij}a_{ji}
\leq 3$ for $i\neq j\in\bfI$.
\item the PBW theorem holds for $\Yhg$.
\end{enumerate}

The first assumption is equivalent to requiring that all rank 2
subalgebras of $\g$ be finite--dimensional, and is needed in
our proof of the $q$--Serre relations.
The second is required for the construction of certain straightening
homomorphisms on $\Yhg$ which are needed in the proof of Theorem
\ref{th:main intro}. We note that, for the Yangians associated to affine
Kac--Moody algebras, the PBW theorem was proved by Guay in
type ${\mathsf A}_n$, for $n\geq 4$ \cite{guay-affine yangian} and,
more recently, by Guay--Nakajima for all simply laced cases \cite
{guay-private}. In particular, the above formulae define a homomorphism from the quantum
toro\"{\i}dal algebra $\tor{\g}$ associated to a simply laced,
simple Lie algebra $\g\ncong\Lsl_2$, to the completion of the affine
Yangian $Y_\hbar\wh{\g}$.

\subsection{} 

We also construct in this paper a homomorphism similar to the one
described in \ref{ss:the Phi} for $\Lg=\Lgl_n$, by relying on the geometric
realisation of $\hloopgl{n}$ obtained by Ginzburg and Vasserot \cite
{ginzburg-vasserot,vasserot-qaffine}. More precisely, fix integers $1
\leq n\leq d$, and let
\[\F = \{0= V_0 \subseteq V_1 \subseteq \cdots \subseteq V_n = \C^d\}\]
be the variety of $n$--step flags in $\C^d$. The cotangent bundle $T^*\F$
may be realised as
\[T^*\F = \{(V_\bullet,x)\in \F\times\End(\C^d)|\,x(V_i) \subset V_{i-1}\}\]
and therefore admits a morphism $T^*\F \rightarrow\calN$ via the second
projection, where $\calN=\{x\in\End(\C^d)|\,x^n=0\}$ is the cone of $n$--step
nilpotent endomorphisms. Define the Steinberg variety $Z = T^*\F \times_\calN
T^*\F$. The group $GL_d\times\nC$ acts on $T^*\F$ and $Z$ and there are
surjective algebra homomorphisms
\begin{gather*}
\Psi_U:\hloopgl{n} \rightarrow K^{GL_d\times \nC}(Z)\\
\Psi_Y:\Yhgl{n} \rightarrow H^{GL_d\times \nC}(Z)
\end{gather*}
see \cite{ginzburg-vasserot,vasserot-qaffine} for the definition of $\Psi_U$.

To understand these more explicitly, one can use the convolution actions
of $K^{GL_d\times\nC}(Z)$ on $K^{GL_d\times \nC}(T^*\F)$ and of $H^
{GL_d\times\nC}(Z)$ on $H_{GL_d\times\nC}(T^*\F)$, which are faithful.
By using the equivariant Chern character, we construct an algebra homomorphism
$$\Phi:\hloopgl{n}\to\wh{\Yhgl{n}}$$ which intertwines these two actions.

\subsection{}

In the sequel to this paper \cite{sachin-valerio-2}, we shall prove that,
for $\g$ semisimple, a modification of the pull--back functor $\Phi^*$
converges for numerical values of $\hbar$, and defines a functor
\[\Rep\ffd(Y_a(\g))\rightarrow\Rep\ffd(U_\epsilon(L\g))\]
where $Y_a(\g)$ is the specialisation of $\Yhg$ at $\hbar=a\in\C
\setminus \R$ and $\epsilon=\exp(\pi i a)$, and defines an equivalence
of an explicit subcategory of $\Rep\ffd(Y_a\g)$ with $\Rep\ffd(U_\epsilon
(L\g))$.

\subsection{} 

It is worth pointing out that most of our results relating $\hloop$ and $\Yhg$
have analogues for the affine and degenerate affine Hecke algebras $\HH$
and $\HH'$ associated to an affine Weyl group $W$ \cite{lusztig-graded-affine},
and were in fact inspired by these and their further study in \cite{cherednik-Lusztig}.
Indeed, in \cite{lusztig-graded-affine}, Lusztig constructs an explicit isomorphism
between appropriate completions of $\HH$ and $\HH'$. In this context, the
isomorphism can be understood in terms of, and in fact obtained from, the
geometric realisations
\begin{align*}
\Xi&:\HH  \stackrel{\sim}{\longrightarrow} K^{G\times \nC}(Z)\\
\Xi'&:\HH' \stackrel{\sim}{\longrightarrow} H^{G\times \nC}(Z)
\end{align*}
where $Z$ is the Steinberg variety corresponding to $W$ \cite{chriss-ginzburg,
ginzburg-Hecke}.

\subsection{Outline of the paper}

In Section \ref{section-definitions}, we review the definition of the quantum loop
algebra $\hloop$ and Yangian $\Yhg$ of a semisimple Lie algebra $\g$. We also
introduce shift homomorphisms of the subalgebras $\Yhb{\pm}$, and straightening
homomorphisms of the subalgebra $\Yhh$.

In Section \ref{section-nec-suff}, we consider assignments mapping the generators of $\hloop$
to $\wh{\Yhg}$. These have the form described in \ref{ss:the Phi}, where the elements $g_
{i,m}^\pm\in\wh{\Yhh}$ are however not necessarily given by formula \eqref{eq:the solution}.
Our main result, Theorem \ref{first-main-theorem}, gives necessary and sufficient conditions for
these elements to give rise to an algebra homomorphism. We call such homomorphisms {\it
of geometric type} since, for $\g$ simply laced, they are related to the Chern character in the
geometric realisation described in \ref{ss:belief}.

The proof that the elements given by \eqref{eq:the solution} satisfy the conditions of
Theorem \ref{first-main-theorem}, and therefore give rise to an algebra homomorphism
$\Phi$, is given in Section \ref{section-existence} (Theorem \ref{third-main-theorem}).
We also prove that the action of $\Phi$ on Drinfeld polynomials exponentiates their
roots (Corollary \ref{compatible-drinfeld-homomorphism}).

In Section \ref{section-uniqueness}, we prove the essential uniqueness of homomorphisms
of geometric type by showing that any two differ by conjugation by an element of the torus
$H$ and an invertible element of $\wh{\Yhh}$ (Theorem \ref{second-main-theorem}).

In Section \ref{section-degeneration}, we show that any homomorphism of geometric
type $\Phi$ induces an isomorphism $\wh{\hloop}\to\wh{\Yhg}$, where $\wh{\hloop}$
is the completion \wrt the evaluation ideal at $z=1$ (Theorem \ref{th:geo iso}). We show
moreover that the associated graded map coincides with Drinfeld's degeneration of
$\hloop$ to $\Yhg$ (Proposition \ref{pr:degeneration}).

Section \ref{section-gln} contains similar results for $\Lg = \Lgl_n$. In addition to
constructing an explicit homomorphism $\Phi:\hloopgl{n}\to\wh{\Yhgl{n}}$ (Theorem
\ref{first-main-theorem-gln}), we review the geometric realisations of these algebras
and show that $\Phi$ intertwines them (Theorem \ref{main-2-gln}).

Appendix \ref{app: Serre} contains a proof of the Serre relations which is required
to complete the proof of Theorem \ref{first-main-theorem}.

\subsection{Acknowledgments}

We are very grateful to Ian Grojnowski from whom we learned that the quantum loop
algebra and Yangian should be isomorphic after appropriate completions. His explanations
and friendly insistence helped us overcome our initial doubts. We are also grateful to
V. Drinfeld for showing us a proof that \fd representations separate elements of the
Yangian and allowing us to reproduce it in Appendix \ref{app: Serre}. We would also
like to thank N. Guay for sharing a preliminary version of the preprint \cite{guay-degeneration}
and E. Vasserot for useful discussions. 

\section{Quantum loop algebras and Yangians}\label{section-definitions}

\subsection{}\label{ss: KMA} 

Let $\g$ be a complex, semisimple Lie algebra and $(\cdot,\cdot)$ a non--degenerate,
invariant bilinear form on $\g$. Let $\h\subset\g$ be a Cartan subalgebra of
$\g$, $\{\alpha_i\}_{i\in\bfI}\subset\h^*$ a basis of simple roots of $\g$ relative
to $\h$ and $a_{ij}=2(\alpha_i,\alpha_j)/(\alpha_i,\alpha_i)$ the entries of the
corresponding Cartan matrix $\bfA$. Set $d_i=(\alpha_i,\alpha_i)/2$, so that
$d_ia_{ij}=d_j a_{ji}$ for any $i,j\in\bfI$. Let $\nu:\h\to\h^*$ be the isomorphism
determined by the inner product $(\cdot,\cdot)$ and set $h_i=\nu^{-1}(\alpha_i)
/d_i$. Choose root vectors $e_i\in\g_{\alpha_i},f_i\in\g_{-\alpha_i}$ such that
$[e_i,f_i]=h_i$. Recall that $\g$ is presented on generators $\{e_i,f_i,h_i\}$
subject to the relations
\begin{gather*}
[h_i,h_j]=0
\\[1.1 ex]
[h_i,e_j]= a_{ij}e_j
\qquad
[h_i,f_j]=-a_{ij}f_j
\\[1.1 ex]
[e_i,f_j]=
\delta_{ij}h_i
\\[1.1ex]
\intertext{for any $i,j\in\bfI$ and, for any $i\neq j\in\bfI$}
\ad(e_i)^{1-a_{ij}}e_j=0\\[1.1 ex]
\ad(f_i)^{1-a_{ij}}f_j=0
\end{gather*}

A closely related, but slightly less standard presentation may be obtained by
setting $t_i=\nu^{-1}(\alpha_i)=d_ih_i$ and choosing, for any $i\in\bfI$, root
vectors $x^\pm_i\in\g_{\pm \alpha_i}$ such that $[x_i^+,x_i^-]=t_i$. Then
$\g$ is presented on $\{x_i^\pm,t_i\}_{i\in\bfI}$ subject to the relations
\begin{gather*}
[t_i,t_j]=0\\
[t_i,x_j^\pm]=\pm d_ia_{ij}x_j^\pm\\
[x_i^+,x_j^-]=\delta_{ij}t_i\\
\ad(x_i^\pm)^{1-a_{ij}}x_j^\pm=0
\end{gather*}

\subsection{}

Throughout this paper, $q$ and $\hbar$ are formal variables related by $q^2
= e^{\hbar}$. For any $i\in\bfI$, we set $q_i=q^{d_i}=e^{\hbar d_i/2}$. We use
the standard notation for Gaussian integers
\begin{gather*}
[n]_q = \frac{q^n - q^{-n}}{q-q^{-1}}\\
[n]_q! = [n]_q[n-1]_q\cdots [1]_q\qquad\qquad
\bin{n}{k}{q} = \frac{[n]_q!}{[k]_q![n-k]_q!}
\end{gather*}

\subsection{The quantum loop algebra \cite{drinfeld-yangian-qaffine}}\label{definition-qla}

Let $\hloop$ be the unital, associative algebra over $\C[[\hbar]]$ topologically
generated by elements $\{E_{i,k},F_{i,k},H_{i,k}\}_{i\in\bfI,k\in\Z}$ subject to the
following relations
\begin{itemize}
\item[(QL1)] For $i,j\in\bfI$ and $r,s\in\Z$
\[[H_{i,r},H_{j,s}]=0\]
\item[(QL2)] For any $i,j\in\bfI$ and $k\in \Z$,
\[[H_{i,0},E_{j,k}]=a_{ij}E_{j,k}
\qquad
[H_{i,0},F_{j,k}]=-a_{ij}F_{j,k}\]
\item[(QL3)] For any $i,j\in\bfI$ and $r\in\Z^\times$,
\[[H_{i,r}, E_{j,k}]= \frac{[ra_{ij}]_{q_i}}{r} E_{j,r+k}
\qquad
[H_{i,r}, F_{j,k}]= -\frac{[ra_{ij}]_{q_i}}{r}F_{j,r+k}\]
\item[(QL4)] For $i,j\in\bfI$ and $k,l\in \Z$
\begin{align*}
E_{i,k+1}E_{j,l} - q_i^{a_{ij}}E_{j,l}E_{i,k+1} &= q_i^{a_{ij}}E_{i,k}E_{j,l+1} - E_{j,l+1}E_{i,k}\\
F_{i,k+1}F_{j,l} - q_i^{-a_{ij}}F_{j,l}F_{i,k+1} &= q_i^{-a_{ij}}F_{i,k}F_{j,l+1} - F_{j,l+1}F_{i,k}
\end{align*}
\item[(QL5)] For $i,j\in\bfI$ and $k,l\in \Z$
\[[E_{i,k}, F_{j,l}] = \delta_{ij} \frac{\psi_{i,k+l} - \phi_{i,k+l}}{q_i - q_i^{-1}}\]
\item[(QL6)] Let $i\not= j\in\bfI$ and set $m = 1-a_{ij}$. For every $k_1,\ldots, k_m\in \Z$
and $l\in \Z$
\begin{gather*}
\sum_{\pi\in \Sym_m} \sum_{s=0}^m (-1)^s\bin{m}{s}{q_i}
E_{i,k_{\pi(1)}}\cdots E_{i,k_{\pi(s)}} E_{j,l}E_{i,k_{\pi(s+1)}}\cdots E_{i,k_{\pi(m)}} = 0\\
\sum_{\pi\in \Sym_m} \sum_{s=0}^m (-1)^s\bin{m}{s}{q_i}
F_{i,k_{\pi(1)}}\cdots F_{i,k_{\pi(s)}} F_{j,l}F_{i,k_{\pi(s+1)}}\cdots F_{i,k_{\pi(m)}} = 0
\end{gather*}
\end{itemize}
where the elements $\psi_{i,r},\phi_{i,r}$ are defined by 
\begin{align*}
\psi_i(z) 
&=
\sum_{r\geq 0} \psi_{i,r}z^{-r} = 
\exp\left(\phantom{-}\frac{\hbar d_i}{2}H_{i,0}\right)
\exp\left(\phantom{-}(q_i - q_i^{-1})\sum_{s\geq 1} H_{i,s}z^{-s}\right)\\
\phi_i(z)
&=
\sum_{r\geq 0} \phi_{i,-r}z^r = 
\exp\left(-\frac{\hbar d_i}{2}H_{i,0}\right)
\exp\left(-(q_i - q_i^{-1})\sum_{s\geq 1} H_{i,-s}z^s\right)
\end{align*}
with $\psi_{i,-k}=\phi_{i,k}=0$ for every $k\geq 1$.

We shall denote by $U^0\subset\hloop$ the commutative subalgebra generated
by the elements $\{H_{i,r}\}_{i\in\bfI,r\in\Z}$. 

\subsection{The Yangian \cite{drinfeld-yangian-qaffine}}\label{definition-yangian}

Let $Y_{\hbar}(\Lg)$ be the unital, associative $\C[\hbar]$--algebra generated
by elements $\{x^{\pm}_{i,r},\xi_{i,r}\}_{i\in\bfI,r\in\N}$, subject to the following
relations

\begin{enumerate}
\item[(Y1)] For any $i,j\in\bfI$ and $r,s\in\N$
\[[\xi_{i,r}, \xi_{j,s}] = 0\]
\item[(Y2)] For $i,j\in\bfI$ and $s\in \N$
\[[\xi_{i,0}, x_{j,s}^{\pm}] = \pm d_ia_{ij} x_{j,s}^{\pm}\]
\item[(Y3)] For $i,j\in\bfI$ and $r,s\in\N$
\[[\xi_{i,r+1}, x^{\pm}_{j,s}] - [\xi_{i,r},x^{\pm}_{j,s+1}] =
\pm\frac{d_ia_{ij}\hbar}{2}(\xi_{i,r}x^{\pm}_{j,s} + x^{\pm}_{j,s}\xi_{i,r})\]
\item[(Y4)] For $i,j\in\bfI$ and $r,s\in \N$
\[
[x^{\pm}_{i,r+1}, x^{\pm}_{j,s}] - [x^{\pm}_{i,r},x^{\pm}_{j,s+1}]=
\pm\frac{d_ia_{ij}\hbar}{2}(x^{\pm}_{i,r}x^{\pm}_{j,s} + x^{\pm}_{j,s}x^{\pm}_{i,r})
\]
\item[(Y5)] For $i,j\in\bfI$ and $r,s\in \N$
\[[x^+_{i,r}, x^-_{j,s}] = \delta_{ij} \xi_{i,r+s}\]
\item[(Y6)] Let $i\not= j\in\bfI$ and set $m = 1-a_{ij}$. For any $r_1,\cdots, r_m\in \N$ and $s\in \N$
\[\sum_{\pi\in\Sym_m}
\left[x^{\pm}_{i,r_{\pi(1)}},\left[x^{\pm}_{i,r_{\pi(2)}},\cdots, \left
[x^{\pm}_{i,r_{\pi(m)}},x^{\pm}_{j,s}\right]\cdots\right]\right]=0\]
\end{enumerate}

\noindent $\Yhg$ is an $\N$--graded algebra by $\deg(\xi_{i,r})=\deg(x_{i,r}^\pm)=r$
and $\deg(\hbar) = 1$.

\subsection{PBW theorem for $\Yhg$}\label{ss:PBW}

For any positive root $\beta$ of $\g$, choose a sequence of simple roots
$\alpha_{i_1},\ldots,\alpha_{i_k}$ such that $\beta=\alpha_{i_1}+\cdots+
\alpha_{i_k}$ and
$$[x^\pm_{i_1}, [x^\pm_{i_2}, \cdots, [x^\pm_{i_{k-1}}, x^\pm_{i_k}]\cdots]]
\in\g_{\pm\beta}$$
are non--zero vectors. For any $r\in\N$, define $x^\pm_{\beta,r}\in\Yhg$
by choosing a partition $r=r_1+\cdots+r_k$ of length $k$ and setting
\[x^\pm_{\beta,r}=
[x^\pm_{i_1,r_1},[x^\pm_{i_2,r_2},\cdots,[x^\pm_{i_{k-1}, r_{k-1}},x^\pm_{i_k,r_k}]\cdots ]]\]

\begin{thm}[\cite{levendorskii-PBW}]\label{PBW-yangian}
Fix a total order on the set $\G=\{\xi_{i,r},x^\pm_{\beta,r}\}_{i\in\bfI,r\in\N,\beta
\in\Sigma_+}$. Then, the ordered monomials in the elements of $\G$ form a
basis of $\Yhg$.
\end{thm}

Let $Y^0,Y^\pm\subset\Yhg$ be the subalgebras generated by the elements
$\{\xi_{i,r}\}_{i\in\bfI,r\in\N}$ (resp. $\{x_{i,r}^\pm\}_{i\in\bfI,r\in\N}$) and $Y^{\geq 0},
Y^{\leq 0}\subset\Yhg$ the subalgebras generated by $Y^0,Y^+$ and $Y^0,
Y^-$ respectively. The following is a direct consequence of Theorem \ref{PBW-yangian}.

\begin{cor}\label{co:PBW}\hfill
\begin{enumerate}
\item $Y^0$ is a polynomial algebra in the generators $\{\xi_{i,r}\}_{i\in\bfI,r\in\N}$.
\item $Y^\pm$ is the algebra generated by elements $\{x_{i,r}^\pm\}_{i\in\bfI,r\in\N}$
subject to the relations (Y4) and (Y6).
\item $Y^{\geq 0}$ (resp. $Y^{\leq 0}$) is the algebra generated by elements
$\{\xi_{i,r},x_{i,r}^\pm\}_{i\in\bfI,r\in\N}$ subject to the relations (Y1)--(Y4) and (Y6).
\item Multiplication induces an isomorphism of vector spaces
\[Y^-\otimes Y^0\otimes Y^+ \rightarrow \Yhg\]
\end{enumerate}
\end{cor}

\subsection{The shift operators $\sigma_i^\pm$}\label{ss:sigma_i}

Fix $i\in\bfI$. By Corollary \ref{co:PBW} (3), the assignment
\[x_{j,r}^\pm\to x_{j,r+\delta_{ij}}^{\pm}\qquad\qquad\xi_{j,r}\to\xi_{j,r}\]
extends to an algebra homomorphism $Y^{\geq 0}\to Y^{\geq 0}$
(resp. $Y^{\leq 0}\to Y^{\leq 0}$) which we shall denote by $\sigma
_i^\pm$.

\subsection{The relations (Y2)--(Y3)}

We rewrite below the defining relations (Y2)--(Y3) of $\Yhg$ in terms
of the shift operators $\sigma_j^\pm$ and the generating series
\begin{equation}\label{eq:xi(u)}
\xi_i(u)=1+\hbar\sum_{r\geq 0}\xi_{i,r}u^{-r-1}\in\Yhg[[u^{-1}]]
\end{equation}

\begin{lem}\label{le:Y23}
The relations (Y2)--(Y3) are equivalent to
\[[\xi_i(u),x_{j,s}^\pm]=
\frac{\pm\hbar d_ia_{ij}}{u-\sigma_j^\pm\pm\hbar d_ia_{ij}/2}\xi_i(u)x_{j,s}^\pm\]
where the rational function on the \rhs is expanded in powers of $u^{-1}$.
\end{lem}
\begin{pf}

Set $a=d_ia_{ij}/2$. Multiplying (Y3) by $\hbar u^{-r-1}$ and summing
over $r\geq 0$ yields
\[u[\xi_i(u)-1-\hbar u^{-1}\xi_{i,0},x_{j,s}^\pm]-[\xi_i(u)-1,x_{j,s+1}^\pm]=
\pm\hbar a\{x_{j,s}^\pm,\xi_i(u)-1\}\]
where $\{x,\xi\}=x\xi+\xi x$. Using (Y2) and $\{x,\xi\}=[x,\xi]+2\xi x$, yields
\begin{equation}\label{eq:num}
(u-\sigma_j^\pm\pm\hbar a)[\xi_i(u),x_{j,s}^\pm]=\pm 2\hbar a\xi_i(u)x_{j,s}^\pm
\end{equation}
as claimed. Conversely, taking the coefficients of $u^0$ and $u^{-r-1}$
in \eqref{eq:num} yields (Y2) and (Y3) respectively.
\end{pf}

\subsection{The relations (Y4) and (Y6)}\label{shift-symbols}

We shall use the following notation
\begin{itemize}
\item for an operator $T \in\End(V)$, $T_{(i)} \in\End(V^{\otimes m})$ is defined as
\[T_{(i)} = 1^{\otimes i-1}\otimes T\otimes 1^{\otimes m-i}\]
\item for an algebra $A$, $ad^{(m)} : A^{\otimes m} \rightarrow\End(A)$ is defined as
\[ad^{(m)}\left(a_1\otimes \cdots \otimes a_m\right) = ad(a_1)\circ \cdots \circ ad(a_m)\]
\end{itemize}

\begin{prop}\label{algebraic-properties-shift}\hfill
\begin{enumerate}
\item The relation (Y4) for $i\neq j$ is equivalent to the requirement that the following
holds for any $A(v_1,v_2)\in\C[[v_1,v_2]]$
\[
A(\sigma_i^{\pm}, \sigma_j^{\pm})
(\sigma_i^{\pm} - \sigma_j^{\pm} \mp a\hbar)x^{\pm}_{i,0}x^{\pm}_{j,0} =
A(\sigma_i^{\pm}, \sigma_j^{\pm})
(\sigma_i^{\pm} - \sigma_j^{\pm} \pm a\hbar)x^{\pm}_{j,0}x^{\pm}_{i,0}
\]
where $a=d_ia_{ij}/2$.
\item The relation (Y4) for $i=j$ is equivalent to the requirement that the following
holds for any $B(v_1,v_2)\in\C[[v_1,v_2]]$ such that $B(v_1,v_2) = B(v_2,v_1)$
\begin{equation}\label{eq:B(u,v)}
\mu\left(B(\sigma^{\pm}_{i, (1)}, \sigma^{\pm}_{i,(2)})
(\sigma^{\pm}_{i,(1)} - \sigma^{\pm}_{i,(2)} \mp d_i\hbar)x^{\pm}_{i,0}\otimes x^{\pm}_{i,0}\right) = 0
\end{equation}
where $\mu:\Yhg^{\otimes 2} \rightarrow\Yhg$ is the multiplication.
\item The relation (Y6) is equivalent to the requirement that the following
holds for any $i\neq j$ and $A\in\C[v_1,\ldots,v_m]^{\Sym_m}$ with
$m=1-a_{ij}$
\[ad^{(m)}\left(A(\sigma^{\pm}_{i,(1)},\ldots, \sigma^{\pm}_{i,(m)})
\left(x^{\pm}_{i,0}\right)^{\otimes m}\right) x^{\pm}_{j,l} = 0\]
\end{enumerate}
\end{prop}
\begin{pf}
(1) The relation (Y4)
\[[x^{\pm}_{i,r+1}, x^{\pm}_{j,s}] - [x^{\pm}_{i,r}, x^{\pm}_{j,s+1}] =
\pm a\hbar (x^{\pm}_{i,r}x^{\pm}_{j,s} + x^{\pm}_{j,s}x^{\pm}_{i,r})\]
may be rewritten as
\[\sigma_i^{\pm r}\sigma_j^{\pm s} \left(\sigma^{\pm}_i -\sigma^{\pm}_j \mp a\hbar\right) x^{\pm}_{i,0}x^{\pm}_{j,0}
=
\sigma_i^{\pm r}\sigma_j^{\pm s} \left(\sigma^{\pm}_i - \sigma^{\pm}_j \pm a\hbar\right) x^{\pm}_{j,0}x^{\pm}_{i,0}\]

(2) If $i=j$, $a=d_i$ and the above may be written as
\begin{multline*}
\mu\left(\sigma^{\pm r}_{i,(1)}\sigma_{i,(2)}^{\pm s} \left(\sigma^{\pm}_{i,(1)} -\sigma^{\pm}_{i,(2)} \mp d_i\hbar\right)
x^{\pm}_{i,0}\otimes x^{\pm}_{i,0}\right) = \\
\mu\left( \sigma_{i,(1)}^{\pm s} \sigma_{i,(2)}^{\pm r}\left(\sigma^{\pm}_{i,(2)} - \sigma^{\pm}_{i,(1)} \pm d_i\hbar\right)
x^{\pm}_{i,0}\otimes x^{\pm}_{i,0}\right)
\end{multline*}
which is equivalent to
\[\mu\left(\left(\sigma_{i,(1)}^{\pm r}\sigma_{i,(2)}^{\pm s} + \sigma_{i,(1)}^{\pm s}\sigma_{i,(2)}^{\pm r}\right)
\left(\sigma^{\pm}_{i,(1)} - \sigma^{\pm}_{i,(2)} \mp d_i\hbar\right)x^{\pm}_{i,0}\otimes x_{i,0}^{\pm}\right) = 0\]

(3) is just the reformulation of (Y6).
\end{pf}

\begin{cor}\label{co:B(u,v)}
If \eqref{eq:B(u,v)} holds for some $B\in\C[[v_1,v_2]]$, then $B(v_1,v_2)
=B(v_2,v_1)$.
\end{cor}
\begin{pf}
By (2) of Proposition \ref{algebraic-properties-shift}, we may assume that
$B(v_1,v_2)=-B(v_2,v_1)$ and therefore that $B=(v_1-v_2)\ol{B}$ where
$\ol{B}$ is symmetric in $v_1\leftrightarrow v_2$. Using the grading on $
\Yhg$, we may further assume that $\ol{B}$ is proportional to $v_1^rv_2
^s+v_1^sv_2^r$ for some $r\geq s\in\N$. An application of (Y4) then yields
\begin{multline*}
\mu\Bigl((\sigma^\pm_{i,(1)}-\sigma^\pm_{i,(2)})
(\sigma_{i,(1)}^{\pm r}\sigma_{i,(2)}^{\pm s} + \sigma_{i,(1)}^{\pm s}\sigma_{i,(2)}^{\pm r})
(\sigma^{\pm}_{i,(1)} - \sigma^{\pm}_{i,(2)} \mp d_i\hbar)x^{\pm}_{i,0}\otimes x^{\pm}_{i,0}\Bigr)\\
=2\Bigl((x_{i,r+2}^\pm x_{i,s}^\pm-x_{i,r+1}^\pm x_{i,s+1}^{\pm}\mp d_i\hbar x_{i,r+1}^\pm x_{i,s}^\pm)
-
(x_{i,r+1}^\pm x_{i,s+1}^\pm-x_{i,r}^\pm x_{i,s+2}^{\pm}\mp d_i\hbar x_{i,r}^\pm x_{i,s+1}^\pm)\Bigr)
\end{multline*}
If $r\geq s+2$, the above is not zero by the PBW Theorem \ref{PBW-yangian}
and $\ol{B}=0$. If $r=s+1$, a further application of (Y4) shows that the second
of the above two parenthesized summands is zero and again $\ol{B}=0$ by
Theorem \ref{PBW-yangian}. Finally, if $r=s$, (Y4) implies that the two
parenthesized summands are opposites of each other and again $\ol{B}=0$.
\end{pf}

\subsection{An alternative system of generators for $Y^0$}\label{new-generators}

The following generators of $Y^0$ were introduced in \cite{levendorskii}.
For any $i\in\bfI$, define the formal power series
\[t_i(u) = \hbar\sum_{r\geq 0} t_{i,r}u^{-r-1}\in Y^0[[u^{-1}]]\]
by
\begin{equation}\label{log-like-yangian}
t_i(u)=\log(\xi_i(u))=\log\left(1 + \hbar\sum_{r\geq 0} \xi_{i,r}u^{-r-1}\right)
\end{equation}
Since \eqref{log-like-yangian} can be inverted, $\{t_{i,r}\}_{i\in\bfI,r\in\N}$
is another system of generators of $Y^0$. These are homogeneous, with
$\deg(t_{i,r})=r$, since $\zeta\cdot\xi_i(u)=\xi_i(\zeta^{-1}u)$, where $\cdot$
is the action of $\zeta\in\C^*$ on $\Yhg$ given by the grading. Moreover,
$t_{i,0}=\xi_{i,0}$ and $t_{i,r}=\xi_{i,r}\mod\hbar$ for any $r\geq 1$ since
$$t_i(u) = \hbar\sum_{r\geq 0} \xi_{i,r}u^{-r-1}\mod\hbar^2$$

To compute the commutation relations between $t_{i,r}$ and $x^{\pm}_{j,s}$,
we introduce the following formal power series (inverse Borel transform of
$t_i(u)$)
\begin{equation}\label{inverse-borel}
B_i(v) = B(t_i(u)) =
\hbar\sum_{r\geq 0} t_{i,r}\frac{v^r}{r!}\in Y^0[[v]]
\end{equation}
\begin{lem}\label{commutation-with-log-like}
For any $i,j\in\bfI$ we have
\begin{equation}\label{eq:B x}
\left[B_i(v), x^{\pm}_{j,s}\right] =
\pm\frac{q_i^{a_{ij}v} - q_i^{-a_{ij}v}}{v} e^{\sigma^{\pm}_j v}\,x^{\pm}_{j,s}
\end{equation}
\end{lem}
\begin{pf}
To simplify notations set $a=\pm\hbar d_ia_{ij}/2$, so that $e^a=q_i^{\pm a_{ij}}$.
By Lemma \ref{le:Y23}
\[\xi_i(u)x^{\pm}_{j,s}\xi_i(u)^{-1}
=\frac{u-\sigma^{\pm}_j+a}{u-\sigma^{\pm}_j-a} x^{\pm}_{j,s}\]
so that
\[[t_i(u), x^{\pm}_{j,s}]=
\log\left(\frac{1-(\sigma^{\pm}_j-a)u^{-1}}{1-(\sigma^{\pm}_j+a)u^{-1}}\right) x^{\pm}_{j,s}\]
Using
\begin{equation}\label{eq:B log}
B\left(\log(1-pu^{-1})\right) = \frac{1-e^{pv}}{v}
\end{equation}
this yields
\[[B_i(v), x^{\pm}_{j,s}]
= 
\left(\frac{1-e^{(\sigma^{\pm}_j-a)v}}{v}-\frac{1-e^{(\sigma^{\pm}_j+a)v}}{v}\right)x^{\pm}_{j,s}\\
=
\frac{e^{av}-e^{-av}}{v} e^{\sigma^{\pm}_j v} x^{\pm}_{j,s}\]
as claimed.
\end{pf}

\begin{rem}
Expanding the \rhs of \eqref{eq:B x} as a power series in $v$ yields the
commutation relations
\[[t_{i,r}, x^{\pm}_{j,s}]=
\pm d_ia_{ij} \sum_{l=0}^{\lfloor r/2\rfloor} \cbin{r}{2l}{}
\frac{(\hbar d_ia_{ij}/2)^{2l}}{2l+1} x^{\pm}_{j,r+s-2l}\]
These relations were obtained in this form in \cite[Lemma 1.4]{levendorskii}.
\end{rem}

\subsection{The operators $\mathbf{\lambda_i^{\pm}(v)}$}\label{rho-lambda}

We introduce below operators which straighten monomials of the form
$x_{i,m}^\pm\xi$, $\xi\in Y^0$, into elements of $Y^0\cdot Y^\pm$.

\begin{prop}\label{raising-lowering-operators}
There are operators $\{\lambda^{\pm}_{i;s}\}_{i\in\bfI,s\in\N}$ on $Y^0$ such
that the following holds
\begin{enumerate}
\item For any $\xi\in Y^0$, the elements $\lambda^{\pm}_{i;s}(\xi)\in Y^0$
are uniquely determined by the requirement that, for any $m\in\N$,
\begin{equation}\label{eq:straight}
x^{\pm}_{i,m}\xi = \sum_{s\geq 0} \lambda^{\pm}_{i;s}(\xi)x^{\pm}_{i,m+s}
\end{equation}
\item For any $\xi,\eta\in Y^0$,
\begin{equation}\label{eq:straight mult}
\lambda^{\pm}_{i;s}(\xi\eta)=
\sum_{k+l = s} \lambda^{\pm}_{i;k}(\xi)\lambda^{\pm}_{i;l}(\eta)
\end{equation}
\item The operator $\lambda_{i;s}:Y^0\to Y^0$ is homogeneous of degree $-s$.
\item Let $\lambda^{\pm}_i(v): Y^0 \rightarrow Y^0[v]$ be given by
$$\lambda^{\pm}_i(v)(\xi)= \sum_{s\geq 0}\lambda^{\pm}_{i;s}(\xi)\,v^s$$
and extend the $\N$--grading on $Y^0$ to $Y^0[v]$ by $\deg(v)=1$.
Then $\lambda^{\pm}_i(v)$ is an algebra homomorphism of degree $0$.
\item $\lambda_{i_1}^{\epsilon_1}(v_1)$ and $\lambda_{i_2}^{\epsilon_2}(v_2)$
commute for any $i_1,i_2\in\bfI$ and $\epsilon_1,\epsilon_2\in\{\pm\}$.
\item For any $i\in\bfI$, $$\lambda^+_i(v)\lambda^-_i(v) = \Id = \lambda^-_i(v)\lambda^+_i(v)$$
\item For any $i,j\in\bfI$,
\[\lambda^{\pm}_j(v_1)\left(B_i(v_2)\right)=B_i(v_2) \mp \frac{q_i^{a_{ij}v_2} - q_i^{-a_{ij}v_2}}{v_2} e^{v_1v_2}\]
\item For any $i\in\bfI$ and $r\in\N$,
\[\lambda^{\pm}_j(v)(t_{i,r}) = t_{i,r}\mp d_ia_{ij}v^r\mod\hbar\]
\end{enumerate}
\end{prop}
\begin{pf}
(1)--(2) by Lemma \ref{le:Y23}, \eqref{eq:straight} holds when $\xi$ is one of the
generators $\xi_{j,r}$ of $Y^0$. Since \eqref{eq:straight} holds for $\xi\eta$ if it
holds for $\xi,\eta\in Y^0$, with $\lambda_{i;s}(\xi\eta)$ given by \eqref{eq:straight mult},
the $\lambda_{i;s}$ can be defined as operators on $Y^0$. The fact they are
uniquely characterised by \eqref{eq:straight} and satisfy \eqref{eq:straight mult}
follows from Corollary \ref{co:PBW}.

(3) the linear independence of the elements on the \rhs of \eqref{eq:straight}
implies that $\deg(\lambda_{i;s}(\xi))=\deg(\xi)-s$ for any homogeneous $\xi\in
Y^0$. (4) is a rephrasing of (2) and (3). (5) and (6) follow from (7) since the
elements $\{t_{i,n}\}$ generate $Y^0$. (7) follows from Lemma \ref
{commutation-with-log-like}. (8) is a direct consequence of (7).
\end{pf}

\begin{rem}
Using the shift operators, the relation \eqref{eq:straight} can be rewritten as
\[x^{\pm}_{i,m}\xi = \lambda_i^{\pm}(\sigma_i^{\pm})(\xi)\,x^{\pm}_{i,m}\]
\end{rem}
 
\section{Homomorphisms of geometric type}\label{section-nec-suff}

Let $\cYhg$ be the completion of $\Yhg$ \wrt its $\N$--grading. In this
section, we define an assignment
\[\Phi:\{H_{i,r}, E_{i,r}, F_{i,r}\}_{i\in\bfI,r\in\Z}\longrightarrow\cYhg\]
and find necessary and sufficient conditions for $\Phi$ to extend to a
homomorphism $\hloop\to\cYhg$.

\subsection{Definition of $\Phi$}\label{definition-assignment}

Define
\begin{align*}
\Phi(H_{i,0}) &= d_i^{-1}t_{i,0}
\intertext{and, for $r\in\Z^\times$}
\Phi(H_{i,r}) &=
\frac{\hbar}{q_i - q_i^{-1}} \sum_{k\geq 0} t_{i,k} \frac{r^k}{k!}=
\frac{B_i(r)}{q_i - q_i^{-1}}
\end{align*}
where $B_i(v)$  is the formal power series \eqref{inverse-borel}. Let $\wt{U}^0$
be the polynomial ring on the generators $\{H_{i,r}\}_{i\in\bfI,r\in\Z}$\footnote{$
\wt{U}^0$ is isomorphic to the subalgebra $U^0$ of $\hloop$ generated by $\{H
_{i,r}\}_{i\in\bfI,r\in\Z}$ by the PBW Theorem for $\hloop$ \cite{beck-pbw}, but
we shall not need this fact.}. The above assignment extends to an homomorphism
$\Phi^0:\wt{U}^0\to\wh{Y^0}$.

Let now $\{g_{i,m}^{\pm}\}_{i\in\bfI,m\in\N}$ be elements of $\widehat{Y^0}$,
and define further
\begin{gather*}
\Phi(E_{i,0}) = \sum_{m\geq 0} \fact{g^+_{i,m}}{m!} x^+_{i,m}\\
\Phi(F_{i,0}) = \sum_{m\geq 0} \fact{g^-_{i,m}}{m!} x^-_{i,m}
\end{gather*}
In terms of the shift operators $\sigma_i^\pm$, the above may be written as
\begin{align}\label{eq: phi EF}
\Phi(E_{i,0})&=g^+_i(\sigma_i^+)x^+_{i,0}\\
\Phi(F_{i,0})&=g^-_i(\sigma_i^-)x^-_{i,0}
\end{align}
where
\[g_i^{\pm}(v) = \sum_{m\geq 0}\fact{g^{\pm}_{i,m}}{m!}v^m\in\widehat{Y^0[v]}\]
with the completion of $Y^0[v]$ taken \wrt the $\N$--grading which extends that
on $Y^0$ by $\deg(v)=1$.

\subsection{Homomorphisms of geometric type}

If $\Phi:\hloop\to\wh{\Yhg}$ is a homomorphism of the above form, we shall say
that it is of {\em geometric type} since its form is related to the Chern character
in the geometric realisations of $\hloop$ and $\Yhg$ discussed in \ref{ss:belief}.

Such a homomorphism has the following properties
\begin{enumerate}
\item It restricts to a homomorphism $\hloopi{i}\to\wh{\Yhgi{i}}$ for any $\bfI$,
where $\hloopi{i}$ $\subset\hloop$ and $\Yhgi{i}\subset\Yhg$ are the subalgebras
generated by $\{E_{i,r},F_{i,r}, H_{i,r}\}_{r\in\Z}$ and $\{x_{i,k}^\pm,\xi_{i,k}\}
_{k\in\N}$ respectively.
\item It restricts to a homomorphism $\hloopb{\pm}\to\wh{\Yhb{\pm}}$, where
$\hloopb{+},\hloopb{-}$ $\subset\hloop$ and $\Yhb{\pm}\subset\Yhg$ are the
subalgebras generated respectively by $\{E_{i,r},H_{i,r}\}_{i\in\bfI,r\in\Z}$,
$\{F_{i,r},H_{i,r}\}_{i\in\bfI,r\in\Z}$ and $\{x^{\pm}_{i,k},\xi_{i,k}\}_{i\in\bfI,k\in\N}$ .
\end{enumerate}
Note however that, unless $g_{i,m}^\pm=0$ for any $i\in\bfI$ and $m\geq 1$,
$\Phi$ does not map the quantum group $U_\hbar\g=\<E_{i,0},F_{i,0},H_{i,0}
\>_{i\in\bfI}\subset\hloop$ to $U\g[[\hbar]]\subset\wh{\Yhg}$.

\subsection{}\label{ss:Phi QL23}

The following result shows that the requirement that $\Phi$ extends to an algebra
homomorphism determines its value on generators $E_{i,k}, F_{i,k}$.

\begin{prop}\label{assignment-arbitrary-nodes}
The assignment $\Phi$ is compatible with relations (QL2)--(QL3) if, and only if
\begin{align}
\Phi(E_{i,k})&=e^{k\sigma_i^+}g_i^+(\sigma_i^+)x^+_{i,0}\label{eq:Phi Ek}\\
\Phi(F_{i,k})&=e^{k\sigma_i^-}g_i^-(\sigma_i^-)x^-_{i,0}\label{eq:Phi Fk}
\end{align}
\end{prop}
\begin{pf}
We only consider the case of the $E's$. Let $i,j\in\bfI$ and $k\in\Z$. By (Y2),
\[\begin{split}
[\Phi(H_{i,0}),\Phi(E_{j,k})]
&=
[d_i^{-1}\xi_{i,0},e^{k\sigma_j^+}g_j^+(\sigma_j^+)x^+_{j,0}]\\
&=
e^{k\sigma_j^+}g_j^+(\sigma_j^+)[d_i^{-1}\xi_{i,0},x^+_{j,0}]\\
&=
a_{ij}\Phi(E_{j,k})
\end{split}\]
so that $\Phi$ is compatible with (QL2). Next, if $r\in\Z^\times$, Lemma \ref
{commutation-with-log-like} yields
\[\begin{split}
[\Phi(H_{i,r}),\Phi(E_{j,k})]
&= \frac{1}{q_i - q_i^{-1}} [B_i(r),e^{k\sigma_j^+}g_j^+(\sigma_j^+) x^+_{j,0}]\\
&= \frac{q_i^{ra_{ij}} - q_i^{-ra_{ij}}}{r(q_i - q_i^{-1})}
e^{r\sigma_j^+}e^{k\sigma_j^+} g_j^+(\sigma_j^+)x^+_{j,0}\\
&= \frac{[ra_{ij}]_{q_i}}{r} \Phi(E_{j,r+k})
\end{split}\]
and $\Phi$ is compatible with (QL3).

Conversely, if $\Phi$ is compatible with (QL3) then $\Phi(E_{i,r})=r/[2r]_{q_i}
[\Phi(H_{i,r}),\Phi(E_{i,0})]$ for $r\neq 0$ and the computation above shows
that this is equal to $e^{r\sigma_i^+}\Phi(E_{i,0})$.
\end{pf}

\subsection{Necessary and sufficient conditions}\label{nec-suff}

Let $\lambda^{\pm}_i(v):Y^0\rightarrow Y^0[v]$ be the homomorphism defined
in Proposition \ref{raising-lowering-operators}.

\begin{thm}\label{first-main-theorem}
The assignment $\Phi$ extends to an algebra homomorphism $\hloop\to\wh
{\Yhg}$ if, and only if the following conditions hold
\begin{enumerate}
\item[(A)] For any $i,j\in\bfI$
\[g_i^+(u)\lambda^+_i(u)(g_j^-(v)) = g_j^-(v)\lambda^-_j(v)(g_i^+(u))\]
\item[(B)] For any $i\in\bfI$ and $k\in \Z$
\[\left.e^{ku}g_i^+(u)\lambda^+_i(u)(g_i^-(u))\right |_{u^m = \xi_{i,m}} =
\Phi^0\left(\frac{\psi_{i,k} - \phi_{i,k}}{q_i - q_i^{-1}}\right)\]
\item[(C)] For any $i,j\in\bfI$ and $a=d_ia_{ij}/2$
\[g_i^{\pm}(u)\lambda^{\pm}_i(u)(g_j^{\pm}(v))
\left(\frac{e^u - e^{v\pm a\hbar}}{u-v\mp a\hbar}\right)=
g_j^{\pm}(v)\lambda^{\pm}_j(v)(g_i^{\pm}(u))
\left(\frac{e^v-e^{u\pm a\hbar}}{v-u\mp a\hbar}\right)\]
\end{enumerate}
\end{thm}
\begin{pf}
By construction and Proposition \ref{assignment-arbitrary-nodes}, $\Phi$ is compatible
with the relations (QL1)--(QL3). The result then follows from Lemmas \ref{le:AB} and
\ref{le:C} below and the proof of the $q$--Serre relations (Proposition \ref{pr:Serre} in
the Appendix).
\end{pf}

\subsection{} 

\begin{lem}\label{le:AB}
$\Phi$ is compatible with the relation (QL5) if, and only if (A) and (B) hold.
\end{lem}
\begin{pf}
Compatibility with (QL5) reads
\[ [\Phi(E_{i,k}),\Phi(F_{j,l})]=
\delta_{ij} \Phi^0\left(\frac{\psi_{i,k+l} - \phi_{i,k+l}}{q_i - q_i^{-1}}\right)\]
for $i,j\in\bfI$ and $k,l\in\Z$. We begin by computing the left--hand side.
To this end, it will be convenient to write formulae \eqref{eq:Phi Ek}--\eqref
{eq:Phi Fk} as
\[\Phi(E_{i,k})=\sum_{m\geq 0}\fact{g_{i,m}^{+,(k)}}{m!}\,x_{i,m}^+
\aand
\Phi(F_{i,k})=\sum_{m\geq 0}\fact{g_{i,m}^{-,(k)}}{m!}\,x_{i,m}^-\]
where $g_{i,m}^{\pm,(k)}\in\wh{Y^0}$ are defined by
$\sum_{m\geq 0}\fact{g_{i,m}^{\pm,(k)}}{m!}v^m=e^{kv}g_i^{\pm}(v)$.
This yields
\[\begin{split}
\Phi(E_{i,k})\Phi(F_{j,l})
&= \sum_{m,n\geq 0} \fact{g_{i,m}^{+,(k)}}{m!}x_{i,m}^+\fact{g_{j,n}^{-,(l)}}{n!} x_{j,n}^-\\
&= \sum_{m,n,s\geq 0} \fact{g_{i,m}^{+,(k)}}{m!} \lambda_{i,s}^+\left(\fact{g_{j,n}^{-,(l)}}{n!}\right) x_{i,m+s}^+ x_{j,n}^-\\
&= \sum_{m,n,s\geq 0} \fact{g_{i,m}^{+,(k)}}{m!} \lambda_{i,s}^+\left(\fact{g_{j,n}^{-,(l)}}{n!}\right)
\left(x_{j,n}^-x_{i,m+s}^++\delta_{ij}\xi_{i,m+n+s}\right)
\end{split}\]
where we used (Y5). Similarly,
$\Phi(F_{j,l})\Phi(E_{i,k})
=\sum_{m,n,s}
\fact{g_{j,m}^{-,(l)}}{m!} \lambda_{j,s}^-\left(\fact{g_{i,n}^{+,(k)}}{n!}\right) x_{j,m+s}^- x_{i,n}^+$. 
Define $R^{(k,l)},L^{(k,l)}\in\wh{Y^0}[[u,v]]$ by
\begin{align*}
R^{(k,l)}
&= 
\sum_{m\geq 0}g_{i,m}^{+,(k)}u^m
\sum_{s\geq 0} \lambda_{i,s}^+u^s
\left(\sum_{n\geq 0} \fact{g_{j,n}^{-,(l)}}{n!} v^n\right)
=
e^{ku}e^{lv}g_i^+(u)\lambda_i^+(u)(g_j^-(v))\\
L^{(k,l)}
&=
\sum_{m\geq 0}g_{j,m}^{-,(l)}v^m
\sum_{s\geq 0} \lambda_{j,s}^-v^s
\left(\sum_{n\geq 0} \fact{g_{i,n}^{+,(k)}}{n!} u^n\right) 
=
e^{ku}e^{lv}g_j^-(v)\lambda_j^-(v)(g_i^+(u))
\end{align*}
By the PBW Theorem \ref{PBW-yangian}, $\Phi$ is compatible with (QL5)
if, and only if $R^{(k,l)}=L^{(k,l)}$ and, for $i=j$,
\[\left.R^{(k,l)}\right|_{u^mv^n=\xi_{i,m+n}}=\Phi^0\left(\frac{\psi_{i,k+l} - \phi_{i,k+l}}{q_i - q_i^{-1}}\right)\]
The first equation is clearly equivalent to (A) and the second to (B).
\end{pf}

\subsection{} 

\begin{lem}\label{le:C}
$\Phi$ is compatible with the relation (QL4) if, and only if (C) holds.
\end{lem}
\begin{pf}
We prove the claim for the $E$'s only. Compatibility with (QL4) reads
\[\Phi(E_{i,k+1})\Phi(E_{j,l}) - q_i^{a_{ij}}\Phi(E_{i,k})\Phi(E_{j,l+1})
=q_i^{a_{ij}} \Phi(E_{j,l})\Phi(E_{i,k+1}) - \Phi(E_{j,l+1})\Phi(E_{i,k})\]
for any $i,j\in\bfI$ and $k,l \in \Z$. Assume first $i\neq j$ and set
$a=d_ia_{ij}/2$,  so that $q_i^{a_{ij}}=e^{a\hbar}$. Since
\[\Phi(E_{i,r})\Phi(E_{j,s}) =
e^{r\sigma^+_i}e^{s\sigma^+_j}
g_i^+(\sigma^+_i)\lambda^+_i(\sigma^+_i)\left(g_j^+(\sigma^+_j)\right)
x^+_{i,0}x^+_{j,0}\]
the above reduces to
\begin{multline*}
e^{k\sigma^+_i}e^{l\sigma^+_j}
g_i^+(\sigma^+_i)\lambda^+_i(\sigma^+_i)(g_j^+(\sigma^+_j))
\left(e^{\sigma^+_i} - e^{\sigma^+_j+a\hbar}\right)x^+_{i,0}x^+_{j,0}\\
 = e^{k\sigma^+_i}e^{l\sigma^+_j}
 g_j^+(\sigma^+_j)\lambda^+_j(\sigma^+_j)(g_i^+(\sigma^+_i))
 \left(e^{\sigma^+_i+a\hbar} - e^{\sigma^+_j}\right)x^+_{j,0}x^+_{i,0}
\end{multline*}
Using (1) of Proposition \ref{algebraic-properties-shift}, we get
\[\begin{split}
\left(e^{\sigma^+_i} - e^{\sigma^+_j+a\hbar}\right)x^+_{i,0}x^+_{j,0}
&=
\frac{e^{\sigma^+_i} - e^{\sigma^+_j+a\hbar}}{\sigma^+_i-\sigma^+_j-a\hbar}
\left(\sigma^+_i-\sigma^+_j-a\hbar\right)x^+_{i,0}x^+_{j,0}\\
&=
\frac{e^{\sigma^+_i} - e^{\sigma^+_j+a\hbar}}{\sigma^+_i-\sigma^+_j-a\hbar}
\left(\sigma^+_i-\sigma^+_j+a\hbar\right)x^+_{j,0}x^+_{i,0}
\end{split}\]
The PBW Theorem \ref{PBW-yangian} then shows that the above is equivalent to (C).

Assume now that $i=j$, then
\[\begin{split}
\Phi(E_{i,r})\Phi(E_{i,s})
&=
\left(e^{r\sigma_i^+} g_i^+(\sigma_i^+) x^+_{i,0}\right)
\left(e^{s\sigma_i^+} g_i^+(\sigma_i^+) x^+_{i,0}\right)\\
&= 
\mu\left(e^{r\sigma^+_{i,(1)}}e^{s\sigma^+_{i,(2)}}g_i^+(\sigma^+_{i,(1)})\lambda^+_i(\sigma^+_{i,(1)})
\left(g_i^+(\sigma^+_{i,(2)})\right) x^+_{i,0}\otimes x^+_{i,0}\right) 
\end{split}\]
The compatibility with (QL4) therefore reduces to
\begin{multline*}
\mu\left(e^{k\sigma^+_{i,(1)}}e^{l\sigma^+_{i,(2)}}
g_i^+(\sigma^+_{i,(1)})\lambda^+_i(\sigma^+_{i,(1)})(g_i^+(\sigma^+_{i,(2)}))
\left(e^{\sigma^+_{i,(1)}} - e^{\sigma^+_{i,(2)} + d_i\hbar}\right) x^+_{i,0}\otimes x^+_{i,0}\right) \\
=
\mu\left(e^{l\sigma^+_{i,(1)}}e^{k\sigma^+_{i,(2)}}
g_i^+(\sigma^+_{i,(1)})\lambda^+_i(\sigma^+_{i,(1)})(g_i^+(\sigma^+_{i,(2)}))
\left(e^{\sigma^+_{i,(2)}+d_i\hbar} - e^{\sigma^+_{i,(1)}}\right)x_{i,0}^+\otimes x^+_{i,0}\right)
\end{multline*}
that is, to
\begin{multline*}
 \mu\left(\left(e^{k\sigma^+_{i,(1)}}e^{l\sigma^+_{i,(2)}} + e^{l\sigma^+_{i,(1)}}e^{k\sigma^+_{i,(2)}}\right)
 \right.\\\left.
 g_i^+(\sigma^+_{i,(1)})\lambda^+_i(\sigma^+_{i,(1)})(g_i^+(\sigma^+_{i,(2)}))
 \left(e^{\sigma^+_{i,(1)}} - e^{\sigma^+_{i,(2)} + d_i\hbar}\right) x^+_{i,0}\otimes x^+_{i,0}\right)  = 0
\end{multline*}
By (2) of Proposition \ref{algebraic-properties-shift} and Corollary \ref{co:B(u,v)}, this
equation is equivalent to the requirement that
\[g_i^+(u)\lambda^+_i(u)(g_i^+(v))\left(\frac{e^u - e^{v+d_i\hbar}}{u-v-d_i\hbar}\right)\]
be symmetric under $u\leftrightarrow v$, which is precisely condition (C) for $i=j$.
\end{pf}

\subsection{}

For later use, we shall need the following
\begin{lem}\label{le:constant term}
Let $\{g_i^{\pm}(u)\}_{i\in\bfI}\subset\wh{Y^0[u]}$ be elements satisfying condition
(B) of Theorem \ref{first-main-theorem}. Then,
\[g_i^{\pm}(u) = \ds \frac{1}{d_i^{\pm}} \mod \wh{Y^0[u]}_+\]
where $\{d_i^{\pm}\}_{i\in\bfI}\subset\C^\times$ satisfy $d_i^+d_i^- = d_i$ for each
$i\in\bfI$. In particular, each $g_i^{\pm}(u)$ is invertible. 
\end{lem}
\begin{pf}
Condition (B) for $k=0$ yields
\[\left.g_i^+(u)\lambda_i^+(u)\left(g_i^-(u)\right)\right|_{u^m = \xi_{i,m}} =
\Phi^0\left(\frac{e^{\frac{\hbar d_i}{2}H_{i,0}}-e^{-\frac{\hbar d_i}{2}H_{i,0}}}{q_i - q_i^{-1}}\right)\]
Computing mod $\hbar$, and {\it a fortiori} mod $\wh{Y^0[u]}_+$, yields
\[\Phi^0\left(\frac{e^{\frac{\hbar d_i}{2}H_{i,0}}-e^{-\frac{\hbar d_i}{2}H_{i,0}}}{q_i - q_i^{-1}}\right)
=
\Phi^0\left(H_{i,0}\right)
=
d_i^{-1}t_{i,0}\]
Write $g_i^\pm(u)=p_i^\pm\mod\wh{Y^0[u]}_+$, where $p_i^\pm\in\C[t_{j,0}]_{j\in\bfI}$.
Computing mod $\wh{Y^0[u]}_+$, we get
\[\left.g_i^+(u)\lambda_i^+(u)\left(g_i^-(u)\right)\right|_{u^m = \xi_{i,m}}=
p_i^+(t_{j,0})\lambda_{i;0}^+(p_i^-(t_{j,0}))\xi_{i,0}=
p_i^+(t_{j,0})p_i^-(t_{j,0}-d_ia_{ij})\xi_{i,0}\]
where we used (8) of  Proposition \ref{raising-lowering-operators}.
Comparing both sides and using $\xi_{i,0}=t_{i,0}$ yields the claim.
\end{pf}

\section{Existence of homomorphisms}\label{section-existence}

In this section, we construct an explicit homomorphism $\hloop\to\wh
{\Yhg}$ by exhibiting a joint solution to equations (A)--(C) of Theorem
\ref{first-main-theorem}. We begin by giving an intrinsic expression for
the \rhs of equation (B) (Proposition \ref{pr:condition-B}). Until \ref
{ss:existence}, we fix $i\in\bfI$, and consider the subalgebras $\wt{U}^0
_i\subset \wt{U}^0$, $Y_i^0\subset Y^0$ generated by $\{H_{i,r}\}_{r\in
\Z}$ and $\{\xi_{i,r}\}_{r\in\N}$ respectively.

\subsection{The functions $\mathbf{G(v)}$ and $\mathbf{\gamma_i(v)}$}\label{ss:G(v)}

Consider the formal power series
\begin{equation}\label{eq:G}
G(v) = \log\left(\frac{v}{e^{v/2} - e^{-v/2}}\right)\in v\Q[[v]]
\end{equation}
Define $\gamma_i(v)\in\wh{Y^0_i[v]}_+$ by
\[\gamma_i(v)=
\hbar\sum_{r\geq 0}\frac{t_{i,r}}{r!} (-\partial_v)^{r+1}\,G(v)\]
Recall that $B_i(v)=\hbar\sum_{r\geq 0} t_{i,r}\frac{v^r}{r!}$ is the inverse
Borel transform of $t_i(u)$. This allows us to write $\gamma_i(u)$ more
compactly as
\begin{equation}\label{gamma-function-sl2-compact}
\gamma_i(v) = -B_i(-\partial_v)G'(v)
\end{equation}

\subsection{}

\begin{prop}\label{pr:condition-B}
The following holds in $\wh{Y_i^0}$ for any $k\in\Z$
\[\Phi^0\left(\frac{\psi_{i,k}-\phi_{i,k}}{q_{i}-q_{i}^{-1}}\right)=
\frac{\hbar}{q_{i}-q_{i}^{-1}}\left.e^{kv}\exp\left(\gamma_i(v)\right)\right|_{v^n = \xi_{i,n}}\]
\end{prop}

\noindent The above identity will be proved in \ref{ss:polynomials}--\ref
{ss:rewriting of B} by injectively mapping both sides to a family of polynomial
rings, and verifying their equality there.

\subsection{Universal Drinfeld polynomials}\label{ss:polynomials}

Fix an integer $m\geq 1$. Following \cite{nakajima-qaffine,varagnolo-yangian},
consider the rings
\begin{gather*}
S(m) = \C[q^{\pm 1},A_1^{\pm 1}, \ldots,A_m^{\pm 1}]^{\Sym_m}\\[1.3ex]
R(m) = \C[\hbar,a_1,\dots,a_m]^{\Sym_m}
\end{gather*}
Define a homomorphism $\D^U:\wt{U}^0_i \rightarrow S(m)$ by
\begin{equation}\label{eq:DU}
\D^U(\psi_i(z)) = \prod_{p=1}^m \frac{q_{i}z - q_{i}^{-1}A_p}{z-A_p} = \D^U(\phi_i(z))
\end{equation}
where the first and second equalities are obtained by expanding
the middle term in powers of $z^{-1}$ and $z$ respectively. Similarly,
define $\D^Y: Y^0_i \rightarrow R(m)$ by
\begin{equation}\label{eq:DY}
\D^Y(\xi_i(u)) = \prod_{p=1}^m \frac{u+d_i\hbar-a_p}{u-a_p}
\end{equation}
The homomorphism $\D^U$ (resp. $\D^Y$) gives the action of
$\psi_i(z),\phi_i(z)$ (resp. $\xi_i(u)$) on the highest weight vector
of the indecomposable $U_{\hbar}(L\Lsl_2)$ (resp. $Y_\hbar(\Lsl
_2)$) module with Drinfeld polynomial $(1-A_1^{-1}z)\cdots(1-A_
m^{-1}z)$ (resp. $(u-a_1)\cdots(u-a_m)$).

\subsection{} 

The following result spells out the image of the generators of $\wt{U}^0_i$
and $Y^0_i$ under $\D^U$ and $\D^Y$ respectively.

\begin{prop}\label{useful-1}\hfill
\begin{enumerate}
\item The following holds: $\D^U(\psi_{i,0})=q_i^m=\D^U(\phi_{i,0})^{-1}$
and, for any $r\in\N^*$
\begin{align}
\D^U(\psi_{i,r})
&=\phantom{-}(q_{i}-q_{i}^{-1}) 
\sum_{p=1}^m A_p^{r\phantom{-}}\prod_{p'\neq p}
\frac{q_{i}A_p - q_{i}^{-1}A_{p'}}{A_p - A_{p'}}
\label{eq:DU psi}
\\
\D^U(\phi_{i,-r})
&=-(q_{i}-q_{i}^{-1})
\sum_{p=1}^m A_p^{-r}\prod_{p'\neq p}
\frac{q_{i}A_p - q_{i}^{-1}A_{p'}}{A_p - A_{p'}}
\label{eq:DU phi}
\end{align}
Moreover, $\D^U$ maps $H_{i,0}$ to $m$ and, for any $r\in\Z^*$,
\begin{equation}\label{eq:DU H}
\D^U((q_{i}-q_{i}^{-1})H_{i,r}) = \frac{1-q_{i}^{-2r}}{r} \sum_{p=1}^m A_p^r
\end{equation}
\item The homomorphism $\D^Y$ maps $\xi_{i,0}$ to $d_i m$ and, for any
$r\in\N$
\begin{equation}\label{eq:DY xir}
\D^Y(\xi_r)=d_i\sum_{p=1}^m a_p^r \prod_{p'\neq p}
\frac{a_p - a_{p'}+d_i\hbar}{a_p - a_{p'}}
\end{equation}
Moreover, for any $r\in\N$,
\begin{equation}\label{eq:DY tr}
\D^Y(t_r)=\frac{1}{r+1}\sum_{p=1}^m\frac{a_p^{r+1}-(a_p-d_i\hbar)^{r+1}}{\hbar}
\end{equation}
\item If $B_i(v)\in Y^0_i[[v]]$ is the series defined by \eqref{inverse-borel},
then
\begin{equation}\label{eq:DY Bt}
\D^Y(B_i(v)) = \frac{1-e^{-d_i\hbar v}}{v} \sum_{p=1}^m e^{a_pv}
\end{equation}
\end{enumerate}
\end{prop}

\begin{pf}
(1) The fact that $\D^U(\psi_{i,0})=q_i^m=\D^U(\phi_{i,0})^{-1}$ follows
by taking the values of the middle term $P(z)$ in \eqref{eq:DU} at $z=\infty$
and $z=0$ respectively. Next, the partial fraction decomposition of $P(z)$
is readily seen to be
\begin{equation}\label{eq:partial P}
\prod_{p=1}^m
\frac{q_{i}z - q_{i}^{-1}A_p}{z-A_p}
=q_{i}^m+(q_{i}-q_{i}^{-1})\sum_{p=1}^mA_p
\left(\prod_{p'\neq p}
\frac{q_{i}A_p-q_{i}^{-1}A_{p'}}{A_p-A_{p'}}\right)\frac{1}{z-A_p}
\end{equation}
The relations \eqref{eq:DU psi}--\eqref{eq:DU phi} follow by expanding
this into powers of $z^{-1}$ and $z$ respectively. For later use, note
that the evaluation of \eqref{eq:partial P} at $z=0$ yields the identity
\begin{equation}\label{eq:for later}
\D^U(\psi_{i,0}-\phi_{i,0})=
q_i^m-q_i^{-m}=
(q_{i}-q_{i}^{-1})\sum_{p=1}^m
\prod_{p'\neq p}
\frac{q_{i}A_p-q_{i}^{-1}A_{p'}}{A_p-A_{p'}}
\end{equation}

Since $\D^U(\psi_{i,0})=q_{i}^m$, it follows that $\D^U(H_{i,0})=m$, and
\[
\D^U\left(\exp\left((q_{i}-q_{i}^{-1})\sum_{s\geq 1} H_{i,s}z^{-s}\right)\right)=
\D^U(\psi_{i,0}^{-1}\psi_i(z))=
\prod_{p=1}^m \frac{z-q_{i}^{-2}A_p}{z-A_p}
\]
taking the $\log$ of both sides and expanding in powers of $z^{-1}$ (resp.
$z$) yields \eqref{eq:DU H} for $r>0$ (resp. $r<0$).

(2) The fact that $\D^Y(\xi_{i,0})=d_i m$ follows by taking the coefficient of
$u^{-1}$ in \eqref{eq:DY}. The partial fraction decomposition of $\D^Y(\xi_i
(u))$ is
\[\prod_{p=1}^m \frac{u+d_i\hbar-a_p}{u-a_p}=
1+
d_i\hbar\sum_{p=1}^m\left(\prod_{p'\neq p}\frac{a_p-a_{p'}+
d_i\hbar}{a_p-a_{p'}}\right)\frac{1}{u-a_p}\]
and \eqref{eq:DY xir} follows by taking the coefficient of $u^{-r-1}$. Taking
the $\log$ of both sides of \eqref{eq:DY} yields
\begin{equation}\label{eq:DY t(u)}
\D^Y(t_i(u)) = 
\sum_p -\log\left(1-a_pu^{-1}\right) + 
\log\left(1-(a_p-d_i\hbar)u^{-1}\right)
\end{equation}
and therefore \eqref{eq:DY tr}.

(3) Follows by applying \eqref{eq:B log} to \eqref{eq:DY}.
\end{pf}

\begin{cor}\label{co:asymptotic inj}
The homomorphism $\D^Y:Y^0_i\to \bigoplus_{m\geq 1}R(m)$ is injective.
\end{cor}
\begin{pf}
This follows from \eqref{eq:DY tr} and the fact that the power sums $\sum
_p a_p^r$ are algebraically independent.
\end{pf}

\subsection{}

Let $\wh{R(m)}$ be the completion of $R(m)$ \wrt the $\N$--grading
defined by $\deg(\hbar)=\deg(a_p)=1$. Since the map $\D^Y:Y^0_i\to
R(m)$ preserves the grading, it extends to a homomorphism $\wh{Y
^0_i}\to\wh{R(m)}$.

\begin{cor}\label{compatible-drinfeld-homomorphism}
Let $\eexp:S(m)\rightarrow R(m)$ be the homomorphism defined by
\[q \longmapsto e^{\hbar/2}\aand A_p \longmapsto e^{a_p}\]
Then, the following diagram commutes
\[
\xymatrix{
\wt{U}^0_i \ar[rr]^{\D^U} \ar[d]_{\Phi^0} && S(m) \ar[d]^{\eexp} \\
\widehat{Y^0_i} \ar[rr]^{\D^Y} && \widehat{R(m)}
}
\]
where $\Phi^0$ is defined in Section \ref{definition-assignment}.
\end{cor}
\begin{pf}
It suffices to check the commutativity on the generators $\{H_{i,r}\}_{r\in\Z}$
of $\wt{U}^0_i$. The statement now follows from \eqref{eq:DU H}, \eqref{eq:DY Bt}
and the fact that, for $r\neq 0$, $\Phi^0(H_{i,r})=\left.B_i(v)/(q_{i}-q_{i}^{-1})
\right|_{v=r}$.
\end{pf}

\subsection{Proof of Proposition \ref{pr:condition-B}}\label{ss:rewriting of B}

By Corollaries \ref{co:asymptotic inj} and \ref{compatible-drinfeld-homomorphism},
it suffices to prove that, for any $k\in\Z$,
\[
\eexp\left(\D^U\left(\frac{\psi_{i,k} - \phi_{i,k}}{q_{i}-q_{i}^{-1}}\right)\right)=
\frac{\hbar}{q_{i}-q_{i}^{-1}}
\D^Y\left(\left.e^{kv}\exp\left(\gamma_i(v)\right)\right |_{v^n=\xi_n}\right)
\]
By \eqref{eq:DU psi}--\eqref{eq:DU phi} and \eqref{eq:for later}, the \lhs is
equal to
\[\eexp\left( \sum_p A_p^k \prod_{p'\neq p} \frac{q_{i}A_p - q_{i}^{-1}A_{p'}}{A_p-A_{p'}} \right)\]
We now compute the right--hand side. By \eqref{gamma-function-sl2-compact}
and \eqref{eq:DY Bt},
\[\D^Y(\gamma_i(v))
=\sum_p e^{-a_p\partial_v}\frac{1-e^{d_i\hbar\partial_v}}{\partial_v}\partial_v G(v)
= \sum_p G(v-a_p) - G(v-a_p+d_i\hbar)\]
where we used $e^{a\partial_v}G(v)=G(v+a)$. Thus,
\[\D^Y\left(e^{kv}\exp\left(\gamma_i(v)\right)\right)=
e^{kv} \prod_p \frac{v-a_p}{v-a_p+d_i\hbar}
\frac{q_{i}e^v - q_{i}^{-1}e^{a_p}}{e^v - e^{a_p}}\]
By \eqref{eq:DY xir}, the substitution $v^n = \xi_{i,n}$ in a formal power series 
$F\in\wh{Y_i^0[v]}$ gives
\[\D^Y(F(v)|_{v^n=\xi_n})=
d_i\sum_p \D^Y(F(a_p))\prod_{p'\neq p} \frac{a_p-a_{p'}+d_i\hbar}{a_p - a_{p'}}\]
Since $\left.(v-a_p)/(e^v-e^{a_p})\right|_{v=a_p}=e^{-a_p}$ and
$\left.(q_{i}e^v - q_{i}^{-1}e^{a_p})/(v-a_p+d_i\hbar)\right|_{v=a_p}=$
$e^{a_p}(q_{i}-q_{i}^{-1})/d_i\hbar$, this implies that
\begin{equation*}
\begin{split}
\D^Y\Bigl(e^{kv} &\exp\left.\left(\gamma_i(v)\right)\right|_{v^n = \xi_{i,n}}\Bigr)\\
&=
\frac{q_{i}-q_{i}^{-1}}{\hbar}
\sum_p e^{ka_p}\prod_{p'\neq p}
\frac{a_p - a_{p'}}{a_p-a_{p'}+d_i\hbar} \frac{q_{i}e^{a_p} - q_{i}^{-1}e^{a_{p'}}}{e^{a_p} - e^{a_{p'}}} \frac{a_p-a_{p'}+d_i\hbar}{a_p-a_{p'}}\\
&=
\frac{q_{i}-q_{i}^{-1}}{\hbar}
\sum_p e^{ka_p}\prod_{p'\neq p}
\frac{q_{i}e^{a_p} - q_{i}^{-1}e^{a_{p'}}}{e^{a_p} - e^{a_{p'}}}\\
\end{split}
\end{equation*}
as claimed.

\subsection{A joint solution to equations (A)--(C)}\label{ss:existence}

Proposition \ref{pr:condition-B} suggests replacing equation (B) of Theorem
\ref{first-main-theorem} by the stronger requirement that, for any $i\in\bfI$
\begin{equation}\label{eq:B tilde}
g_i^+(v)\lambda^+_i(v)(g_i^-(v)) =
\frac{\hbar}{q_{i}-q_{i}^{-1}}\exp\left(\gamma_i(v)\right)
\tag{$\wt{B}$}
\end{equation}

Using equation (A) for $j=i$ and $u=v$ then shows that
\[g_i^-(v)\lambda^-_i(v)(g_i^+(v)) =
\frac{\hbar}{q_{i}-q_{i}^{-1}}\exp\left(\gamma_i(v)\right)\]
Applying the operator $\lambda^-_i(v)$ to the first of these equations, and
using $\lambda^+_i(v)\lambda^-_i(v)=\Id$ (Proposition \ref{raising-lowering-operators})
and the second equation, yields $\lambda^-_i(v)(\gamma_i(v))=\gamma_i(v)$.
Similarly, applying $\lambda^+_i(v)$ to the second of these equations yields
$\lambda^+_i(v)(\gamma_i(v))=\gamma_i(v)$. This suggests in turn that a
solution to the above equations may be given by $g_i^\pm(v)=g_i(v)$, where
\begin{equation}\label{eq:intro solution}
g_i(v)=
\left(\frac{\hbar}{q_i-q_i^{-1}}\right)^{1/2}
\exp\left(\frac{\gamma_i(v)}{2}\right)
\end{equation}
We now show that this is indeed the case.

\begin{thm}\label{third-main-theorem}
The series $g_i^\pm(v)=g_i(v)$ satisfy the conditions (A),($\wt{B}$)
and (C) of Theorem \ref{first-main-theorem}, and therefore give rise
to a homomorphism $\Phi:\hloop\to\cYhg$.
\end{thm}

\subsection{}

We shall need the following

\begin{lem}\label{le:lambda g}
Let $i,j\in\bfI$, and set $a=d_ia_{ij}/2$. Then,
\[\lambda^{\pm}_i(u)\left(g_j(v)\right)=
g_j(v) \exp\left(\pm\frac{G(v-u+a\hbar) - G(v-u-a\hbar)}{2}\right)\]
where $G(v)$ is given by \eqref{eq:G}.
\end{lem}
\begin{pf}
By Proposition \ref{raising-lowering-operators},
\[\lambda^{\pm}_i(u)(B_j(v))=
B_j(v) \mp \frac{e^{a\hbar v} - e^{-a\hbar v}}{v} e^{uv}\]
Since $\gamma_j(v)=-B_j(-\partial_v)\partial_v G(v)$, we get
\[\begin{split}
\lambda^{\pm}_i(u)\gamma_j(v)
&=
\gamma_j(v)\pm
\frac{e^{a\hbar\partial_v}-e^{-a\hbar\partial_v}}{\partial_v}e^{-u\partial_v} \partial_v G(v)\\
&=
\gamma_j(v)\pm\left(G(v-u+a\hbar)-G(v-u-a\hbar)\right)
\end{split}\]
The claim follows by exponentiating.
\end{pf}

\subsection{Proof of condition (A)}

We need to prove that for every $i,j\in\bfI$, we have
\[g_i(u)\lambda^+_i(u)(g_j(v)) = g_j(v)\lambda^-_j(v)(g_i(u))\]
By Lemma \ref{le:lambda g}, this is equivalent to
\begin{multline*}
g_i(u)g_j(v) \exp\left(\frac{G(v-u+a\hbar) - G(v-u-a\hbar)}{2}\right)\\
= g_i(u)g_j(v) \exp\left(\frac{G(u-v-a\hbar) - G(u-v+a\hbar)}{2}\right)
\end{multline*}
The result now follows since $G$ is an even function.

\subsection{Proof of condition ($\wt{B}$)}

Lemma \ref{le:lambda g} implies that
\begin{equation*}
\begin{split}
g_i(u)\lambda^+_i(u)(g_i(u))
&= 
g_i(u)^2 \exp\left(\frac{G(d_i\hbar) - G(-d_i\hbar)}{2}\right)\\
&= 
g_i(u)^2\\
&= \frac{\hbar}{q_i - q_i^{-1}} \exp\left(\gamma_i(u)\right)
\end{split}
\end{equation*}
where the second equality holds because $G$ is even.

\subsection{Proof of condition (C)}

Let $i,j\in\bfI$ and set $a = d_ia_{ij}/2$. We need to prove that
\[
g_i(u)\lambda^{\pm}_i(u)(g_j(v))\frac{e^u-e^{v\pm a\hbar}}{u-v\mp a\hbar}=
g_j(v)\lambda^{\pm}_j(v)(g_i(u))\frac{e^v-e^{u\pm a\hbar}}{v-u\mp a\hbar}
\]
By Lemma \ref{le:lambda g} and the fact that $G$ is even, we get the following
equivalent assertion
\[\exp\left(G(v-u\pm a\hbar) - G(v-u\mp a\hbar)\right)
\frac{e^u-e^{v\pm a\hbar}}{u-v\mp a\hbar}=
\frac{e^v - e^{u\pm a\hbar}}{v-u\mp a\hbar}\]
Using the definition of $G$, the above becomes the equality
\[
\left(\frac{v-u\pm a\hbar}{e^{v\pm a\hbar} - e^u}\right)
\left( \frac{e^v - e^{u\pm a\hbar}}{v-u\mp a\hbar}\right)
\left(\frac{e^u-e^{v\pm a\hbar}}{u-v\mp a\hbar}\right)=
\frac{e^v - e^{u\pm a\hbar}}{v-u\mp a\hbar}
\]

\section{Uniqueness of homomorphisms}\label{section-uniqueness}

The aim of this section is to prove that homomorphisms of geometric
type are unique up to conjugation and scaling.

\subsection{}

Let $\G$ be the set of solutions $\fg=\{g_i^{\pm}(u)\}_{i\in\bfI}$ of equations
(A)--(C) of Theorem \ref{first-main-theorem}. Given a collection $\fr=\{r_i^
{\pm}(u)\}_{i\in\bfI}$ of invertible elements  of $\widehat{Y^0[u]}$, set
$$\fr\cdot\fg = \{r_i^{\pm}(u)\cdot g_i^{\pm}(u)\}_{i\in\bfI}$$

\begin{lem}\label{gauge-1}
Let $\fg\in\G$. Then, $\fr\cdot\fg\in \G$ if and only if the following
holds
\begin{itemize}
\item[(A$_0$)] For any $i,j\in\bfI$,
$$r_i^+(u)\lambda^+_i(u)(r_j^-(v)) = r_j^-(v)\lambda^-_j(v)(r_i^+(u))$$
\item[(B$_0$)] For any $i\in\bfI$,
$$r_i^+(u)\lambda^+_i(u)(r_i^-(u)) = 1  = r_i^-(u)\lambda^-_i(u)(r_i^+(u))$$
\item[(C$_0^\pm$)] For any $i,j\in\bfI$,
$$r_i^{\pm}(u)\lambda^{\pm}_i(u)(r_j^{\pm}(v))=
r_j^{\pm}(v)\lambda^{\pm}_j(v)(r_i^{\pm}(u))$$
\end{itemize}
\end{lem}
\begin{pf} Let $\fh=\fr\cdot\fg$. The following assertions are straightforward
to check
\begin{itemize}
\item $\fh$ satisfies (A) if and only if $\fr$ satisfies (A$_0$).
\item $\fh$ satisfies (B) if $\fr$ satisfies (B$_0$).
\item $\fh$ satisfies (C) if and only if $\fr$ satisfies (C$^\pm_0$).
\end{itemize}
There remains to prove that if $\fh$ lies in $\G$, then $\fr$ satisfies (B$_0$).

We claim that (A$_0$) and (C$_0^\pm$) imply that $c_i(u)=r^+_i(u)\lambda^+
_i(u)(r^-_i(u))$ lies in $\C[[\hbar,u]]$. Assuming this, write $c_i(u)=\sum_n c_i
^{(n)}u^n$, where $c_i^{(n)}\in\C[[\hbar]]$. Then,
\begin{equation*}
\begin{split}
\left.\left(h^+_i(u)\lambda^+_i(u)(h^-_i(u))\right)\right|_{u^m = \xi_{i,m}}
&=
\left.\left(c_i(u)g^+_i(u)\lambda^+_i(u)(g^-_i(u))\right)\right|_{u^m = \xi_{i,m}}\\[1.1 ex]
&=
\sum_{n\geq 0}c_i^{(n)}\left.\left(g^+_i(u)\lambda^+_i(u)(g^-_i(u))\right)\right|_{u^m = \xi_{i,m+n}}\\[1.1 ex]
&=
c_i(\sigma_i^0)\left.\left(g^+_i(u)\lambda^+_i(u)(g^-_i(u))\right)\right|_{u^m = \xi_{i,m}}
\end{split}
\end{equation*}
where $\sigma_i^0:Y^0\to Y^0$ is the algebra homomorphism defined by $
\sigma_i^0(\xi_{j,m})=\xi_{j,m+\delta_{ij}}$. Since both $\fh$ and $\fg$ satisfy
$(B)$ with $k=0$, this yields
$$\Phi^0\left(\frac{e^{\frac{\hbar d_i}{2}H_{i,0}}-e^{-\frac{\hbar d_i}{2}H_{i,0}}}{q_i - q_i^{-1}}\right)=
c_i(\sigma_i^0)\Phi^0\left(\frac{e^{\frac{\hbar d_i}{2}H_{i,0}}-e^{-\frac{\hbar d_i}{2}H_{i,0}}}{q_i - q_i^{-1}}\right)$$
An inductive argument using the $\C[\hbar]$--linear $\N$--grading on $Y^0
$ given by $\deg(\xi_{j,m})=m$ and $\deg(\hbar)=0$ then shows that
$c_i^{(0)}=1$ and $c_i^{(n)}=0$ for any $n\geq 1$, so that $\fr$ satisfies (B$_0$).

To prove our claim, set
$$c_i(u)=r^+_i(u)\lambda^+_i(u)(r^-_i(u)) = r_i^-(u)\lambda^-_i(u)(r_i^+(u))$$
so that $r^-_i(u)=c_i(u)\lambda^-_i(u)(r_i^+(u))^{-1}$. By $(A_0)$, the following
holds for every $i,j\in\bfI$
$$r_i^+(u)\lambda^+_i(u)\left(c_j(v)\lambda^-_j(v)(r_j^+(v))^{-1}\right) =
c_j(v)\lambda^-_j(v)(r_j^+(v)^{-1})\lambda^-_j(v)(r_i^+(u))$$
Since $\lambda^+_i(u)$ and $\lambda^-_j(v)$ commute, we get
\[\begin{split}
\lambda^+_i(u)(c_j(v))\left(r_i^+(u)\lambda^-_j(v)(r^+_j(v))\right)
&=
c_j(v)\lambda^-_j(v)\left(r_i^+(u)\lambda^+_i(u)(r_j^+(v))\right)\\
&=
c_j(v)\lambda^-_j(v)\left(r_j^+(v)\lambda^+_j(v)(r_i^+(u))\right)\\
&=
c_j(v)\left(r_i^+(u)\lambda^-_j(v)(r_j^+(v))\right)
\end{split}\]
where the second equality uses (C$_0^+$) and the third one
$\lambda^-_j(v)\lambda^+_j(v) = 1$. We have therefore proved that
$$\lambda^+_i(u)(c_j(v))=c_j(v) \quad \mbox{for every }  i,j\in\bfI$$
By definition of the operators $\lambda_i^{\pm}$, this implies that
the coefficients of $c_j(v)$ lie in the centre of $Y_{\hbar}(\Lg)$,
which is trivial.
\end{pf}

\subsection{}

The uniqueness of homomorphisms of geometric type relies on the
following

\begin{prop}\label{pr:surjectivity}
Let $\{r_i^+(u)\}_{i\in\bfI}\subset 1+\wh{Y^0[u]}_+$ be a collection
of invertible elements satisfying condition $(C_0^+)$ of Lemma \ref{gauge-1}.
Then, there exists an element $\xi\in 1+\wh{Y^0}_+$ such that, for any
$i\in\bfI$
$$r_i^+(u)=\xi\cdot\lambda_i^+(u)(\xi)^{-1}$$
Moreover, if $\zeta \in \wh{Y^0}^{\times}$ is any element such that
$r_i^+(u) = \zeta\cdot \lambda_i^+(u)(\zeta)^{-1}$, then $\zeta = c\xi$
for some $c\in\C[[\hbar]]^{\times}$.
\end{prop}

\noindent
The proof of Proposition \ref{pr:surjectivity} is given in \S \ref{ss:start surj}--\S
\ref{ss:end surj}.

\subsection{}\label{ss:start surj}

We begin by linearising the problem. Set
$$\ol{r}_i(u)=\log(r_i(u))\in\wh{Y^0[u]}_+$$
By condition $(C_0^+)$, the following holds for any $i,j\in\bfI$
\begin{equation}\label{eq:c0}
(\lambda_i^+(u)-1)(\ol{r}_j(v)) = (\lambda_j^+(v)-1)(\ol{r}_i(u))
\end{equation}
and we need to show that
\begin{equation}\label{eq:goal}
\ol{r}_i(u)=(\lambda_i^+(u)-1)\eta
\end{equation}
for some $\eta\in\wh{Y^0}_+$.

\subsection{Rank $1$ case}

We assume first that $|\bfI|=1$ and accordingly drop the subscript $i$ from our
formulae. We shall prove \eqref{eq:goal} by working with an adapted system
of generators of $Y^0$.

Recall that, by Proposition \ref{raising-lowering-operators},
$$(\lambda^+(u)-1) B(v)=-\frac{e^{\hbar v} - e^{-\hbar v}}{v}e^{uv}$$
Define $\displaystyle{B'(v)=\sum_{k\geq 0}\frac{v^k}{k!}\,t_k'}$ by equating the
coefficients of $v$ in
$$B'(v)=
-\frac{v}{e^{\hbar v}-e^{-\hbar v}} B(v)=
-\frac{\hbar v}{e^{\hbar v}-e^{-\hbar v}}\sum_{n\geq 0}  \frac{v^n}{n!}\,t_n$$
The elements $\{t_k'\}_{k\in\N}$ give another system of generators of $Y^0$
which are homogeneous, with $\deg(t_k')=k=\deg(t_k)$ for any $k\in\N$, and
satisfy
\begin{equation}\label{eq:lambda on t'}
\lambda^+(u)(t_k') = t_k' + u^k
\end{equation}

\subsection{}

Since the operator $\lambda^+(u):Y^0\to Y^0[u]$ is homogeneous with respect
to the $\N$--grading extending that on $Y^0$ by $\deg(u)=1$, it suffices to prove
\eqref{eq:goal} when $\ol{r}(u)$ is homogeneous of degree $n\in\N$. Moreover,
since $\lambda^+(u)$ is $\C[\hbar]$--linear and the formulae \eqref{eq:lambda on t'}
do not involve $\hbar$, we may further assume that the coefficients of $\ol{r}(u)$
lie in the $\C$--subalgebra $\ol{Y^0}\subset Y^0$ generated by the $\{t_k'\}$.

An element of $\ol{Y^0}[u]_n$ has the form

\begin{equation}\label{general-form-r-1}
\bar{r}(u) = \sum_{|\mu| \leq n} a_{\mu} t_{\mu}' u^{n-|\mu|}
\end{equation}
where $a_{\mu} \in \C[t_0']$ and, for a partition $\mu$ of length $l$, we define
$t_{\mu}' = t_{\mu_1}'\cdots t_{\mu_l}'$.
The proof of the existence of $\eta\in\ol{Y^0}_n$ such that $(\lambda^+(u)-1)
(\eta)=\ol{r}(u)$ proceeds in two steps:
\begin{enumerate}
\item\label{it:eta 1} Show that, modulo elements of the form $(\lambda^+(u)-1)(\eta)$,
\begin{equation}\label{general-form-r-2}
\bar{r}(u) = \sum_{|\mu|<n} a_{\mu} t_{\mu}'u^{n-|\mu|}
\end{equation}
where $a_{\mu}\in\C$ do not depend on $t_0'$.
\item\label{it:eta 2} Show that any $\ol{r}(u)$ of the form \eqref{general-form-r-2} is equal
to $(\lambda^+(u)-1)(\eta)$ for some $\eta\in\ol{Y^0}_n$.
\end{enumerate}

\subsection{Proof of \eqref{it:eta 1}}

For $\bar{r}(u)\in\ol{Y^0}_n$ of the form \eqref{general-form-r-1}, choose $b_{\nu}
\in\C[t_0']$ for every $\nu\vdash n$ such that
$$b_{\nu}(t_0'+1)-b_{\nu}(t_0') = a_{\nu}(t_0')$$
Then
$$\bar{r}(u)-(\lambda^+(u)-1) \left(\sum_{\nu\vdash n} b_{\nu} t_{\nu}'\right)
=
\sum_{|\mu|<n} a_{\mu}' t_{\mu}' u^{n-|\mu|}$$
for some $a_\mu'\in\C[t_0']$, so that we may assume that $\ol{r}(u)$ is of
the form \eqref{general-form-r-1} with $a_\mu=0$ for any $\mu\vdash n$.

Write now
\begin{multline*}
(\lambda^+(v)-1)\ol{r}(u)\\
=
\sum_{|\mu|<n}\left(a_\mu(t_0+1)-a_\mu(t_0)\right)t'_\mu u^{n-|\mu|}+
\sum_{|\mu|<n}
a_\mu(t_0+1)
\left(\sum_{\nu\subsetneq \mu}c(\nu,\mu)t_\nu' v^{|\mu|-|\nu|}\right)u^{n-|\mu|}
\end{multline*}
where $c(\nu,\mu)$ is the number of ways of obtaining $\nu$ by removing rows
from $\mu$. By \eqref{eq:c0}, the above expression is symmetric in $u$ and $v$.
Its value at $u=0$, which is $0$, must therefore equal its value at $v=0$, thus
leading to
$$\sum_{|\mu|<n}\left(a_\mu(t_0+1)-a_\mu(t_0)\right)t'_\mu u^{n-|\mu|}=0$$
which implies that $a_\mu\in\C$ for any $\mu$.

\subsection{Proof of \eqref{it:eta 2}}

Let $\bar{r}(u)$ be of the form \eqref{general-form-r-2}. For any $0\leq l\leq n$,
write
$$\ol{r}_l(u)=\sum_{\substack{|\mu|<n \\ l(\mu) = l}}a_{\mu}t_{\mu}' u^{n-|\mu|}$$
so that $\bar{r}(u) = \sum_l \bar{r}_l(u)$. We proceed by induction on the largest
positive integer $k$ such that $\bar{r}_k(u) \neq 0$. If $k=0$, then $\bar{r}(u)=cu
^n=(\lambda^+(u)-1)(ct_n')$.

Assume now that $k>0$ and let $D(u):\ol{Y^0}\rightarrow\ol{Y^0}[u]$ be
the differential operator $D(u)=\sum_{m\geq 1} u^m\partial_{t_m'}$. Since
$(\lambda^+(u)-1)(t_k') = u^k$ we get, for any partition $\mu$
$$(\lambda^+(u)-1)(t_{\mu}') = D(u)(t_{\mu}') +\text{terms of smaller length}$$
Thus, \eqref{eq:c0} implies that 
$$D(u)\bar{r}_k(v) = D(v)\bar{r}_k(u)$$
This cross--derivative condition implies the existence of $\eta\in\ol{Y^0}_n$ such
that $r_k(u)=D(u)\eta$. This implies that $\bar{r}(u)-(\lambda^+(u)-1)(\eta)$
has smaller $k$.

This completes the proof of the existence part of Proposition \ref{pr:surjectivity}
when $\Lg$ is of rank $1$.

\subsection{Arbitrary rank}

The argument for arbitrary $\Lg$ rests on the following

\begin{lem}
There exist generators $\{\varpi_{i,r}\}_{i\in\bfI,r\in \N}$ of $Y^0$
which are homogeneous, with $\deg(\varpi_{i,r})=r$ and such that
$$\lambda^\pm_i(u)\varpi_{j,r} = \varpi_{j,r}\pm\delta_{i,j} u^r$$
\end{lem}
\begin{pf}
By Proposition \ref{raising-lowering-operators}, the generating
series $B_j(v)=\hbar\sum_{r\geq 0}t_{j,r}v^r/r!$ satisfy
$$(\lambda_i^{\pm}(u)-1)\hbar^{-1}B_j(v)=\mp Q_{ij}(v)e^{uv}$$
where $Q_{ij}(v)=2\sinh(\hbar d_ia_{ij}v/2)/\hbar v$. Since $Q_{ij}
=d_ia_{ij}$ mod $\hbar$, the matrix $Q=(Q_{ij})$ is invertible. Set
$B'_i(v)=-
\hbar^{-1}\sum_jQ^{-1}_{ij}B_j(v)$. Then $(\lambda_i^{\pm} (u)-1)
B'_j(v)=\pm\delta_{ij} e^{uv}$ which, in terms of the expansion
$B'_i(v)=\sum\varpi_{i,r} v^r/r!$ yields the required transformation
property.

Since $\deg(v)=1$, the stated homogeneity of the $\varpi_{i,r}$ is
equivalent to $\Ad(\zeta)(B_i'(v))=B_i'(\zeta^2v)$ where $\Ad(\zeta)$
denotes the action of $\zeta\in\C^\times$ on $\Yhg[[v]]$ corresponding
to the $\N$--grading. This in turn follows from the fact that $\Ad(\zeta)
(\hbar^{-1}B_j(v))=B_j(\zeta^2v)$ and $\Ad(\zeta)Q(v)=Q(\zeta^2 v)$.
\end{pf}

Using the generators $\varpi_{i,r}$, the proof of the existence part of
Proposition \ref{pr:surjectivity} in higher rank follows the same argument
as the one used for proving the sufficiency of the cross derivative condition
(here the existence of a primitive for any $i\in\bfI$ is guaranteed by the
rank $1$ case).

\subsection{Uniqueness of $\xi$}\label{ss:end surj}

Let $\zeta\in \wh{Y^0}^{\times}$ be an element such that $r_i^+(u)=\zeta
\cdot\lambda_i^+(u)(\zeta)^{-1}$ for each $i\in\bfI$. Then
\begin{align*}
\lambda_i^+(u)(\zeta\xi^{-1}) &= \lambda_i^+(u)(\zeta)\lambda_i^+(u)(\xi)^{-1} \\
&= r_i^+(u)^{-1} \zeta r_i^+(u) \xi^{-1} \\
&= \zeta\xi^{-1}
\end{align*}
By Proposition \ref{raising-lowering-operators} (6), we get that $\lambda_i^{\pm}
(u)(\zeta\xi^{-1}) = \zeta\xi^{-1}$ for each $i\in\bfI$. By definition of the operators
$\lambda_i^{\pm}(u)$, this implies that the coefficients of $\zeta\xi^{-1}$ lie in the
centre of $\Yhg$, which is trivial. This completes the proof of the last assertion of
Proposition \ref{pr:surjectivity}.

\subsection{Torus action}

The adjoint action of $\h$ on $\Yhg$ exponentiates to one of the algebraic torus
$H=\Hom_{\Z}(Q,\C^\times)$ where $Q\subset\h^*$ is the root lattice. This action
preserves homomorphisms of geometric type and acts on the corresponding
formal power series by $\zeta\cdot\{g_i^\pm(u)\}=\{\zeta_i^{\pm 1}g_i^\pm(u)\}$
where $H\ni\zeta\to\zeta_i=\zeta(\alpha_i)$ is the $i$th coordinate function on
$H$.

\subsection{Uniqueness of homomorphisms of geometric type}

\begin{thm}\label{second-main-theorem}
Let $\Phi, \Phi':\hloop\to\wh{\Yhg}$ be two homomorphisms of geometric
type. Then, there exists $\zeta\in H$ and  $\xi\in 1+\wh{Y^0}_+$ such that
$$\Phi'=\Ad(\xi)\circ (\zeta\cdot\Phi)$$
Moreover, $\zeta$ is unique and $\xi$ is unique up to multiplication by
$c\in \C[[\hbar]]^{\times}$.
\end{thm}
\begin{pf}
Let $\{g_i^{\pm}(u)\}$, $\{h_i^{\pm}(u)\}\subset \wh{Y^0[[u]]}$ be elements
of $\G$ corresponding to $\Phi$ and $\Phi'$ respectively. By Lemma
\ref{le:constant term}, we may use the action of $H$ to assume that
$g_i^{\pm}(u) = h_i^{\pm}(u)\mod \wh{Y^0[u]}_+$. By Lemma \ref
{gauge-1}, the elements $r_i^{\pm}(u)=h_i^\pm(u)\cdot g_i^\pm(u)^{-1}
\in 1+\wh{Y^0[u]}_+$ satisfy conditions (A$_0$)--(C$_0^\pm$). By
Proposition \ref{pr:surjectivity}, we may find an element $\xi\in 1+\wh
{Y^0}_+$ such that $r_i^+(u) = \xi\cdot\lambda^+_i(u)(\xi^{-1})$. It
follows that for any $i\in\bfI$
\[\begin{split}
\Phi'(E_{i,0})
&= h_i^+(\sigma_i^+)x^+_{i,0}\\
&= r_i^+(\sigma_i^+)g_i^+(\sigma_i^+)x^+_{i,0} \\
&= \xi\lambda^+_i(\sigma_i^+)(\xi^{-1}) g_i^+(\sigma_i^+)x^+_{i,0} \\
&= \xi g_i^+(\sigma_i^+) x^+_{i,0} \xi^{-1} 
\end{split}\]
Moreover, for any $r\in\Z$, 
$$\Phi'(E_{i,r})=
e^{r\sigma_i^+}\Phi'(E_{i,0})=
e^{r\sigma_i^+}\xi\Phi(E_{i,0})\xi^{-1}=
\xi\Phi(E_{i,r})\xi^{-1}$$
By $(B_0)$, $r_i^-(u)=\lambda^-_i(u)(r_i^+(u)^{-1})=\xi\lambda^-_i(u)(\xi^{-1})$
and it follows similarly that $\Phi'(F_{i,r})=\xi\Phi(F_{i,r})\xi^{-1}$ for any $i\in\bfI$
and $r\in\Z$. Since $\Phi$ and $\Phi'$ coincide on $U^0$ and $\Ad(\xi)(\eta)=
\eta$ for any $\eta\in Y^0$ it follows that $\Phi'=\Ad(\xi)\circ\Phi$. The last
assertion of Proposition \ref{pr:surjectivity} implies the uniqueness of $\xi$
up to multiplication by an element of $\C[[\hbar]]^{\times}$.
\end{pf}

\section{Isomorphisms of geometric type}\label{section-degeneration}

We prove in this section that any homomorphism of geometric type
$\Phi:\hloop\to\wh{\Yhg}$ extends to an isomorphism of completed
algebras and induces Drinfeld's degeneration of $\hloop$ to $\Yhg$.

\subsection{Classical limit}\label{classical-limit}

The specialisations of the quantum loop algebra $\hloop$ and Yangian
$\Yhg$ at $\hbar=0$ are the enveloping algebras $U(\Lg[z,z^{-1}])$ and
$U(\Lg[s])$ respectively. Specifically, if $\{e_i,f_i,h_i\}_{i\in\bfI}$ are the
generators of $\g$ given in Section \ref{ss: KMA}, the assignments
\[e_i\otimes z^k\rightarrow E_{i,k},\qquad
f_i\otimes z^k\rightarrow F_{i,k},\qquad
h_i\otimes z^r\rightarrow H_{i,r}\]
and
\[e_i\otimes s^r\rightarrow\frac{1}{\sqrt{d_i}}\,x^+_{i,r},\qquad
f_i\otimes s^r\rightarrow\frac{1}{\sqrt{d_i}}\,x^-_{i,r},\qquad
h_i\otimes s^r\rightarrow\frac{1}{d_i}\,\xi_{i,r}\]
extend respectively to isomorphisms
$$\cloopvar\isom\hloop/\hbar\hloop\aand\current\isom\Yhg/\hbar\Yhg$$

\begin{prop}\label{classical-limit-phi}
Let $\Phi:\hloop\rightarrow\widehat{\Yhg}$ be the homomorphism given
by Theorem \ref{third-main-theorem}. Then, the specialisation of $\Phi$
at $\hbar=0$ is the homomorphism
$$\exp^*:\cloopvar\longrightarrow U(\g[[s]])\subset\wh{\currentvar}$$
given on $\g[z,z^{-1}]$ by $\exp^*(X\otimes z^k)=X\otimes e^{ks}$.
\end{prop}
\begin{pf}
Since $\Phi(H_{i,0})=d_i^{-1}t_{i,0}$ and, for $r\neq 0$,
\[\Phi(H_{i,r}) = \frac{\hbar}{q_i - q_i^{-1}} \sum_{k\geq 0} t_{i,k} \frac{r^k}{k!}\]
setting $\hbar=0$ yields $\left.\Phi\right|_{\hbar=0}(h_i\otimes z^0)=h_i\otimes
s^0$, and
\[\left.\Phi\right|_{\hbar=0}(h_i\otimes z^r)=
\frac{1}{d_i}\sum_{k\geq 0} d_ih_i\otimes \frac{s^k r^k}{k!} = h_i\otimes e^{rs}\]

Further, since $g_i^+(u)=\frac{1}{\sqrt{d_i}}\mod\hbar$ by \eqref
{eq:intro solution}, we get
\[\left.\Phi\right|_{\hbar=0}(e_i\otimes z^r)=
\frac{1}{\sqrt{d_i}}e^{r\sigma_i^+} \sqrt{d_i}e_i\otimes s^0=
\sum_{k\geq 0} e_i\otimes\frac{s^k r^k}{k!}=
e_i\otimes e^{rs}\]
where we used the fact that, in the classical limit, the operator $\sigma_i^+$
corresponds to multiplication by $s$. Similarly, $\left.\Phi\right|_{\hbar=0}
(f_i\otimes z^r)=f_i\otimes e^{rs}$.
\end{pf}

\subsection{}\label{ss:geo iso}

Let $\J\subset U_{\hbar}L\Lg$ be the kernel of the composition
\[U_{\hbar}L\Lg\xrightarrow{\hbar\to 0}\cloop\xrightarrow{z\to 1}U\Lg\]
and let
\[\wh{\hloop}=\lim_{\longleftarrow}\hloop/\J^n\]
be the completion of $\hloop$ \wrt the ideal $\J$.

\begin{thm}\label{th:geo iso}
Let $\Phi:\hloop\to\wh{Y_\hbar\Lg}$ be a homomorphism of geometric type.
Then,
\begin{enumerate}
\item $\Phi$ maps $\J$ to the ideal $\wh{\Yhg}_+=\prod_{n\geq 1}\Yhg_n$.
\item The corresponding homomorphism
$$\wh{\Phi}:\wh{\hloop}\rightarrow\wh{\Yhg}$$
is an isomorphism.
\end{enumerate}
\end{thm}
\begin{pf}
(1) Note first that $\J$ is generated by $\hbar\hloop$ and the elements
$\{H_{i,r}-H_{i,s},E_{i,r}-E_{i,s},F_{i,r}-F_{i,s}\}_{i\in\bfI,r,s\in\Z}$ since its
image in $\cloopvar$ is generated by the classes of these elements. Note
next that, for $r,s\neq 0$
\begin{align*}
\Phi(H_{i,r}-H_{i,s})
&=
\frac{\hbar}{q_i - q_i^{-1}}\sum_{k\geq 1}\frac{r^k - s^k}{k!} t_{i,k}\\
\intertext{while}
\Phi(H_{i,r}-H_{i,0})
&=
\frac{\hbar}{q_i - q_i^{-1}}\sum_{k\geq 1}\frac{r^k}{k!}t_{i,k}+
\left(\frac{\hbar}{q_i - q_i^{-1}}-d_i^{-1}\right)t_{i,0}
\end{align*}
which lies in $\prod_{n\geq 1}\Yhg_n$ since $\hbar/(q_i - q_i^{-1})=d_i
^{-1}\mod\hbar$. Finally, for $r,s\in\Z$,
$$\Phi(E_{i,r}-E_{i,s})=
(e^{r\sigma^+_i}-e^{s\sigma^+_i}) g_i^+(\sigma^+_i)e_{i,0}\in\J$$
and similarly $\Phi(F_{i,r}-F_{i,s})\in\J$.

(2) By Theorem \ref{second-main-theorem}, it suffices to prove this
for the explicit homomorphism given by Theorem \ref{third-main-theorem}.
The result then follows Proposition \ref{pr:complete flat} below and the fact
that, by Proposition \ref{classical-limit-phi}, the specialisation of $\wh{\Phi}$
at $\hbar=0$ is an isomorphism $\wh{\cloopvar}\to\wh{\currentvar}$.
\end{pf}

\subsection{}\label{ss:completion} 

Let $J\subset\cloopvar$ be the kernel of evaluation at $z=1$ and $\wh
{\cloop}$ the completion of $\cloop$ \wrt $J$.

\begin{prop}\label{pr:complete flat}\hfill
\begin{enumerate}
\item 
$\wh{\hloop}$ is a flat deformation of $\wh{\cloopvar}$.
\item 
$\wh{\Yhg}$ is a flat deformation of $\wh{\current}$ over $\C[[\hbar]]$.
\end{enumerate}
\end{prop}
\begin{pf}
(1) Set, for brevity $\U=\hloop$ and $U=\cloopvar$. We claim that $\wh
{\U}$ is a flat deformation of $\wh{\U}/\hbar\wh{\U}$, and that $\wh{\U}
/\hbar\wh{\U}\cong \wh{U}$.

To prove the first assertion it suffices to show, by \cite[Prop XVI.2.4]
{Kassel}, that $\wh{\U}$ is a separated, complete and torsion--free
$\C[[\hbar]]$--module. To show that it is separated, note that $\hbar
\in\J$, so that $\hbar^k\wh{\U}\subset\ds{\lim_{\stackrel{\longleftarrow}
{n>k}}}\J^k/\J^n$ and
\[\bigcap_{k\geq 0}\hbar^k\wh{\U}=\{0\}\]
To show completeness, note that
\[\wh{\U}/\hbar^k\wh{\U}=
\lim_{\stackrel{\longleftarrow}{n}}\,(\U/\J^n)/(\hbar^k\U/\hbar^k\U\cap\J^n)=
\lim_{\stackrel{\longleftarrow}{n}}
\left\{\begin{array}{ll}
\U/\J^n			&\text{if $n\leq k$}\\
\U/\hbar^k\U+\J^n	&\text{if $n>k$}\\
\end{array}\right.\]
from which it readily follows that the map
$$\wh{\U}\longrightarrow\lim_{\stackrel{\longleftarrow}{k}}\wh{\U}/\hbar^k\wh{\U}$$
is surjective.
Finally, to prove that $\wh{\U}$ is torsion--free, note that the kernel of
multiplication by $\hbar$ on $\U/\J^n$ is $\hbar^{-1}(\hbar\U\cap\J^n)/
\J^n$. We claim that $\hbar\U\cap\J^n=\hbar\J^{n-1}$, which implies
that the kernel of $\hbar$ on $\wh{\U}$ is $\lim_n\J^{n-1}/\J^n=\{0\}$.
To prove the claim, use the flatness of $\U$ to identify it with the $\C
[[\hbar]]$--module $U[[\hbar]]$, so that $\J=J\oplus\hbar U[[\hbar]]$.
Let $a_1,\ldots,a_n\in\J$ and write $a_i=a_i^0+\hbar\ol{a}_i$, where
$a_i^0\in J$ and $\ol{a}_i\in U[[\hbar]]$. Then
$$a_1\cdots a_n=
a_1^0a_2\cdots a_n\mod\hbar\J^{n-1}=
\cdots=
a_1^0\cdots a_n^0\mod\hbar\J^{n-1}$$
from which the claim follows.

The fact that $\wh{\U}/\hbar\wh{\U}\cong \wh{U}$ follows by taking 
limits in the sequence
\[
0\rightarrow \hbar\lp\U/\J^{n-1}\rp
\rightarrow \U/\J^n\rightarrow U/J^n\rightarrow 0
\]
The latter is exact since, under the natural surjection $\U\to U$, the
ideal $\J^n$ is mapped to $J^n$ with kernel $\hbar\U\cap \J^n=\hbar
\J^{n-1}$.\\

(2) Since $\wh{\Yhg}$ is the completion of $\Yhg$ \wrt the ideal $\Yhg_+$
of elements of positive degree, it follows as in (1) that it is a separated and
complete $\C[[\hbar]]$--module. The lack of torsion of $\Yhg$ implies that
$\hbar\Yhg\cap\Yhg_+^n=\hbar\Yhg_+^{n-1}$ and therefore that $\wh{\Yhg}$
is torsion--free. Thus, $\wh{\Yhg}$ is a flat deformation of
$$\wh{\Yhg}/\hbar\wh{\Yhg}\cong\wh{\Yhg/\hbar\Yhg}\cong\wh{\current}$$
as claimed.
\end{pf}

\subsection{Drinfeld's degeneration}\label{phi-degeneration-principle}

Consider the descending filtration
\begin{equation}\label{eq:filter}
U_\hbar(L\Lg)=\J^0\supset\J\supset\J^2\supset\cdots
\end{equation}
defined by the powers of $\J$ and let $\gr_\J(\hloop)=\bigoplus_{n\geq 0}
\J^n/\J^{n+1}$ be its associated graded. 

\begin{thm}[\cite{drinfeld-quantum-groups,guay-degeneration}]\label{th:degeneration-principle}
Let $\{d_i^{\pm}\}_{i\in\bfI}\subset\C^\times$ be such that $d_i^+d_i^-=d_i$.
Then, the following assignment extends uniquely to an isomorphism of graded
algebras $\Yhg\isom\gr_\J(\hloop)$
\begin{gather*}
\xi_{i,0} 	\longmapsto d_iH_{i,0}\in\hloop/\J\\[1.1ex]
x^+_{i,0}	\longmapsto d_i^+E_{i,0}\in\hloop/\J,\qquad
x^-_{i,0}	\longmapsto d_i^-F_{i,0}\in\hloop/\J\\[1.1ex]
x^+_{i,1}	\longmapsto d_i^+\,(E_{i,1}-E_{i,0})\in\J/\J^2,\qquad
x^-_{i,1}	\longmapsto d_i^-\,(F_{i,1}-F_{i,0})\in\J/\J^2
\end{gather*}
\end{thm}

\begin{rem} The fact that $\hloop$ degenerates to $\Yhg$ is stated, without
proof, in \cite[\S 6]{drinfeld-quantum-groups}. The formulae above and the
proof that they define an isomorphism $\Yhg\cong\gr_\J(\Uhg)$ are given
in \cite{guay-degeneration}.
\end{rem}

\subsection{Relation to Drinfeld's degeneration}

By Theorem \ref{th:geo iso}, a homomorphism of geometric type $\Phi$
induces a homomorphism
$$\gr(\Phi):\gr_\J(\hloop)\longrightarrow\Yhg=\gr_{\wh{\Yhg}+}{\wh{\Yhg}}$$
Let $\{g_i^{\pm}(v)\}\subset \wh{Y^0[v]}^{\times}$ be the elements defining
$\Phi$. By Lemma \ref{le:constant term},
\begin{equation}\label{eq:gi mod +}
g_i^{\pm}(v) = \frac{1}{d_i^{\pm}} \mod \wh{Y^0[v]}_+
\end{equation}
for some $d_i^{\pm}\in\C^\times$ such that $d_i^+d_i^- = d_i$.\\

\begin{prop}\label{pr:degeneration}
$\gr(\Phi)$ is the inverse of the degeneration isomorphism $\imath:\Yhg
\isom\hloop$ given by Theorem \ref{th:degeneration-principle}.
\end{prop}
\begin{pf}
It suffices to verify the claim on the generators $\{\xi_{i,0},x^\pm_{i,0},x^
\pm_{i,1}\}_{i\in\bfI}$ of $\Yhg$. Now,
\[\gr(\Phi)\circ\imath(\xi_{i,0})
= \gr(\Phi)(d_iH_{i,0})
= \xi_{i,0}\]
and
\begin{align*}\gr(\Phi)\circ\imath(x^+_{i,0})
&= d_i^+\,\Phi(E_{i,0})\mod\wh{\Yhg}_+\\
&= d_i^+\,g_i^+(\sigma_i^+) x^+_{i,0}\mod\wh{\Yhg}_+\\
&= x^+_{i,0}\end{align*}
by \eqref{eq:gi mod +}. Moreover,
\begin{align*}
\gr(\Phi)\circ\imath(x^+_{i,1})
&= d_i^+\,\Phi(E_{i,1} - E_{i,0})\mod\wh{\Yhg}_{\geq 2}\\
&= d_i^+(e^{\sigma^+_i} - 1)g_i^+(\sigma^+_i)x^+_{i,0}\mod\wh{\Yhg}_{\geq 2}\\
&= x^+_{i,1}
\end{align*}
And similarly $\gr(\Phi)\circ\imath(x^-_{i,r})=x^-_{i,r}$ for $r=0,1$.
\end{pf}

\section{Geometric solution for $\Lgl_n$}\label{section-gln}

In this section, we construct a homomorphism of geometric type
for $\Lgl_n$, and show that it intertwines the geometric realisations
of the corresponding loop algebra and Yangian constructed by
Ginzburg--Vasserot \cite{ginzburg-vasserot,vasserot-qaffine}.

\subsection{The quantum loop algebra \cite{ding-frenkel}}\label{section-definitions-gln}

Throughout this section, we fix $n\geq 2$ and mostly follow the notation of
\cite{vasserot-qaffine}. Set $\bfI=\{1,\ldots,n-1\}$ and $\bfJ=\{1,\ldots,n\}$.
Then, $\hloopgl{n}$ is topologically generated over $\C[[\hbar]]$ by elements
$\{E_{i,r},F_{i,r},D_{j,r}\}_{i\in\bfI,j\in\bfJ,r\in\Z}$. To describe the relations,
introduce the formal power series
\[E_i(z) = \sum_{r\in\Z} E_{i,r}z^{-r}\qquad\qquad
F_i(z) = \sum_{r\in\Z} F_{i,r}z^{-r}\]
and
\[\Theta_j^\pm(z)=
\sum_{s\geq 0}\Theta^\pm_{j,\pm s}z^{\mp s}=
\exp\left(\pm\frac{\hbar D_{j,0}}{2}\right)
\exp\left(\pm(q-q^{-1})\sum_{s\geq 1}D_{j,\pm s}z^{\mp s}\right)\]
The relations are
\begin{itemize}
\item[(QL1-gl)] For any $j,j'\in\bfJ$ and $r,s\in\Z$,
\[[D_{j,r},D_{j',s}]=0\]
\item[(QL2-gl)] For any $i\in\bfI$ and $j\in\bfJ$,
\begin{align*}
\Theta_j^\pm(z)E_i(w)\Theta_j^\pm(z)^{-1}
&=
\vartheta_{c_{ji}}(q^{c_{ji}}z/w)E_i(w)\\[1.1 ex]
\Theta_j^\pm(z)F_i(w)\Theta_j^\pm(z)^{-1}
&=
\vartheta_{c_{ji}}(q^{c_{ji}}z/w)^{-1}F_i(w)
\end{align*}
where $c_{ji}=-\delta_{ji}+\delta_{j\,i+1}$, $\vartheta_m(\zeta)=\ds{\frac{q^m\zeta-1}
{\zeta-q^m}}$, and the \rhs is expanded in powers of $z^{\mp 1}$.\footnote{note
that the expansions in $z^ {\pm 1}$ are related by the symmetry $\vartheta_m(\zeta
^{-1})=\vartheta_{-m}(\zeta)$.} 
%
\item[(QL3-gl)] For any $i,i'\in\bfI$,
\begin{align*}
E_i(z)E_{i'}(w)&=\vartheta_{a_{ii'}}(q^{i-i'}z/w)E_{i'}(w)E_i(z)\\
F_i(z)F_{i'}(w)&=\vartheta_{a_{ii'}}(q^{i-i'}z/w)^{-1}F_{i'}(w)F_i(z)
\end{align*}
where $a_{ii'}=2\delta_{ii'}-\delta_{|i-i'|,1}$ are the entries of the Cartan matrix
of $\Lsl_n$ and the equalities are understood as holding after both side have
been multiplied by the denominator of the function $\vartheta_m$.
\item[(QL4-gl)] For any $i,i'\in\bfI$,
\[(q-q^{-1})[E_i(z), F_{i'}(w)] =
\delta_{i,i'} \delta(z/w)
\left(\frac{\Theta^+_{i+1}(z)}{\Theta^+_i(z)} - \frac{\Theta^-_{i+1}(z)}{\Theta^-_i(z)}\right)
\]
where $\delta(\zeta)=\sum_{r\in\Z}\zeta^r$ is the formal delta function.
\item[(QL5-gl)] For any $i,i'\in\bfI$ such that $|i-i'|=1$,
\begin{multline*}
E_i(z_1)E_i(z_2)E_{i'}(w) - (q+q^{-1})E_i(z_1)E_{i'}(w)E_i(z_2) + \\
E_{i'}(w)E_i(z_1)E_i(z_2) + (z_1\leftrightarrow z_2) = 0
\end{multline*}
\begin{multline*}
F_i(z_1)F_i(z_2)F_{i'}(w) - (q+q^{-1})F_i(z_1)F_{i'}(w)F_i(z_2) + \\
F_{i'}(w)F_i(z_1)F_i(z_2) + (z_1\leftrightarrow z_2) = 0
\end{multline*}
For any $i,i'\in\bfI$ such that $|i-i'|\geq 2$,
\begin{gather*}
E_i(z)E_{i'}(w)=E_{i'}(w)E_i(z)\\
F_i(z)F_{i'}(w)=F_{i'}(w)F_i(z)
\end{gather*}
\end{itemize}
In terms of the generators $\{E_{i,r},F_{i,r},D_{j,r}\}$, the relations (QL1-gl)--(QL4-gl) read
\begin{gather*}
[D_{j,0},E_{i,k}]=c_{ji}E_{i,k}\qquad [D_{j,r},E_{i,k}]=q^{-c_{ji}r}\frac{[c_{ji}r]}{r}E_{i,k+r}\\
[D_{j,0},F_{i,k}]=-c_{ji}F_{i,k}\qquad [D_{j,r},F_{i,k}]=-q^{-c_{ji}r}\frac{[c_{ji}r]}{r}F_{i,k+r}\\
q^iE_{i,k+1}E_{i',l}-q^{a_{ii'}+i'}E_{i,k}E_{i',l+1}=q^{a_{ii'}+i}E_{i',l}E_{i,k+1}-q^{i'}E_{i',l+1}E_{i,k}\\
q^{a_{ii'}+i}F_{i,k+1}F_{i',l}-q^{i'}F_{i,k}F_{i',l+1}=q^iF_{i',l}F_{i,k+1}-q^{a_{ij}+i'}F_{i',l+1}F_{i,k}\\
[E_{i,k},F_{i',l}]=\delta_{ii'}\frac{P^+_{i,k+l}-P^-_{i,k+l}}{q-q^{-1}}
\end{gather*}
where
$P_i^\pm(z)=\sum_{s\geq 0}P^\pm_{i,\pm s}z^{\mp s}=\Theta_{i+1}^\pm(z)/\Theta_i^\pm(z)$

We denote by $U^0\subset\hloopgl{n}$ the commutative subalgebra generated by
the elements $D_{j,r}$.

\subsection{The Yangian $\Yhgl{n}$}

The following definition can be found in \cite[\S 3.1]{molev-yangian}.
$\Yhgl{n}$ is the algebra over $\C[\hbar]$ generated by elements
$\{e_{i,r},f_{i,r},\theta_{j,r}\}_{i\in\bfI,j\in\bfJ,r\in\N}$ subject to the
following relations\footnote{our conventions are adapted to \cite
{ginzburg-vasserot,vasserot-qaffine}. They differ from those of
\cite{molev-yangian} by the permutation $e_{i,r}\leftrightarrow
f_{i,r}$ and the relabelling $\theta_{i,r}\leftrightarrow h_{i,r}$.}
\begin{itemize}
\item[(Y1-gl)] For any $j,j'\in\bfJ$ and $r,s\in\N$,
\[[\theta_{j,r},\theta_{j',s}]=0\]
\item[(Y2-gl)] For any $j\in\bfJ$ and $i\in\bfI$,
\begin{align*}
[\theta_{j,0},e_{i,s}] &= \phantom{-}c_{ji}e_{i,s} \\ 
[\theta_{j,0},f_{i,s}] &= -c_{ji}f_{i,s} \\ 
\end{align*}
\begin{align*}
[\theta_{j,r+1},e_{i,s}]-[\theta_{j,r},e_{i,s+1}] &= \phantom{-}\hbar c_{ji}e_{i,s}\theta_{j,r}\\
[\theta_{j,r+1},f_{i,s}]-[\theta_{j,r},f_{i,s+1}] &= -\hbar c_{ji}\theta_{j,r}f_{i,s}
\end{align*}
where $c_{ji}=-(\delta_{ji}-\delta_{j,i+1})$.
\item[(Y3-gl)] For any $i\in\bfI$,
\begin{align*}
[e_{i,r+1},e_{i,s}] - [e_{i,r}, e_{i,s+1}] &= \phantom{-}\hbar(e_{i,r}e_{i,s} + e_{i,s}e_{i,r})\\
[f_{i,r+1},f_{i,s}] - [f_{i,r}, f_{i,s+1}] &= -\hbar(f_{i,r}f_{i,s} + f_{i,s}f_{i,r})
\end{align*}
For any $i\in\bfI\setminus\{n-1\}$ and $r,s\in\N$,
\begin{align*}
[e_{i,r+1}, e_{i+1,s}] - [e_{i,r}, e_{i+1,s+1}]         &= -\hbar e_{i+1,s}e_{i,r}\\
[f_{i,r+1}, f_{i+1,s}] - [f_{i,r}, f_{i+1,s+1}]                 &= \phantom{-}\hbar f_{i,r}f_{i+1,s}
\end{align*}
\item[(Y4-gl)] For any $i,i'\in\bfI$,
\[
[e_{i,r},f_{i',s}] = \delta_{i,i'} p_{i,r+s}
\]
where $p_i(v) = 1 + \hbar\sum_{r\geq 0} p_rv^{-r-1} = \theta_{i+1}(v)\theta_i(v)^{-1}$.
\item[(Y5-gl)] For any $i,i'\in\bfI$ such that $|i-i'|=1$, and $r_1,r_2,s\in\N$,
\begin{align*}
[e_{i,r_1}, [e_{i,r_2}, e_{i',s}]] + [e_{i,r_2},[e_{i,r_1}, e_{i',s}]]  &= 0\\
[f_{i,r_1}, [f_{i,r_2}, f_{i',s}]] + [f_{i,r_2},[f_{i,r_1}, f_{i',s}]]         &= 0
\end{align*}
For $i,i'\in\bfI$ such that $|i-i'|>1$, and $r,s\in\N$
\[[e_{i,r}, e_{i',s}] = 0 = [f_{i,r}, f_{i',s}]\]
\end{itemize}
The Yangian $\Yhgl{n}$ is $\N$--graded by $\deg(e_{i,r})=\deg(f_{i,r})=\deg
(\theta_{j,r})=r$ and $\deg(\hbar)=1$.

\subsection{Shift homomorphisms}

Let $Y^0\subset\Yhgl{n}$ be the commutative subalgebra generated
by the elements $\{\theta_{j,r}\}$ and $Y^+,Y^-\subset\Yhgl{n}$ the
subalgebras generated by $Y^0$ and the elements $\{e_{i,r}\}$ (resp.
$\{f_{i,r}\}$), $i\in\bfI,r\in\N$.

For any $i\in\bfI$, define, as in Section \ref{ss:sigma_i}, a $Y^0$--linear homomorphism $\sigma_i^\pm$ of $Y^\pm$ by $e_{i',r}\to e_{i',r+\delta
_{ii'}}$ (resp. $f_{i',r}\to f_{i',r+\delta_{ii'}}$). The definition of $\sigma_i^
\pm$ relies on the PBW theorem for $\Yhgl{n}$, which is proved in
\cite{olshanskii}.

\subsection{Alternative generators for $Y^0$}

Define an alternative family of generators $\{d_{j,r}\}_{j\in\bfJ,r\in\N}$
of $Y^0$ by
\[d_j(u) = \hbar\sum_{r\geq 0} d_{j,r}u^{-r-1} = \log(\theta_j(u))\]
Set $\ds{B_j(v) = \hbar \sum_{r\geq 0} d_{j,r}\frac{v^r}{r!}\in Y^0[[v]]}$.
The following commutation relations are proved exactly as their
counterparts in Lemma \ref{commutation-with-log-like}.

\begin{lem}\label{commutation-with-log-like-gln}
The following holds for any $j\in\bfJ$ and $i\in\bfI$,
\begin{align*}
[B_j(v), e_{i,s}] 
&=
\phantom{-}
\frac{1-e^{-c_{ji}\hbar v}}{v}e^{\sigma_i^+ v}\,e_{i,s}\\[1.1 ex]
[B_j(v), f_{i,s}] 
&=
-
\frac{1-e^{-c_{ji}\hbar v}}{v}e^{\sigma_i^- v}\,f_{i,s}
\end{align*}
\end{lem}

\subsection{The operators $\mathbf{\lambda_i^\pm(v)}$}

The following result is analogous to, and proved in the same
way as, Proposition \ref{raising-lowering-operators}.

\begin{prop}\label{raising-lowering-operators-gln}
There are operators $\{\lambda^{\pm}_{i;s}\}_{i\in\bfI,s\in\N}$
on $Y^0$ such that the following holds.
\begin{enumerate}
\item For any $\xi\in Y^0$,
\begin{align*}
e_{i,r}\xi	&= \sum_{s\geq 0} \lambda^+_{i;s}(\xi) e_{i,r+s}\\
f_{i,r}\xi	&= \sum_{s\geq 0} \lambda^-_{i;s}(\xi) f_{i,r+s}
\end{align*}
\item The operator $\lambda^{\pm}_i(v):Y^0 \rightarrow Y^0[v]$ given by
\[\lambda^{\pm}_i(v)(\xi) = \sum_{s\in\N}\lambda^{\pm}_{i;s}(\xi)v^s\]
is an algebra homomorphisms of degree 0 \wrt the $\N$--grading on $Y^0[v]$
extending that on $Y^0$ by $\deg(v)=1$.
\item The operators $\lambda^\epsilon_i(v)$ and $\lambda^{\epsilon'}_{i'}(v')$ commute for
any $i,i'\in\bfI$ and $\epsilon,\epsilon'\in\{\pm\}$. Moreover,
\[\lambda^+_i(v)\lambda^-_i(v) = Id\]
\item For any $i\in\bfI$ and $j\in\bfJ$,
\begin{equation}\label{eq: raising-lowering-gln}
(\lambda^{\pm}_i(v_1)-1)B_j(v_2)=
\pm \frac{e^{-c_{ji}\hbar v_2}-1}{v_2} e^{v_1v_2}
\end{equation}
\end{enumerate}
\end{prop}

\subsection{}

Let $\{g_i^{\pm}(u)\}_{i\in\bfI}$ be a collection of elements in $\wh{Y^0
[u]}$. Define, as in Section \ref{definition-assignment}, an assignment
$\Phi:\{E_{i,r},F_{i,r},D_{j,r}\}
\to\wh{\Yhgl{n}}$ by
\begin{align*}
\Phi(D_{j,0})	&=\theta_{j,0}\\
\Phi(D_{j,r}) 	&= \frac{B_j(r)}{q-q^{-1}} \mbox{ for } r\not= 0\\
\Phi(E_{i,k})	&= e^{k\sigma_i^+}g_i^+(\sigma_i^+)\,e_{i,0}\\
\Phi(F_{i,k})	&= e^{k\sigma_i^-}g_i^-(\sigma_i^-)\,f_{i,0}
\end{align*}
and denote the restriction of $\Phi$ to $U^0$ by $\Phi^0$. For any $i\in
\bfI$, set
\begin{gather*}
\xi_i(u)=
1+\hbar\sum_{r\geq 0}\xi_{i,r}u^{-r-1}=\theta_{i+1}(u)\theta_i(u)^{-1}\in Y^0[[u^{-1}]]\\
P_i^\pm(z)
=\sum_{s\geq 0}P_{i,\pm s}^{\pm}z^{\mp s}=
\Theta^\pm_{i+1}(z)\Theta^\pm_i(z)^{-1}\in U^0[[z^{\mp 1}]]
\end{gather*}

\begin{thm}\label{first-main-theorem-gln}
The assignment $\Phi$ extends to an algebra homomorphism if and only
if the following conditions hold.
\begin{itemize}
\item[(A)] For any $i,i'\in\bfI$,
\[g_i^+(u)\lambda^+_i(u)(g_{i'}^-(v)) = g_{i'}^-(v)\lambda^-_{i'}(v)(g_i^+(u))\]
\item[(B)] For any $i\in \bfI$ and $k\in\Z$,
\[e^{ku}g_i^+(v)\lambda^+_i(v)(g_i^-(v))|_{v^m = \xi_{i,m}}
= \Phi^0\left(\frac{P_{i,k}^+ - P_{i,k}^-}{q-q^{-1}}\right)\]
\item[(C0)] For any $i,i'\in\bfI$ such that $|i-i'|>1$,
\[g^{\pm}_i(u)\lambda^{\pm}_i(u)(g^{\pm}_{i'}(v))=
g^{\pm}_{i'}(v)\lambda^{\pm}_{i'}(v)(g^{\pm}_i(u))\]
\item[(C1)] For any $i\in \bfI$
\[
g_i^{\pm}(u)\lambda^{\pm}_i(u)(g_i^{\pm}(v)) \frac{e^u - e^{v\pm\hbar}}{u-v\mp\hbar} =
g_i^{\pm}(v)\lambda_i^{\pm}(v)(g_i^{\pm}(u))\frac{e^v-e^{u\pm\hbar}}{v-u\mp\hbar}
\]
\item[(C2)] For any $i\in\bfI\setminus\{n-1\}$,
\[
g_i^{\pm}(u)\lambda^{\pm}_i(u)(g_{i+1}^{\pm}(v)) \left(\frac{e^u - e^v}{u-v}\right)^{\pm 1} = g_{i+1}^{\pm}(v)\lambda^{\pm}_{i+1}(v)(g_i^{\pm}(u)) \left(\frac{e^{u-\hbar/2}-e^{v+\hbar/2}}{u-v-\hbar}\right)^{\pm 1}
\]
\end{itemize}
\end{thm}

\noindent
The proof of Theorem \ref{first-main-theorem-gln} is given in \ref{ss:start gl}--\ref{ss:end gl}. It follows the same lines as that of Theorem \ref{first-main-theorem}, with the exception of the $q$--Serre relations (QL5-gl) which are proved directly.

\subsection{}\label{ss:start gl}

A proof similar to that of Lemmas \ref{le:AB} and \ref{le:C} yields the following
\begin{enumerate}
\item $\Phi$ is compatible with the relation (QL4-gl) if, and only if (A) and (B) hold.
\item $\Phi$ is compatible with the relation (QL3-gl) if, and only if (C0)--(C2) hold.
\end{enumerate}

\noindent
By virtue of condition (C0), $\Phi$ is compatible with the the $q$--Serre relations
(QL5-gl) whenever $|i-i'|>1$. We therefore need only consider the case $|i-i'|=1$.
We shall in fact restrict to $i'=i+1$ since the case $i'=i-1$ is dealt with similarly.

\subsection{}

The essential ingredient is the following analogue of Lemma \ref{app-lem: symmetric}.
We leave it to the reader to carry out the construction of the auxiliary algebra $\ol{Y}$
(see $\S$ \ref{app-sec: auxiliary}), the operators $\ol{\sigma}_{i,(1)},\ol{\sigma}_{i,(2)}$
and $\ol{\sigma}_{i'}$ on $\ol{Y}_{2\alpha_i+\alpha_{i'}}$ ($\S$ \ref{app-sec: shift}) and
the map $p_{ii'}:\ol{Y}_{2\alpha_i+\alpha_{i'}} \rightarrow \Yhgl{n}$.

\begin{lem}\label{le:shift-gln}
The kernel of $p_{ii'}$ is the $\C[\hbar]$--linear span of the following elements
\begin{enumerate}
\item For any $A(u_1,u_2,v)\in\ol{Y}^0[u_1,u_2,v]$
\begin{align*}
&A(\ol{\sigma}_{i,(1)},\ol{\sigma}_{i,(2)},\ol{\sigma}_{i'})
\left((\ol{\sigma}_{i,(2)}-\ol{\sigma}_{i'})\ol{e}_{i,0}^2\ol{e}_{i',0} - 
(\ol{\sigma}_{i,(2)}-\ol{\sigma}_{i'}-\hbar)\ol{e}_{i,0}\ol{e}_{i',0}\ol{e}_{i,0}\right)\\[1.1ex]
&A(\ol{\sigma}_{i,(1)}, \ol{\sigma}_{i,(2)},\ol{\sigma}_{i'})
\left((\ol{\sigma}_{i,(1)}-\ol{\sigma}_{i'})\ol{e}_{i,0}\ol{e}_{i',0}\ol{e}_{i,0}-
(\ol{\sigma}_{i,(1)}-\ol{\sigma}_{i'}-\hbar)\ol{e}_{i',0}\ol{e}_{i,0}^2\right)
\end{align*}
\item For any $B(u_1,u_2,v)=B(u_2,u_1,v)\in\ol{Y}^0[u_1,u_2,v]$
\begin{align*}
&B(\ol{\sigma}_{i,(1)},\ol{\sigma}_{i,(2)},\ol{\sigma}_{i'})
(\ol{\sigma}_{i,(1)}-\ol{\sigma}_{i,(2)}-\hbar)\ol{e}_{i,0}^2\ol{e}_{i',0}\\[1.1ex]
&B(\ol{\sigma}_{i,(1)},\ol{\sigma}_{i,(2)},\ol{\sigma}_{i'})
(\ol{\sigma}_{i,(1)}-\ol{\sigma}_{i,(2)}-\hbar)\ol{e}_{i',0}\ol{e}_{i,0}^2
\end{align*}
\item For any $B(u_1,u_2,v)=B(u_2,u_1,v)\in\ol{Y}^0[u_1,u_2,v]$
\[B(\ol{\sigma}_{i,(1)}, \ol{\sigma}_{i,(2)},\ol{\sigma}_{i'})
\left(\ol{e}_{i,0}^2\ol{e}_{i',0} - 2\ol{e}_{i,0}\ol{e}_{i',0}\ol{e}_{i,0}+\ol{e}_{i',0}\ol{e}_{i,0}^2\right)
\]
\end{enumerate}
\end{lem}

\begin{cor}\label{co:symmetric}
The kernel of $p_{ii'}$ is stable under the action of $A(\ol{\sigma}_{i,(1)},\ol{\sigma}_
{i,(2)}, \ol{\sigma}_{i'})$, for any $A(u_1,u_2,v)=A(u_2,u_1,v)\in\ol{Y}^0[u_1,u_2,v]$.
\end{cor}

\begin{rem}
In the next sections we use the following convention for notational convenience:
for each $\ol{X} \in \ol{Y}_{2\alpha_i+\alpha_{i'}}$ and $X = p_{ii'}(\ol{X})$, we set
\[
A(\sigma_{i,1},\sigma_{i,2},\sigma_{i'})(X)=
p_{i,i'}\left(A(\ol{\sigma}_{i,(1)}, \ol{\sigma}_{i,(2)},\ol{\sigma}_{i'})(\ol{X})\right)
\]
\end{rem}

\subsection{}

We shall only prove the $q$--Serre relations for the case of the $E$'s and consequently drop the superscript $+$. We need to show that the following
holds for any $k_1,k_2,l\in\Z$
\begin{multline*}
\Phi(E_{i,k_1})\Phi(E_{i,k_2})\Phi(E_{j,l}) -
(q+q^{-1})\Phi(E_{i,k_1})\Phi(E_{j,l})\Phi(E_{i,k_2})\\
+ \Phi(E_{j,l})\Phi(E_{i,k_1})\Phi(E_{i,k_2}) + (k_1 \leftrightarrow k_2) = 0
\end{multline*}
As in $\S$ \ref{app-sec: reduction}, an application of Corollary \ref{co:symmetric}
shows that this reduces to showing that
\[\Phi(E_{i,0})^2\Phi(E_{j,0})-(q+q^{-1})\Phi(E_{i,0})\Phi(E_{j,0})\Phi(E_{i,0})+\Phi(E_{j,0})\Phi(E_{i,0})^2=0\]

\subsection{}

With Corollary \ref{co:symmetric} in mind, we seek to factor a common
symmetric function out of each of the above summands. This is achieved
by the following result.

\begin{lem}\label{le:factor symm}
There exists $H(u_1,u_2,v)\in\ol{Y}^0[[u_1,u_2,v]]$ symmetric in $u_1
\leftrightarrow u_2$, such that
\begin{align*}
\Phi(E_{i,0})^2\Phi(E_{i',0})
&=
H(\sigma_{i,1},\sigma_{i,2},\sigma_{i'})
\PP_0(\sigma_{i,1},\sigma_{i,2},\sigma_{i'})e_{i,0}e_{i,0}e_{i',0}\\
\Phi(E_{i,0})\Phi(E_{i',0})\Phi(E_{i,0})
&=
H(\sigma_{i,1},\sigma_{i,2},\sigma_{i'})
\PP_1(\sigma_{i,1},\sigma_{i,2},\sigma_{i'})e_{i,0}e_{i',0}e_{i,0}\\
\Phi(E_{i',0})\Phi(E_{i,0})^2
&=
H(\sigma_{i,1},\sigma_{i,2},\sigma_{i'})
\PP_2(\sigma_{i,1},\sigma_{i,2},\sigma_{i'})e_{i',0}e_{i,0}e_{i,0}
\end{align*}
where $\PP_0,\PP_1,\PP_2\in\C[[u_1,u_2,v]]$ are given in terms of the function
\[P(x,y)=\frac{e^x-e^y}{x-y}\in\C[[x,y]]^{\Sym_2}\]
by
\begin{align*}
\PP_0 &= P(u_1+\hbar,u_2)P(u_1-\hbar/2,v+\hbar/2)P(u_2-\hbar/2, v+\hbar/2)\\
\PP_1 &= P(u_1+\hbar,u_2)P(u_1-\hbar/2,v+\hbar/2)P(u_2, v)\\
\PP_2 &= P(u_1+\hbar,u_2)P(u_1,v)P(u_2, v)
\end{align*}
\end{lem}

\begin{pf}

Define $G_{ab}(x,y)\in\ol{Y}^0[[x,y]]$ by
$\lambda_a(x)(g_b(y)) = g_b(y)G_{ab}(x,y)$. 
Then, in obvious notation,
\[\begin{split}
\Phi(E_{a,0})\Phi(E_{b,0})\Phi(E_{c,0})
&=
g_a(\sigma_a)e_{a,0}g_b(\sigma_b)e_{b,0}g_c(\sigma_c)e_{c,0}\\
&=
g_a(\sigma_a)\lambda_a(\sigma_a)(g_b(\sigma_b))
\lambda_a(\sigma_a)\circ\lambda_b(\sigma_b)(g_c(\sigma_c))e_{a,0}e_{b,0}e_{c,0}\\
&=
g_a(\sigma_a)g_b(\sigma_b)g_c(\sigma_c)
G_{ab}(\sigma_a,\sigma_b)
G_{ac}(\sigma_a,\sigma_c)
\lambda_a(\sigma_a)(G_{bc}(\sigma_b,\sigma_c))
e_{a,0}e_{b,0}e_{c,0}
\end{split}\]
We record for later use the symmetry in the interchange $a\leftrightarrow b$
of the term
\begin{equation}\label{eq:a<-->b}
\begin{split}
G_{ac}(\sigma_a,\sigma_c)
\lambda_a(\sigma_a)(G_{bc}(\sigma_b,\sigma_c))
&=
\lambda_a(\sigma_a)\circ\lambda_b(\sigma_b)(g_c(\sigma_c))/g_c(\sigma_c)\\
&=
G_{bc}(\sigma_b,\sigma_c)
\lambda_b(\sigma_b)(G_{ac}(\sigma_a,\sigma_c))
\end{split}
\end{equation}
where the second equality follows from the commutativity of $\lambda_a(\sigma_a)$
and $\lambda_b(\sigma_b)$.

Set now $F=g_i(\sigma_{i,1})g_i(\sigma_{i,2})g_{i'}(\sigma_{i'})\in\ol{Y}^0
[[\sigma_{i,1},\sigma_{i,2},\sigma_{i'}]]^{\Sym_2}$. Then, the above yields
\begin{align*}
\Phi(E_{i,0})^2\Phi(E_{i',0})
&= F\,
G_{ii}(\sigma_{i,1},\sigma_{i,2})
G_{ii'}(\sigma_{i,1},\sigma_{i'})
\lambda_i(\sigma_{i,1})(G_{ii'}(\sigma_{i,2},\sigma_{i'}))
e_{i,0}^2e_{i',0}\\
\Phi(E_{i,0})\Phi(E_{i',0})\Phi(E_{i,0}) 
&= F\,
G_{ii'}(\sigma_{i,1},\sigma_{i'})
G_{ii}(\sigma_{i,1}, \sigma_{i,2})
\lambda_i(\sigma_{i,1})(G_{i'i}(\sigma_{i'},\sigma_{i,2})) e_{i,0}e_{i',0}e_{i,0}\\
\Phi(E_{i',0})\Phi(E_{i,0})^2
&= F\,
G_{i'i}(\sigma_{i'}, \sigma_{i,1})
G_{i'i}(\sigma_{i'},\sigma_{i,2})
\lambda_{i'}(\sigma_{i'})(G_{ii}(\sigma_{i,1},\sigma_{i,2})) e_{i',0}e_{i,0}^2
\end{align*}

We claim that $G_{ii}(u_1,u_2)=\overline{G}_{ii}(u_1,u_2)P(u_1+\hbar, u_2)$,
where $\ol{G}$ is symmetric in $u_1,u_2$. Indeed, by condition (C1)
\[G_{ii}(u_1,u_2) P(u_1,u_2+\hbar) = G_{ii}(u_2,u_1)P(u_2,u_1+\hbar)\]
whence the result with $\overline{G}_{ii}(u_1,u_2)=G_{ii}(u_1,u_2)/P(u_1
+\hbar,u_2)$. It follows that
\begin{equation*}
\Phi(E_{i,0})^2\Phi(E_{i',0})
= \ol{H}(\sigma_{i,1},\sigma_{i,2},\sigma_{i'})
P(\sigma_{i,1}+\hbar,\sigma_{i,2})
e_{i,0}^2e_{i',0}
\end{equation*}
where
\[\ol{H}(u_1,u_2,v)=
g_i(u_1)g_i(u_2)g_{i'}(v)\ol{G}_{ii}(u_1,u_2)G_{ii'}(u_1,v)\lambda_i(u_1)(G_{ii'}(u_2,v))
\in{Y}^0[[u_1,u_2,v]]\]
is symmetric in $u_1,u_2$ by \eqref{eq:a<-->b}.

Next, assuming that $i'=i+1$, condition (C2) yields
\[G_{ii'}(u,v)P(u,v) = G_{i'i}(v,u)P(u-\hbar/2, v+\hbar/2)\]
so that
\begin{multline*}
\Phi(E_{i,0})\Phi(E_{i',0})\Phi(E_{i,0})=\\
\ol{H}(\sigma_{i,1},\sigma_{i,2},\sigma_{i'})
P(\sigma_{i,1}+\hbar,\sigma_{i,2})
\frac{P(\sigma_{i,2},\sigma_{i'})}{P(\sigma_{i,2}-\hbar/2, \sigma_{i'}+\hbar/2)}
e_{i,0}e_{i',0}e_{i,0}
\end{multline*}

Finally, using \eqref{eq:a<-->b} again, with $a=i,b=i',c=i$ and the previous calculation
yields
\begin{multline*}\label{eq:third}
\Phi(E_{i',0})\Phi(E_{i,0})^2=\\
\ol{H}(\sigma_{i,1},\sigma_{i,2},\sigma_{i'})
P(\sigma_{i,1}+\hbar,\sigma_{i,2})
\frac{P(\sigma_{i,1},\sigma_{i'})}{P(\sigma_{i,1}-\hbar/2, \sigma_{i'}+\hbar/2)}\frac{P(\sigma_{i,2},\sigma_{i'})}{P(\sigma_{i,2}-\hbar/2, \sigma_{i'}+\hbar/2)}e_{i',0}e_{i,0}^2
\end{multline*}
as claimed.
\end{pf}

\subsection{}\label{ss:end gl}

By Lemma \ref{le:factor symm} and Corollary \ref{co:symmetric}, we are reduced to proving the following
\begin{multline*}
\mathcal{S}^q = \PP_0(\sigma_{i,1},\sigma_{i,2},\sigma_{i'})e_{i,0}^2e_{i',0} - (q+q^{-1}) \PP_1(\sigma_{i,1},\sigma_{i,2},\sigma_{i'})e_{i,0}e_{i',0}e_{i,0}\\
+ \PP_2(\sigma_{i,1},\sigma_{i,2},\sigma_{i'})e_{i',0}e_{i,0}^2 = 0
\end{multline*}

{\bf Step 1.} We first observe that
\[P(u_1+\hbar, u_2)-\frac{1+e^{\hbar}}{2}P(u_1,u_2)
\in(u_1-u_2-\hbar)\C[[\hbar,u_1,u_2]]^{\Sym_2}\]
This allows us to use (2) of Lemma \ref{le:shift-gln} to obtain
\[\mathcal{S}^q = \PP_0'(\sigma_{i,1},\sigma_{i,2},\sigma_{i'})e_{i,0}^2e_{i',0} - 2 \PP_1'(\sigma_{i,1},\sigma_{i,2},\sigma_{i'})e_{i,0}e_{i',0}e_{i,0}
+ \PP_2'(\sigma_{i,1},\sigma_{i,2},\sigma_{i'})e_{i',0}e_{i,0}^2\]
where
\begin{align*}
\PP_0' &= e^{\hbar/2}P(u_1,u_2)P(u_1-\hbar/2,v+\hbar/2)P(u_2-\hbar/2, v+\hbar/2)\\
\PP_1' &= P(u_1+\hbar,u_2)P(u_1-\hbar/2,v+\hbar/2)P(u_2, v)=\PP_1\\
\PP_2' &= e^{\hbar/2}P(u_1,u_2)P(u_1,v)P(u_2, v)
\end{align*}

{\bf Step 2.} We use next (3) of Lemma \ref{le:shift-gln} with $B = \PP_2'$ to get
\[\mathcal{S}^q = (\PP_0'-\PP_2')e_{i,0}^2e_{i',0} - 2(\PP_1'-\PP_2')e_{i,0}e_{i',0}e_{i,0}\]
One can easily check that $\PP_1' - \PP_2'$ is divisible by $u_2-v-\hbar$, which together
with (1) of Lemma \ref{le:shift-gln}, with $A = \displaystyle 2 \frac{\PP_1'-\PP_2'}{u_2-v-\hbar}$,
yields
\[\mathcal{S}^q = \left(\PP_0'-\PP_2' - 2\frac{\PP_1'-\PP_2'}{u_2-v-\hbar}(u_2-v)\right)e_{i,0}^2e_{i',0}\]

{\bf Step 3.} Finally we can directly verify that the function $$\mathcal{F} = \displaystyle \PP_0'-\PP_2' - 2\frac{\PP_1'-\PP_2'}{u_2-v-\hbar}(u_2-v)$$ is divisible by $u_1-u_2-\hbar$. Moreover the quotient $\displaystyle \frac{\mathcal{F}}{u_1-u_2-\hbar}$ is symmetric in $u_1,u_2$. This allows us to use (2) of Lemma \ref{le:shift-gln} to conclude that $\mathcal{S}^q=0$.


\subsection{The variety $\F$}

Fix integers $1\leq n\leq d\in\N$, let
\[\F=\left\{ 0=V_0 \subseteq V_1 \subseteq\cdots\subseteq V_n = \C^d \right\}\]
be the variety of $n$--step flags in $\C^d$, and $T^*\F$ its cotangent bundle.
We describe below the $GL_d(\C)\times\nC$--equivariant $K$--theory and
cohomology of $T^*\F$ following \cite{ginzburg-vasserot,vasserot-qaffine}.

The connected components of $\F$ are parametrised by the set $\PP$ of
partitions of $[1,d]$ into $n$ intervals, \ie
\[\PP=\{ \dd = (0=d_0\leq d_1\leq \cdots \leq d_n=d) \}\]
where $\dd\in\PP$ labels the component $\F_\dd$ consisting of flags such
that $\dim V_k=d_k$. The symmetric group $\Sym_d$ acts on the rings
\begin{gather*}
S = \C[q^{\pm 1}, X_1^{\pm 1}, \ldots, X_d^{\pm 1}]\\
R = \C[\hbar, x_1,\ldots, x_d]
\end{gather*}
by permuting the variables $X_1,\ldots,X_d$ and $x_1,\ldots,x_d$ and fixing
$q,\hbar$ respectively. For any $\dd\in\PP$, denote by
\[\Sdd=\Sym_{d_1-d_0}\times\cdots\times\Sym_{d_n-d_{n-1}}\subset\Sym_d\]
the subgroup preserving the corresponding partition. Then, the following holds
\begin{gather*}
K^{GL_d(\C)\times \nC}(T^*\F)\cong\bigoplus_{\dd\in\PP} S^\Sdd\\
H_{GL_d(\C) \times \nC}(T^*\F)\cong\bigoplus_{\dd\in\PP} R^\Sdd
\end{gather*}
where $K^{\C^\times}(pt)=\C[q,q^{-1}]$ and $H_{\C^\times}(pt)=\C[\hbar]$.

\subsection{}
For any partition $\dd\in\PP$ and $i\in\bfI$, set
\[\dd_i^\pm = (0=d_0\leq\cdots\leq d_{i-1}\leq d_i\pm 1\leq d_{i+1}\leq\cdots\leq d_n=d)\]
if the \rhs makes sense as a partition.

If $\dd,\dd'\in\PP$ are two partitions, and $P$ is one of the rings $R,S$,
we denote by $\sigma(\dd,\dd')$ the symmetrisation operator
\[\sigma(\dd,\dd'):P^{\Sdd\cap\Sddp}\to P^\Sddp,\qquad
\sigma(\dd,\dd')(p)=\sum_{\tau\in\Sddp/\Sdd\cap\Sddp}\tau(p)\]

\subsection{$\hloopgl{n}$--action \cite{ginzburg-vasserot,vasserot-qaffine}}
\label{geometric-action-qaffine}

Consider the following operators acting on
\[S(\PP)=\bigoplus_{\dd\in\PP}S^\Sdd\]
\begin{enumerate}
\item For any $j\in\bfJ$, $\Psi_U(\Theta^\pm_j(z))$ acts on $S^\Sdd$ as
multiplication by
\[\prod_{k=1}^{d_{j-1}}\frac{qz - q^{-1}X_k}{z-X_k}
\prod_{k=d_j+1}^{d}\frac{z-X_k}{q^{-1}z - qX_k}
\in S^\Sdd[[z^{\mp 1}]]\]

\item For any $i\in\bfI$, the operators 
\begin{gather*}
\Psi_U(E_{i}(z)) : S^\Sdd \rightarrow S^{\Sym(\dd_i^+)}[[z,z^{-1}]]\\
\Psi_U(F_{i}(z)) : S^\Sdd \rightarrow S^{\Sym(\dd_i^-)}[[z,z^{-1}]]
\end{gather*}
act by $0$ if $\dd_i^\pm$ is not defined, and by
\begin{gather*}
\Psi_U(E_i(z))p=\sigma(\dd,\dd_i^+)\left(
\delta(X_{d_i+1}/z)\prod_{k\in I_i} \frac{qX_{d_i+1} - q^{-1}X_k}{X_{d_i+1} - X_k}p
\right)\\
\Psi_U(F_i(z))p=\sigma(\dd,\dd_i^-)\left(
\delta(X_{d_i}/z)\prod_{k\in I_{i+1}} \frac{q^{-1}X_{d_i} - qX_k}{X_{d_i} - X_k}p
\right)
\end{gather*}
otherwise, where $I_i$ is the interval $[d_{i-1}+1,\ldots,d_i]$.
\end{enumerate}

\noindent The following result is due to Ginzburg and Vasserot and
is proved in \cite[\S 2.2]{vasserot-qaffine}.
\begin{thm}\label{th:GV}
The assignment $\Psi_U$ extends to an algebra homomorphism
\[\Psi_U : \hloopgl{n} \rightarrow\End_{\C[q,q^{-1}]}(S(\PP))\]
\end{thm}

\subsection{$\Yhgl{n}$--action}\label{geometric-action-yangian}

Consider the following operators acting on
\[R(\PP) = \bigoplus_{\dd\in\PP} R^\Sdd\]

\begin{enumerate}
\item For any $j\in\bfJ$, $\Psi_Y(\theta_j(u))$ acts on $R^\Sdd$ as
multiplication by
\[\prod_{k=1}^{d_{j-1}}\frac{u-x_k+\hbar}{u-x_k}\prod_{k=d_j+1}^{d}
\frac{u-x_k}{u-x_k-\hbar}\in R^\Sdd[[u^{-1}]]\]

\item For any $i\in \bfI$, 
\begin{gather*}
\Psi_Y(e_i(u)) : R^\Sdd \rightarrow R^{\Sym(\dd_i^+)}[[u^{-1}]]\\
\Psi_Y(f_i(u)) : R^\Sdd \rightarrow R^{\Sym(\dd_i^-)}[[u^{-1}]]
\end{gather*}
act as zero if $\dd_i^\pm$ is not defined, and by
\begin{align*}
\Psi_Y(e_i(u))p
&=
\hbar\sigma(\dd,\dd_i^+)
\left(\frac{1}{u-x_{d_i+1}}\prod_{k\in I_i} \frac{x_{d_i+1} - x_k+\hbar}{x_{d_i+1} - x_k}p\right)\\
\Psi_Y(f_i(u))p
&=
\hbar\sigma(\dd,\dd_i^-)
\left(\frac{1}{u-x_{d_i}}\prod_{k\in I_{i+1}} \frac{x_{d_i} - x_k-\hbar}{x_{d_i} - x_k}p\right)
\end{align*}
otherwise.
\end{enumerate}

\noindent
The following result is proved in a similar way to Theorem \ref{th:GV}
\begin{prop}
The assignment $\Psi_Y$ extends to an algebra homomorphism
\[\Psi_Y : \Yhgl{n} \rightarrow\End_{\C[\hbar]}(R(\PP))\]
\end{prop}

\begin{rem}
The above formulae are degenerations of those of the previous section obtained
by setting $z=e^{tu},q=e^{t\hbar/2},X_k=e^{tx_k}$ and letting $t\to 0$.
\end{rem}

\subsection{}

\begin{lem}\label{useful-1-gln}
The homomorphism $\Psi_Y$ maps the centre $\calZ$ of $\Yhgl{n}$
surjectively to $\C[\hbar, x_1,\ldots, x_d]^{\Sym_d}$. In particular, there
exists an element
\[\Delta(u) = 1 + \hbar\sum_{r\geq 0} \Delta_r u^{-r-1}\in\calZ[[u^{-1}]]\]
such that
\[\Psi_Y(\Delta(u))=\prod_{k=1}^d \frac{u-x_k-\hbar}{u-x_k}\]
\end{lem}
\begin{pf}
By \cite[Cor. 1.11.8]{molev-yangian}, $\calZ$ is generated by the coefficients
of the element
\[qdet(u) = \theta_1(u)\theta_2(u-\hbar)\cdots \theta_n(u-(n-1)\hbar)\in\Yhgl{n}[[u^{-1}]]\]
It readily follows from \ref{geometric-action-yangian} that
\[\Psi_Y(qdet(u))=\prod_{k=1}^d \frac{u-x_k}{u-x_k - (n-1)\hbar}\]
By \eqref{eq:B log}, $L(v)=B(\log(qdet(u)))\in\calZ[[v]]$ therefore satisfies
\[\Psi_Y(L(v))=
\sum_{k=1}^d\frac{e^{(x_k+(n-1)\hbar)v}-e^{x_k v}}{v}=
\sum_{r\geq 1}\Bigl(p_r(\{x_k+(n-1)\hbar\})-p_r(\{x_k\})\Bigr)\frac{v^{r-1}}{r!}\]
which yields the surjectivity since the power sums $p_r(x_1,\ldots,x_d)=\sum_k
x_k^r$ generate $\C[x_1,\ldots, x_d]^{\Sym_d}$.
\end{pf}

\subsection{}

We will need the following 

\begin{lem}\label{existence-todd}
For any $i\in\bfI$, there exists $\Td^{\pm}_i(v) = \sum_{r\geq 0}\Td^{\pm}_{i,r}
\fact{v^r}{r!}\in\widehat{Y^0[v]}$ 
such that
\begin{align*}
\Psi_Y\left(\Td^+_i(v)\right)
&= \prod_{k\in I_i} \frac{v-x_k}{1-e^{-v+x_k}} \frac{1-e^{-v+x_k-\hbar}}{v-x_k+\hbar}\\
\Psi_Y\left(\Td^-_i(v)\right)
&= \prod_{k\in I_{i+1}} \frac{v-x_k}{1-e^{-v+x_k}} \frac{1-e^{-v+x_k+\hbar}}{v-x_k-\hbar}
\end{align*}
\end{lem}
The proof of this lemma is given in Section \ref{section-proof-todd}.

\subsection{A compatible assignment}\label{definition-assignment-gln}

Let $\Phi:\{E_{i,0}, F_{i,0},D_{j,r}\}_{i\in\bfI,j\in\bfJ,r\in\Z}\rightarrow\wh{\Yhgl{n}}$
be the assignment defined by
\begin{align*}
\Phi(D_{j,0}) &= \theta_{j,0}\\
\Phi(D_{j,r}) &= \left.\frac{B_j(v)}{q-q^{-1}}\right|_{v=r}\\
\Phi(E_{i,0}) &= \sum_{s\geq 0} e_{i,s} \fact{\Td^+_{i,s}}{s!}\, q^{-\Delta_0 - \theta_{i,0}}\\
\Phi(F_{i,0}) &= \sum_{s\geq 0} f_{i,s} \fact{\Td^-_{i,s}}{s!}\, q^{\Delta_0 + \theta_{i+1,0}}
\end{align*}
where $\Delta_0$ is given by Lemma \ref{useful-1-gln}. Extend $\Phi$ to the generators
$E_{i,r},F_{i,r}$, $r\in\Z$ by defining, as in Section \ref{ss:Phi QL23}, 
\[\Td^{\pm,(r)}_i(v)=
\sum_{m\geq 0}\Td^{\pm,(r)}_{i,m}v^m=
e^{rv}\Td^\pm_i(v)\in\wh{Y^0}[[v]]\]
and setting
\begin{align*}
\Phi(E_{i,r}) &= \sum_{s\geq 0} e_{i,s} \fact{\Td^{+,(r)}_{i,s}}{s!}\,q^{-\Delta_0 - \theta_{i,0}}\\
\Phi(F_{i,r}) &= \sum_{s\geq 0} f_{i,s} \fact{\Td^{-,(r)}_{i,s}}{s!}\,q^{\Delta_0 + \theta_{i+1,0}}
\end{align*}

\subsection{}\label{ss:pf main-2-gln}

Let $\wh{R(\PP)}$ be the completion with respect to the $\N$--grading given by
$\deg(x_k)=\deg(\hbar)=1$. Define a homomorphism $\eexp:S(\PP)\rightarrow
\wh{R(\PP)}$ mapping each $S^\Sdd$ to $\wh{R^\Sdd}$ by
\[q\longmapsto e^{\hbar/2}
\aand
X_k \longmapsto e^{x_k}\]

\begin{thm}\label{main-2-gln}
The assignment $\Phi$ above intertwines the geometric realisations of $\hloopgl{n}$
and $\Yhgl{n}$ on $S(\PP)$ and $R(\PP)$ respectively. In other words, the following
holds for any $X\in
\{E_{i,r}, F_{i,r},D_{j,r}\}_{i\in\bfI,j\in\bfJ,r\in\N}$ and $\pi\in S(\PP)$.
\[\eexp(X\cdot \pi) = \Phi(X)\cdot \eexp(\pi)\]
\end{thm}
\begin{pf}
Consider first the case $X=D_{j,r}$, $j\in\bfJ,r\in\Z$. By definition of $\Psi_U$ and $\Psi_
Y$, $D_{j,0}=2/\hbar\log(\Theta_{j,0})$ and $\theta_{j,0}$ act on $S^\Sdd$ and $R^\Sdd$
respectively as multiplication by $d-(d_j-d_{j-1})$. Further, \eqref{eq:B log} yields
\[\Psi_Y(B_j(v))=\frac{1}{v}\left(
(1-e^{-\hbar v})\sum_{k=1}^{d_{j-1}}e^{x_kv} +
(e^{\hbar v}-1)\sum_{k=d_j+1}^d e^{x_kv}\right)\]
Similarly, taking $\log$ in
\[\Psi_U\left(\exp\left((q-q^{-1})\sum_{s\geq 1}D_{j,s}z^{-s}\right)\right)
=
\prod_{k=1}^{d_{j-1}}\frac{z - q^{-2}X_k}{z-X_k}
\prod_{k=d_j+1}^d\frac{z-X_k}{z - q^2X_k}\]
yields
\[\Psi_U((q-q^{-1})D_{j,r}) =
\frac{1}{r}\left((1-q^{-2r}) \sum_{k=1}^{d_{j-1}}X_j^r +
(q^{2r}-1)\sum_{k=d_j+1}^d X_k^r\right)\]
Thus, $\eexp(\Psi_U(D_{j,r})\pi) = \Psi_Y\left(B_j(r)\pi/(q-q^{-1})\right)$ for any
$\pi\in S^\Sdd$.

We turn next to $X=E_{i,0}$. Let $\pi\in S^\Sdd$ and set $p = \eexp(\pi)\in R^
\Sdd$. Since $\Delta_0$ acts on $R^\Sdd$ as multiplication by $-d$ by Lemma
\ref{useful-1-gln}, and $\theta_{j,0}$ acts as multiplication by $d-(d_j-d_{j-1})$,
we get
\[\begin{split}
\Phi(E_{i,0})(p)
&= \sum_{s\geq 0} e_{i,s}\left(\fact{\Td_{i,s}^+}{s!}\,p\right) q^{d_j-d_{j-1}} \\ 
&=
\sigma(\dd,\dd_i^+)\left(\sum_{s\geq 0} \Td^+_{i,s} \fact{x_{d_i+1}^s}{s!}p
\prod_{k\in I_i} \frac{x_{d_i+1} - x_k + \hbar}{x_{d_i+1} - x_k}\right)q^{d_j-d_{j-1}}\\
&=
\sigma(\dd,\dd_i^+)\left(\Td^+_i(x_{d_i+1})\,p
\prod_{k\in I_i} \frac{x_{d_i+1} - x_k + \hbar}{x_{d_i+1} - x_k}\right)q^{d_j-d_{j-1}}\\
& =
\sigma(\dd,\dd_i^+)
\left(p\prod_{k\in I_i}
\frac{e^{x_{d_i+1}} - e^{x_k-\hbar}}{e^{x_{d_i+1}} - e^{x_k}}\right)q^{d_j-d_{j-1}}\\
&=
\sigma(\dd,\dd_i^+)\left(p\eexp\prod_{k\in I_i} \frac{qX_{d_i+1} - q^{-1}X_k}{X_{d_i+1} - X_k}\right)\\
&=
\eexp\left(E_{i,0}\,\pi\right)
\end{split}\]
The proof for the rest of the generators is identical.
\end{pf}

\subsection{Proof of Lemma \ref{existence-todd}}\label{section-proof-todd}

Let $\Delta(u)$ be the formal series given in Lemma \ref{useful-1-gln},
and set
\[\zed(u) = \hbar\sum_{r\geq 0} \zed_ru^{-r-1} = \log(\Delta(u))\]
For any $j\in\bfJ$, define $y_j(u)\in Y^0[[u^{-1}]]$ by
\begin{equation}\label{t-series-gln}
y_j(u)=
\zed(u+(j-1)\hbar) + d_j(u) + \sum_{s=1}^{j-1}
\left(d_{j-s}(u+s\hbar) - d_{j-s}(u+(s-1)\hbar)\right)
\end{equation}
A computation similar to the one given in \ref{ss:pf main-2-gln} shows that,
for any $j\in\bfJ$,
\begin{equation}\label{eq:useful-2-gln}
\Psi_Y\left(B(y_j(u))\right) = \frac{1-e^{\hbar v}}{v} \sum_{k\in I_j} e^{x_kv}
\end{equation}

Set now\footnote{Note the difference between $J(v)$ and the function $G(v)$
used in Section \ref{ss:G(v)} for constructing the solutions for simple Lie algebras:
$J(v) = G(v) + \frac{v}{2}$.}
\[J(v) = \log\left(\frac{v}{1-e^{-v}}\right)\in\Q[[v]]\]
and, for any $i\in\bfI$, define
\begin{align}
\td_i^+(v) &= B(y_i(-\partial))J'(v+\hbar)
\label{small-todd-+}\\
\td_i^-(v) &= -B(y_{i+1}(-\partial))J'(v)
\label{small-todd--}
\end{align}
where $\partial = d/dv$. We claim that $\Td^{\pm}_i(u) = \exp(td^{\pm}_i(u))$ satisfy
the conditions of the Lemma. By \eqref{eq:useful-2-gln} we have
\[\Psi_Y\left(td^+_i(v)\right)=
\left(\frac{1-e^{-\hbar\partial}}{-\partial} \sum_{k\in I_i} e^{-x_k\partial}\right) \partial J(v+\hbar)\]
Using $e^{-p\partial}f(v) = f(v-p)$, we get
\[
\begin{split}
\Psi_Y(td^+_i(v)) &= \sum_{k\in I_i}J(v-x_k) - J(v-x_k+\hbar)\\
 &= \sum_{k\in I_i} \log\left(\frac{v-x_k}{1-e^{-v+x_k}} \frac{1-e^{-v+x_k-\hbar}}{v-x_k+\hbar}\right)\\
 &= \log\left(\prod_{k\in I_i} \frac{v-x_k}{1-e^{-v+x_k}} \frac{1-e^{-v+x_k+\hbar}}{v-x_k+\hbar}\right)
\end{split}
\]
The proof for the $-$ case is same.

\subsection{Standard form of $\Phi$}\label{algebraic-gln}

We rewrite below the assignment $\Phi$ in a form in which Theorem \ref
{first-main-theorem-gln} can be applied, and use this to prove that $\Phi$
extends to an algebra homomorphism $\hloopgl{n}\to\wh{\Yhgl{n}}$.

\begin{lem}\label{useful-3-gln}
Let $y_j(v)$ be given by \eqref{t-series-gln}, and $\lambda^{\pm}_i(u)$ the
operators of Proposition \ref{raising-lowering-operators-gln}. Then,
\[(\lambda^{\pm}_i(u)-1)B(y_j(v)) = 
\pm(\delta_{i,j} - \delta_{j,i+1})
\frac{e^{\hbar v}-1}{v} e^{uv} \]
\end{lem}
The proof of this lemma essentially follows from Proposition \ref{raising-lowering-operators-gln}.

\begin{cor}\label{raising-lowering-todd}
For any $i,i'\in\bfI$, we have
\begin{gather*}
\frac{\lambda^+_i(u)(\Td_{i'}^+(v))}{\Td_{i'}^+(u)}
=
\frac{\Td_{i'}^+(v)}{\lambda^-_i(u)(\Td_{i'}^+(v))}
=
\left(\frac{v-u+\hbar}{1-e^{-v+u-\hbar}}
\frac{1-e^{-v+u}}{v-u}\right)^{\delta_{i,i'} - \delta_{i,i'-1}}\\
\frac{\lambda^+_i(u)(\Td_{i'}^-(v))}{\Td_{i'}^-(u)}
=
\frac{\Td_{i'}^-(v)}{\lambda^-_i(u)(\Td_{i'}^-(v))}
=\left(\frac{v-u}{1-e^{-v+u}}\frac{1-e^{-v+u+\hbar}}{v-u-\hbar}\right)
^{\delta_{i,i'} - \delta_{i,i'+1}}
\end{gather*}
\end{cor}

\noindent
It follows that
\begin{gather*}
\Phi(E_{i,k})=e^{k\sigma_i^+}g_i^+(\sigma_i^+)e_{i,0}\\
\Phi(F_{i,k})=e^{k\sigma_i^-}g_i^-(\sigma_i^-)f_{i,0}
\end{gather*}
where
\begin{equation}\label{eq: g+- gln}
\begin{split}
g_i^+(u) &= q^{-\Delta_0 - \theta_{i,0}} \frac{\hbar}{q-q^{-1}} \Td_i^+(u)\\
g_i^-(u) &= q^{\Delta_0 + \theta_{i+1,0}} \frac{\hbar}{q-q^{-1}}\Td_i^-(u)
\end{split}
\end{equation}

\subsection{}

We record the action of the operators $\lambda^{\pm}_{i'}(u)$ on $g_i^{\pm}(v)$ 
using Corollary \ref{raising-lowering-todd}
\begin{gather}
\lambda^+_i(u)(g_i^+(v))
= \lambda^-_{i-1}(u)(g_i^+(v)) = g_i^+(v) \frac{v-u+\hbar}{e^{v+\hbar/2} - e^{u-\hbar/2}} \frac{e^v - e^u}{v-u}
\label{g+-gln}\\
\lambda^+_{i+1}(u)(g_i^-(v))
= \lambda^-_i(u)(g_i^-(v)) = g_i^-(v)\frac{v-u-\hbar}{e^{v-\hbar/2} - e^{u+\hbar/2}} \frac{e^v - e^u}{v-u}
\label{g--gln}
\end{gather}
Using the fact that $\lambda^+_i(u)\lambda^-_i(u) = \Id$, we get four more equations
from these. Moreover, $\lambda^{\pm}_{i'}(u)(g_i^+(v))=g_i^+(v)$ for $i'\neq i,i-1$ and
$\lambda^{\pm}_{i'}(u)(g_i^-(v))=g_i^-(v)$ for $i'\neq i,i+1$.
 
\subsection{}

\begin{thm}
The series $g_i^\pm(u)$ satisfy the conditions (A),(B),(C0)--(C2) of Theorem
\ref{first-main-theorem-gln} and therefore give rise to an algebra homomorphism
$\Phi:\hloopgl{n}\to\wh{\Yhgl{n}}$.
\end{thm}

\subsection{Proof of  (A)}

We need to prove that for every $i, i'\in \bfI$, we have
\[g_i^+(u)\lambda^+_i(u)(g_{i'}^-(v)) = g_{i'}^-(v)\lambda_{i'}^{-}(v)(g_i^+(u))\]
If $i\neq i',i'+1$, both sides are equal to $g_i^+(u)g_{i'}^-(v)$. For $i=i'$, by
\eqref{g+-gln}--\eqref{g--gln}, the left and right--hand sides are respectively
equal to
\[\frac{v-u}{e^v - e^u} \frac{e^{v-\hbar/2} - e^{u+\hbar/2}}{v-u-\hbar}
\aand
\frac{u-v}{e^u - e^v} \frac{e^{u+\hbar/2} - e^{v-\hbar/2}}{u-v+\hbar}\]
The case $i=i'+1$ follows in the same way.

\subsection{Proof of (B)}

Let $i\in \bfI$. By \eqref{g--gln},
\[\begin{split}
g_i^+(u)\lambda^+_i(u)(g_i^-(u))
&= g_i^+(u)g_i^-(u) \frac{q-q^{-1}}{\hbar} \\
&= q^{\theta_{i+1,0} - \theta_{i,0}} \frac{\hbar}{q-q^{-1}} \Td_i^+(u)\Td_i^-(u)\\
&= \frac{\hbar}{q-q^{-1}} \exp \left(\frac{\hbar(\theta_{i+1,0}-\theta_{i,0})}{2} + td_i^+(u) + td_i^-(u)\right)
\end{split}\]
By definition of $td_i^{\pm}$,
\[\begin{split}
td_i^+(u) + td_i^-(u) &= -By_{i+1}(-\partial) J'(u) + By_i(-\partial) J'(u+\hbar)\\
 &= \left.-B(y_{i+1}(v) - y_i(v+\hbar))\right|_{v=-\partial} J'(u)
\end{split}\]
where the second equality follows from $J'(u+\hbar)=e^{\hbar\partial}J'(u)$ and the fact
that $e^{pv}B(f(u))=B(f(u+p))$. Next, the definition of $y_i$ yields
\[y_{i+1}(v) - y_i(v+\hbar) = d_{i+1}(v) - d_i(v)\]
hence
\[td^+_i(u) + td^-_i(u) = \hbar\sum_{r\geq 0} (-1)^{r+1} \frac{d_{i+1,r} - d_{i,r}}{r!} J^{(r+1)}(u)\]
which implies that
\[\frac{\hbar(\theta_{i+1,0}-\theta_{i,0})}{2} + td^+_i(u) + td^-_i(u) =  \hbar\sum_{r\geq 0} (-1)^{r+1} \frac{d_{i+1,r} - d_{i,r}}{r!} G^{(r+1)}(u)\]
and the proof of (B) follows from Proposition \ref{pr:condition-B}

\subsection{Proof of (C0)--(C2)}

The condition (C0) follows from the fact that $(\lambda_i^\pm(u)-1)(g_{i'}^\pm(v))
=0$ if $|i-i'|>1$. Since the proof of (C1) is the same as the one given in the verification
of (A), we are left with checking (C2). We need to show that, for any $i\in\bfI\setminus
\{n-1\}$, 
\[g_i^+(u)\lambda_i^+(u)g_{i+1}^+(v)\frac{e^u - e^v}{u-v} =
g_{i+1}^+(v)\lambda^{\pm}_{i+1}(v)g_i^+(u)\frac{e^{u-\hbar/2} - e^{v+\hbar/2}}{u-v-\hbar}\]

Using \eqref{g+-gln}--\eqref{g--gln}, the left and right--hand sides are respectively
equal to
\[\frac{u-v}{e^u - e^v} \frac{e^{u-\hbar/2} - e^{v+\hbar/2}}{u-v-\hbar} \frac{e^u - e^v}{u-v}
\aand
\frac{e^{u-\hbar/2} - e^{v+\hbar/2}}{u-v-\hbar}\]

\subsection{Isomorphism between completions}\label{iso-completed-gln}

We give below an analogue of Theorem \ref{th:geo iso} for $\g=\gl_n$.
We begin by defining the appropriate completion of $\hloopgl{n}$, which
differs from the one used in \ref{ss:completion} due to the fact that $\gl
_n$ is not semisimple.

For each $r\geq 0$, $t\in \Z$ and $X=E_i,F_i,\Th_j$, where $i\in\bfI$
and $j\in \bfJ$, consider the element
\[X_{r;t}=\sum_{s=0}^r (-1)^s \cbin{r}{s}{} X_{s+t}\]
where $\Theta_{j,l}=(\Theta_{j,l}^+-\Theta_{j,l}^{-})/(q-q^{-1})$. Note
that $X_{r;t}=x\otimes z^t (1-z)^r\mod\hbar$ where $x\in\gl_n$ is such
that $X=x\mod\hbar$. Let $\K_r$ be the two--sided ideal of $\hloopgl{n}$
generated by the elements $\{X_{r';t}\}_{r'\geq r, t\in\Z}$, and $\hbar$
if $r=1$. Finally, let $\J_n\subset\hloopgl{n}$ be the ideal
\[\J_m=
\sum_{\substack{m_1,\ldots,m_k\geq 1\\m_1+\cdots+m_k=m}}\K_{m_1}\cdots\K_{m_k}\]
Then, $\J_m$ is a descending filtration, $\J_m\J_{m'}\subset\J_{m+m'}$,
and the completion
\[\wh{\hloopgl{n}}=\lim_{\longleftarrow}\hloopgl{n}/\J_m\]
is a flat deformation of the completion of $U(\gl_n[z,z^{-1}])$ at $z=1$.

\begin{rem}
The reason we use $\J_m$ instead of the powers of the evaluation
ideal $\J_1$ as in \ref{ss:completion} can be seen at the classical
level. Indeed, the correct filtration of $U(\gl_n[z,z^{-1}])$ is given by
the $J_m=U((z-1)^m\gl_n[z,z^{-1}])$ which contains, but does not
equal $J_1^m$. For example, if $I$ is the identity matrix, then $I
\otimes(z-1)^2\in J_2\setminus\bigcup_{m\geq 1}J_1^m$.
\end{rem}

For each $m\in \N$, set
\[\wh{\Yhgl{n}}_{\geq m}=\prod_{m'\geq m} (\Yhgl{n})_{m'}
\subset\wh{\Yhgl{n}}\]
where $(\Yhgl{n})_m$ is the subspace of degree $m$. The proof
of the following result is similar to that of Theorem \ref{th:geo iso}
and therefore omitted.

\begin{thm}\label{thm: iso-completed-gln}
The homomorphism $\Phi$ maps $\J_m$ to $\wh{\Yhgl{n}}_m$ for 
any $m\in\N$, and  induces an isomorphism of completed algebras 
\[\wh{\Phi} : \wh{\hloopgl{n}} \isom \wh{\Yhgl{n}}\]
\end{thm}
\appendix

\section{Proof of the Serre relations}\label{app: Serre}

\subsection{}

Let $\g$ be a complex, semisimple Lie algebra. The aim of this appendix is to prove
the following

\begin{prop}\label{pr:Serre}
Let $\Phi$ be the assignment $\{E_{i,k},F_{i,k},H_{i,k}\}\to\wh{\Yhg}$ given in Sections
\ref{definition-assignment}--\ref{ss:Phi QL23}, and assume that the relations (A) and
(B) of Theorem \ref{first-main-theorem} hold. Then, $\Phi$ is compatible with the $q
$--Serre relations (QL6).
\end{prop}

For $i\neq j\in\bfI$, set $m=1-a_{ij}$. Define, for any $\ul{k}=(k_1,\ldots,k_m)\in\Z^m$
and $l\in\Z$
\begin{multline}\label{app-eq: serre-arbitrary}
\mathcal{S}_{ij}^q(\underline{k},l)
=
\sum_{\pi\in \Sym_m} \sum_{s=0}^m (-1)^s \left[\begin{array}{c} m \\ s\end{array}\right]_{q_i}
\Phi(E_{i,k_{\pi(1)}}) \cdots \Phi(E_{i,k_{\pi(m-s)}}) \Phi(E_{j,l}) \\
\Phi(E_{i,k_{\pi(m-s+1)}}) \cdots \Phi(E_{i,k_{\pi(m)}})\in\wh{\Yhg}
\end{multline}
and let $\mathcal{S}_{ij}^q = \mathcal{S}_{ij}^q(\underline{0},0)$, explicitly given as follows
\begin{equation}\label{app-eq: serre-zero}
\mathcal{S}_{ij}^q=
\sum_{s=0}^m (-1)^s \left[\begin{array}{c} m\\ s\end{array}\right]_{q_i}
\left(\Phi(E_{i,0})\right)^{m-s} \Phi(E_{j,0}) \left(\Phi(E_{i,0})\right)^s
\end{equation}
Our aim is to show that $\mathcal{S}_{ij}^q(\underline{k},l)=0$. Let us outline the
main steps of the proof.
\begin{enumerate}
\item We first reduce the proof of $\Cs_{ij}^q(\ul{k},l)=0$ to $\Cs_{ij}^q=0$. This is
achieved in Lemma \ref{app-lem: reduction}.
\item By a standard argument using the representation theory of $U_\hbar\Lsl_2$,
we deduce in Lemma \ref{app-lem: serre-on-fd} that $\Cs_{ij}^q$ acts by zero on
any \fd representation of $\Yhg$.
\item Finally, we show that these representations separate points in $\Yhg$, and
hence that $\Cs_{ij}^q=0$. \S \ref{app-section: co-poisson} and \ref{app-section: ideals}
are devoted to the proof of this fact (Corollary \ref{app-corr: biideals}) which was
communicated to us by V. G. Drinfeld.
\end{enumerate}

\subsection{The algebra $\ol{Y}$}\label{app-sec: auxiliary}

Define an auxiliary algebra $\ol{Y}$ to be the unital, associative $\C[\hbar]$--algebra
generated by $\{\ol{\xi}_{i,r},\ol{x}_{i,r}\}_{i\in\bfI, r\in \N}$ subject to the following
relations
\begin{enumerate}
\item For every $i,j\in\bfI$ and $r,s\in \N$
\[[\ol{\xi}_{i,r}, \ol{\xi}_{j,s}]=0\]
\item For every $i,j\in\bfI$ and $s\in \N$
\[[\ol{\xi}_{i,0}, \ol{x}_{j,s}] = d_ia_{ij}\ol{x}_{j,s}\]
\item For every $i,j\in\bfI$ and $r,s\in\bfI$
\[[\ol{\xi}_{i,r+1}, \ol{x}_{j,s}] - [\ol{\xi}_{i,r}, \ol{x}_{j,s+1}]=
\frac{d_ia_{ij}\hbar}{2}(\ol{\xi}_{i,r}\ol{x}_{j,s} + \ol{x}_{j,s}\ol{\xi}_{i,r})\]
\end{enumerate}

We denote by $\ol{Y}^0\subset\ol{Y}$ the commutative subalgebra generated by
$\{\ol{\xi}_{i,r}\}_ {i\in\bfI, r\in \N}$ and by $\ol{Y}^{>0}$ the subalgebra of $\ol{Y}$
generated by $\{\ol{x}_{i,r}\}_{i\in\bfI, r\in \N}$. The latter is a free $\C[\hbar]$--algebra
over this set of generators. Moreover, by Corollary \ref{co:PBW}, $\ol{Y}\cong\ol{Y}
^0\otimes\ol{Y}^{>0}$.

\subsection{The operators $\ol{\sigma}_{i,(k)}$ and $\ol{\sigma}_j$}\label{app-sec: shift}

The algebra $\ol{Y}$ has a grading by the root lattice $Q$ given by
\[\ddeg(\ol{\xi}_{i,r})=0\aand\ddeg(\ol{x}_{i,r})=\alpha_i\]
Fix henceforth $i\neq j\in\bfI$, set $m=1-a_{ij}$ and let $\ol{Y}_{m\alpha_i+\alpha_j}$
be the homogeneous component of $\ol{Y}$ of degree $m\alpha_i+\alpha_j$.

Define operators $\ol{\sigma}_j,\ol{\sigma}_{i,(k)}$ on $\ol{Y}_{m\alpha_i+\alpha_j}$ as
follows. Since $\ol{Y}_{m\alpha_i+\alpha_j}\cong \ol{Y}^0\otimes \ol{Y}^{>0}_{m\alpha_i
+\alpha_j}$ and $\ol{Y}^{>0}$ is free, we have
\[\ol{Y}^{>0}_{m\alpha_i+\alpha_j}\cong
\ol{Y}^0\otimes\bigoplus_{s=0}^m
\ol{Y}(i)^{\otimes m-s}\otimes \ol{Y}(j) \otimes \ol{Y}(i)^{\otimes s}\]
where, for $a=i,j$, $\ol{Y}(a)=\ol{Y}^{>0}_{\alpha_a}$ is spanned by $\{\ol{x}_{a,r}\}_{r\in\N}$.
Let $\ol{\sigma}_a$ denote the $\C[\hbar]$--linear map on $\ol{Y}(a)$ given by
$\ol{\sigma}_a(x_{a,r}) = x_{a,r+1}$. For any $k = 1,\ldots, m$, define the $Y^0$--linear
operator $\ol{\sigma}_{i,(k)}$ on $\ol{Y}_{m\alpha_i+\alpha_j}$ by letting it act on the
summand $\ol{Y}(i)^{\otimes m-s}\otimes \ol{Y}(j)\otimes \ol{Y}^{\otimes s}(i)$ as
\[\ol{\sigma}_{i,(k)} 
=\left\{\begin{array}{ll}
1^{\otimes k-1}\otimes\ol{\sigma}_i\otimes 1^{\otimes m+1-k} & \mbox{if } k\leq m-s \\
1^{\otimes k} \otimes \ol{\sigma}_i\otimes 1^{\otimes m-k}& \mbox{otherwise} 
\end{array}\right.\]
Similarly, let $\sigma_j\in\End_{\ol{Y}^0}(\ol{Y}_{m\alpha_i+\alpha_j})$ be given by
$1^{\otimes m-s}\otimes \ol{\sigma}_j\otimes 1^{\otimes s}$ on $\ol{Y}(i)^{\otimes
m-s}\otimes\ol{Y}(j)\otimes \ol{Y}(i)^{\otimes s}$.

\subsection{The projection $\mathbf{p_{ij}}$}

Let $p:\ol{Y}\to\Yhg$ be the algebra homomorphism obtained by sending $\ol{\xi}_
{a,r}\longmapsto \xi_{a,r}$ and $\ol{x}_{a,r} \longmapsto x_{a,r}^+$ for every $a\in
\bfI$ and $r\in \N$, and let $p_{ij}$ be the restriction of $p$ to $\ol{Y}_{m\alpha_i+
\alpha_j}$. The following holds by Proposition \ref{algebraic-properties-shift}.

\begin{lem}\label{app-lem: symmetric}
The kernel of $p_{ij}$ is the $\C[\hbar]$--linear span of the following elements
\begin{enumerate}
\item For any $0\leq s\leq m-1$ and $A(u_1,\ldots,u_m,w)\in\ol{Y}^0[u_1,\ldots,u_m,w]$
\begin{multline*}
A(\ol{\sigma}_{i,(1)},\ldots, \ol{\sigma}_{i,(m)},\ol{\sigma}_j)
\Bigl((\ol{\sigma}_{i,(m-s)}-\ol{\sigma}_j-a\hbar)\ol{x}_{i,0}^{m-s}\ol{x}_{j,0}\ol{x}_{i,0}^s \\
- (\ol{\sigma}_{i,(m-s)}-\ol{\sigma}_j+a\hbar)\ol{x}_{i,0}^{m-s-1}\ol{x}_{j,0}\ol{x}_{i,0}^{s+1}\Bigr)
\end{multline*}
where $a=d_ia_{ij}/2$.
\item For any $0\leq s\leq m$, $k\in \{1,\ldots, m-1\}\setminus \{m-s\}$ and $A(u_1,\ldots,u_m,w)
\in\ol{Y}^0[u_1,\ldots, u_m,w]^{(k\, k+1)}$ 
\[A(\ol{\sigma}_{i,(1)},\ldots, \ol{\sigma}_{i,(m)},\ol{\sigma}_j)(\ol{\sigma}_{i,(k)}-\ol{\sigma}_{i,(k+1)} 
-d_i\hbar)\ol{x}_{i,0}^{m-s}\ol{x}_{j,0}\ol{x}_{i,0}^s\]
\item For every $A(u_1,\ldots, u_m,w)\in \ol{Y}^0[u_1,\ldots, u_m,w]^{\Sym_m}$
\[
A(\ol{\sigma}_{i,(1)},\ldots, \ol{\sigma}_{i,(m)},\ol{\sigma}_j)
\left(\sum_{s=0}^m (-1)^s
\left(\begin{array}{c} m \\ s\end{array}\right)
\ol{x}_{i,0}^{m-s}\ol{x}_{j,0}\ol{x}_{i,0}^s \right)
\]
\end{enumerate}
\end{lem}

\begin{cor}\label{app-corr: symmetric}
Let $X\in \Ker(p_{ij})$ and $A(u_1,\ldots, u_m,w) \in \ol{Y}^0[u_1,\ldots, u_m,w]^{\Sym_m}$.
Then,
\[A(\ol{\sigma}_{i,(1)},\ldots, \ol{\sigma}_{i,(m)}, \ol{\sigma}_j)X \in \Ker(p_{ij})\]
\end{cor}

\subsection{Reduction step}\label{app-sec: reduction}

Let $\ol{\Cs_{ij}^q}(\ul{k},l),\ol{\Cs_{ij}^q}$ denote the elements of $\ol{Y}_{m\alpha_i+\alpha_j}$
defined by the same expressions as \eqref{app-eq: serre-arbitrary}--\eqref{app-eq: serre-zero}.
Then,
\[\ol{\Cs_{ij}^q}(\ul{k},l) =
\left(\sum_{\pi\in \Sym_m}
e^{k_{\pi(1)}\ol{\sigma}_{i,(1)}}\cdots e^{k_{\pi(m)}\ol{\sigma}_{i,(m)}}e^{l\ol{\sigma}_j}\right)
\ol{\Cs_{ij}^q}\]
Using Corollary \ref{app-corr: symmetric}, we obtain the following
\begin{lem}\label{app-lem: reduction}
$\Cs_{ij}^q=0$ implies $\Cs_{ij}^q(\ul{k},l) = 0$ for every $k_1,\cdots, k_m,l\in \Z$.
\end{lem}

\subsection{}

By a \fd representation of $\wh{\Yhg}$ we shall mean a \fing topologically free
$\C[[\hbar]]$--module endowed with a $\C[[\hbar]]$--linear action of $\wh{\Yhg}$.

\begin{lem}\label{app-lem: serre-on-fd}
Let $\V$ be a \fd representation of $\wh{\Yhg}$. Then, $\Cs_{ij}^q$ acts by zero on $\V$.
\end{lem}
\begin{pf} Let $\U_i$ be the subalgebra of $\wh{\Yhg}$ generated by
\[\E_i=\Phi(E_{i,0})\qquad \F_i=\Phi(F_{i,0})\qquad \H_i=\Phi(H_{i,0}) \]
By Lemma \ref{le:AB}, $\{\E_i, \F_i, \H_i\}$ satisfy the defining relations of the quantum
group $U_{\hbar_i}\Lsl_2$, where $\hbar_i=d_i\hbar/2$. We use the following notation
of $q$--adjoint operator (see \cite [\S 4.18]{jantzen}) which gives a representation of $
\U_i$ on any algebra containing it
\begin{gather*}
\ad_q(\E_i)(X) = \E_iX - \K_iX\K_i^{-1}\E_i\\
\ad_q(\F_i)(X) = \F_iX\K_i - X\F_i\K_i\\
\ad_q(\H_i)(X) = [\H_i,X] 
\end{gather*}
where $\K_i=q_i^{\H_i}=e^{\hbar_i\H_i}$.
Let $\rho : \wh{\Yhg} \rightarrow \End(\V)$ be the representation. Then,
\begin{gather*}
\ad_q(\rho(\F_i))\rho(\E_j) = 0\\
\ad_q(\rho(\H_i))\rho(\E_j) = a_{ij}\rho(\E_j)
\end{gather*}
where the first identity follows from Lemma \ref{le:AB}. 
Thus, as a $\U_i$--module, $\End(\V)$ contains the lowest weight vector $\rho(\E_j)$
of weight $a_{ij}$. By the representation theory of $U_{\hbar_i}\Lsl_2$, we get
\[\ad_q(\rho(\E_i))^m\rho(\E_j) = \rho\left(\ad_q(\E_i)^m\E_j\right) = 0\]
and the assertion follows from the well--known identity (see \cite[Lemma 4.18]{jantzen})
\[\ad_q(\E_i)^m\E_j = \sum_{s=0}^m (-1)^s \left[\begin{array}{c} m \\ s\end{array}\right]_{q_i} \E_i^{m-s}\E_j\E_i^s\]
\end{pf}

\subsection{}

Let $I_\hbar\subset\Yhg$ be the ideal defined by
\[I_\hbar=\bigcap_{(\V,\rho)}\,\Ker(\rho)\]
where $\V$ runs over all \fd {\it graded} modules over $\Yhg$, that is \fing
torsion--free $\C[\hbar]$--modules admitting a $\C[\hbar]$--linear action
$\rho:\Yhg\to\End(\V)$ and a $\Z$--grading compatible with that on $\Yhg$.

\begin{lem}\label{app-lem: sufficiency}
$I_\hbar=0$ implies $\Cs_{ij}^q=0$.
\end{lem}
\begin{pf}
The action of $\Yhg$ on any \fd graded module $\V$ extends to one of $\wh
{\Yhg}$ on the completion $\wh{\V}$ of $\V$ \wrt its grading. By Lemma \ref
{app-lem: serre-on-fd}, $\Cs_{ij}^q$ acts by $0$ on $\wh{\V}$ and therefore
so do its homogeneous components $\Cs_{ij;n}^q\in\Yhg$, $n\geq 0$ on
$\V$. Thus, $\Cs_{ij;n}^q\in I_\hbar$ for any $n$ and $\Cs_{ij}^q=0$.
\end{pf}

\subsection{}\label{app-section: co-poisson}

The following result, and its proof are due to V. Drinfeld \cite{drinfeld-personal-communication}
\begin{prop}\label{pr:Drinfeld}
The ideal $I_\hbar\subset\Yhg$ is trivial.
\end{prop}
\begin{pf}
It suffices to show that $I=I_\hbar/\hbar I_\hbar$ is trivial. Indeed, if $I_\hbar=\hbar I_
\hbar$, then $I_\hbar=\bigcap_k\hbar^k I_\hbar\subset \bigcap_k\hbar^k\Yhg=0$. By
definition of $I_\hbar$, $I_\hbar\cap\hbar\Yhg=\hbar I_{\hbar}$
so that $I$ embeds into $
\Yhg/\hbar\Yhg=U(\g[s])$. Since graded representations are stable under tensor
product, $I_\hbar$ is a Hopf ideal of $\Yhg$, that is 
\[\Delta(I_\hbar) \subset\Yhg\otimes I_\hbar + I_\hbar\otimes\Yhg\]
It follows that $I$ is a co--Poisson Hopf ideal of $U(\g[s])$. By Corollary \ref{app-corr: biideals}
below, any such ideal is either trivial or equal to $U(\g[s])$. Since $\Yhg$ possesses
non--trivial \fd graded representations, for example the action on $\C[\hbar]$ given
by the counit, $I_\hbar$ is a proper ideal of $\Yhg$ and is therefore equal to zero.
\end{pf}

\subsection{}\label{app-section: ideals}

Recall that a co--Poisson Hopf algebra $A$ is a Hopf algebra together with a Poisson
cobracket $\delta:A\rightarrow A\wedge A$ satisfying the following compatibility condition
(see \cite[\S 6.2]{chari-pressley} for details):
\[\delta(xy) = \delta(x)\Delta(y) + \Delta(x)\delta(y)\]
For a Lie algebra $\ga$, there is a one to one correspondence between co--Poisson structures
on $U\ga$ and Lie bialgebra structures on $\ga$ \cite[Proposition 6.2.3]{chari-pressley}. Moreover,
there is a one to one correspondence between co--Poisson Hopf ideals of $U\ga$ and Lie bialgebra
ideals of $\ga$.

The Lie bialgebra structure on $\Lg[s]$ is given by
\begin{gather}
\delta:\Lg[s] \rightarrow \Lg[s]\otimes \Lg[s] \cong (\Lg\otimes \Lg)[s,t]\nonumber\\
\delta(f)(s,t)=\left(\ad(f(s))\otimes 1 + 1\otimes\ad(f(t))\right)\left(\frac{\Omega}{s-t}\right)
\label{app-eq: delta}
\end{gather}
where $\Omega\in \Lg\otimes \Lg$ is the Casimir tensor. Note that $\delta$ lowers the degree by $1$.

Let $\ga\subset\Lg[s]$ be the Lie bialgebra ideal corresponding to co--Poisson Hopf ideal $I\subset
\current$. By the discussion given in previous paragraph,
\begin{equation}\label{app-eq: biideal}
\delta(\ga) \subset \ga\otimes \Lg[s] + \Lg[s] \otimes \ga
\end{equation}

\begin{lem}\label{app-lem: ideals}
Let $\ga \subset \Lg[s]$ be an ideal. Then $\ga$ is of the form $\ga = \Lg\otimes g\C[s]$ for some
polynomial $g\in\C[s]$.
\end{lem}

\begin{pf}
Let $S \subset \C[s]$ be the set of all polynomials $f$ such that there exists some non zero $x\in \Lg$
for which $x\otimes f \in \ga$. We claim that $S$ is an ideal of $\C[s]$. Let $f\in S$ and $g\in \C[s]$.
Let $0\neq x\in \Lg$ be such that $x\otimes f\in \ga$, and choose $y\in \Lg$ such that $[x,y]\not= 0$.
Then \[[x,y]\otimes fg=[x\otimes f, y\otimes g]  \in\ga\]
and hence $fg\in S$.

Now for any $f\in S$, the set $\{x\in \Lg : x\otimes f \in \ga\}$ is an ideal in $\Lg$, which is non--zero
and hence equal to $\Lg$. This proves that $\ga = \Lg \otimes S$. Since $\C[s]$ is a principal ideal
domain, the lemma is proved.
\end{pf}

\begin{cor}\label{app-corr: biideals}
Let $\ga \subset \Lg[s]$ be a Lie bialgebra ideal. Then either $\ga = 0$ or $\ga = \Lg[s]$.
\end{cor}

\begin{pf}
Let $g\in \C[s]$ be such that $\ga = \Lg\otimes (g) \subset \Lg[s]$. By \eqref{app-eq: delta},
we know that the Lie cobracket $\delta$ lowers the degree by $1$. Using \eqref{app-eq: biideal},
we conclude that $g$ is a constant polynomial.
\end{pf}

\end{document}